\documentclass[11pt,a4paper]{amsart}
\usepackage[all]{xy}
\usepackage{amssymb}
\usepackage{amsmath}
\usepackage{bm}
\usepackage[colorlinks,linkcolor=blue]{hyperref}
\usepackage{color,xcolor,dsfont}
\usepackage{enumerate}
\usepackage{epsfig}
\usepackage{extarrows}
\usepackage{framed}
\usepackage{graphicx}
\usepackage{marvosym}
\usepackage{mathrsfs}
\usepackage{setspace} 
\usepackage{tikz}
\usepackage{tikz-cd}

\newcommand{\scha}{\mathcal{A}}

\newcommand{\schc}{\mathcal{C}}
\newcommand{\schd}{\mathcal{D}}
\newcommand{\sche}{\mathcal{E}}
\newcommand{\schf}{\mathcal{F}}
\newcommand{\schk}{\mathcal{K}}
\newcommand{\scht}{\mathcal{T}}
\newcommand{\schp}{\mathcal{P}}

\newcommand{\she}{\mathscr{E}}
\newcommand{\shf}{\mathscr{F}}
\newcommand{\shg}{\mathscr{G}}
\newcommand{\shh}{\mathscr{H}}
\newcommand{\shi}{\mathscr{I}}

\newcommand{\sho}{\mathscr{O}}

\newcommand{\shr}{\mathscr{R}}

\newcommand{\shu}{\mathscr{U}}
\newcommand{\shv}{\mathscr{V}}

\newcommand{\stable}{\textsf{\textup{st}}}
\DeclareMathOperator{\Stab}{Stab}

\newcommand{\deriver}{\textup{\textsf{R}}}

\newcommand{\Hoch}{\textup{\textsf{HH}}}
\newcommand{\Kuc}{\cate{K\!u}}

\DeclareMathOperator{\id}{id}
\DeclareMathOperator{\GL}{GL}

\DeclareMathOperator{\NS}{NS}
\DeclareMathOperator{\Aut}{Aut}
\DeclareMathOperator{\Bir}{Bir}

\DeclareMathOperator{\Hom}{Hom}

\DeclareMathOperator{\Fix}{Fix}
\DeclareMathOperator{\Pic}{Pic}

\DeclareMathOperator{\rank}{rank}

\DeclareMathOperator{\Amp}{Amp}
\DeclareMathOperator{\Mov}{Mov}
\DeclareMathOperator{\Pos}{Pos}

\newcommand{\bz}{\mathbb{Z}}
\newcommand{\bq}{\mathbb{Q}}
\newcommand{\br}{\mathbb{R}}
\newcommand{\bc}{\mathbb{C}}

\newcommand{\bp}{\mathbb{P}}

\newcommand{\bv}{\mathbf{v}}
\newcommand{\ch}{\textup{\textsf{ch}}}
\newcommand{\td}{\textup{\textsf{td}}}
\newcommand{\Gr}{\textup{\textsf{Gr}}}
\newcommand{\forg}{\textup{\textsf{forg}}}
\newcommand{\sfinf}{\textup{\textsf{inf}}}

\newcommand{\HHom}{\shh\! om}

\newcommand{\defi}[1]{{\textbf{\emph{#1}}}}
\newcommand{\cate}[1]{{\textbf{\textup{#1}}}}
\newcommand{\funct}[1]{{\textup{\textbf{\textsf{#1}}}}}

\newcommand{\GLp}{\operatorname{GL^+}}
\newcommand{\grp}{{\tilde{\GLp}}(2,\br)}

\theoremstyle{plain}
\newtheorem{theorem}{Theorem}[section]
\newtheorem{corollary}[theorem]{Corollary}
\newtheorem{lemma}[theorem]{Lemma}
\newtheorem{proposition}[theorem]{Proposition}

\theoremstyle{definition}
\newtheorem{definition}[theorem]{Definition}
\newtheorem{remark}[theorem]{Remark}
\newtheorem{example}[theorem]{Example}

\author{Ziqi Liu}
\title{On Two Families of Enriques Categories over K3 Surfaces}
\address{Dipartimento di Matematica “F. Enriques”, Università degli Studi di Milano, Via Cesare Saldini 50, 20133 Milano, Italy.}
\email{ziqi.liu@unimi.it}



\thanks{The author is a member of the INdAM group GNSAGA (2025-2026).}

\begin{document}

\maketitle
\thispagestyle{plain}

\begin{abstract}
This article studies the moduli spaces of semistable objects related to two families of Enriques categories over K3 surfaces, coming from quartic double solids and special Gushel--Mukai threefolds. In particular, some classic geometric constructions are recovered in a modular way, such as the double EPW sextic and cube associated with a general Gushel--Mukai surface, and the Beauville's birational involution on the Hilbert scheme of two points on a quartic K3 surface. In addition, we describe the singular loci in some moduli spaces of semistable objects and an explicit birational involution on O'Grady's hyperkähler tenfold. Also, the appendix investigates the general theory of Enriques categories over K3 surfaces and provides a criterion for when an equivariant category of a K3 surface is an Enriques category.
\end{abstract}

\tableofcontents

\section{Introduction}

\subsection{Background and motivation}
The quotient of a complex K3 surface by a fixed-point free non-symplectic involution on it is an Enriques surface and every Enriques surface arises in this way. The study of moduli spaces of semistable objects on Enriques surfaces and the associated K3 surfaces produces rich geometric content, such as \cite{Beck20,Nu16,NY20}. However, these considerations only relate to fairly special K3 surfaces. 

On the other hand, the bounded derived category $\cate{D}^b(S)$ of a certain polarized K3 surface $S$ admits a $\bz/2\bz$ group action generated by a non-symplectic autoequivalence $\Pi$ on it such that the equivariant category $\cate{D}^b(S)^{\Pi}$ shares some properties of the bounded derived category of an Enriques surface. Such an equivariant category is called an \emph{Enriques category} \cite{BP23,PPZ26}. Unlike the Enriques surfaces, Enriques categories can be found over large families of K3 surfaces, such as the family of all smooth quartic surfaces. In this case, the study of moduli spaces of semistable objects on Enriques categories over K3 surfaces will provide insights into the geometry related to much more general K3 surfaces, compared to those on the Enriques surfaces.

The aim of this article is to explore this idea through studying moduli spaces of semistable objects related to two families of Enriques categories over K3 surfaces.

\subsection{Enriques categories over K3 surfaces}
The general theory of Enriques categories is developed in \cite{BP23,PPZ26}. In this article, we only treat geometric ones over K3 surfaces: a \emph{geometric triangulated category} is an admissible subcategory of the bounded derived category of a smooth proper variety and an \emph{Enriques category} $\sche$ over a K3 surface $S$ is a triangulated category which is equivalent to the equivariant category of the bounded derived category $\cate{D}^b(S)$ by a suitable involution $\Pi\colon\cate{D}^b(S)\rightarrow\cate{D}^b(S)$. The equivariant structure $\sche\cong \cate{D}^b(S)^{\Pi}$ in turn induces a unique residual involution $\funct{U}\colon\sche\rightarrow\sche$ such that $\cate{D}^b(S)\cong\sche^{\funct{U}}$ and vice versa.

Similar to ample line bundles on projective varieties, one has a notion of stability condition \cite{Bri07} which can be seen as the polarization on triangulated categories. This article will only consider \emph{numerical} stability conditions on some geometric category $\scht\subset\cate{D}^b(X)$ i.e.~those with central charge factoring through the \emph{numerical Grothendieck group} $N(\scht)\subset N(X)$. 

Under the equivariant construction (see Theorem~\ref{En-equivariant-stability}), a $\funct{U}$-fixed stability condition $\sigma_{\sche}$ on the geometric Enriques category $\sche$ induces a $\Pi$-fixed stability condition $\sigma_{S}$ on $\cate{D}^b(S)$ and vice versa. Once the two stability conditions are in addition \emph{proper}, meaning that the moduli stacks of semistable objects with a given numerical class and phase are proper algebraic stacks, one can establish the following theorem generalizing \cite[Proposition 7.2]{Nu16}.

\begin{theorem}[{{\cite[Theorem 1.8]{PPZ26}}}]
Consider a geometric Enriques category $\sche$ over a K3 surface $S$ characterized by the involutions $\funct{U}$ on $\sche$ and $\Pi$ on $\cate{D}^b(S)$, a $\funct{U}$-fixed proper stability condition $\sigma_{\sche}$ on $\sche$, and the induced proper stability condition $\sigma_{S}$ on $\cate{D}^b(S)$.

Then for any $\lambda\in N(\sche)$, the moduli stack $\mathfrak{M}_{\sigma_{\sche}}(\lambda)$ is smooth at any point corresponding to a $\sigma_{\sche}$-stable object that is not fixed by the $\funct{U}$-action on $\sche$. Suppose in addition that $\lambda$ is fixed by the group automorphism $\funct{U}_*$, then the forgetful functor $\cate{D}^b(S)\rightarrow\sche$ induces a surjective morphism 
$$\bigsqcup_{\forg(v)=\lambda}\mathfrak{M}_{\sigma_S}(v)\rightarrow\mathfrak{M}_{\sigma_{\sche}}(\lambda)^{\funct{U}}$$
where the target is the fixed locus of the $\funct{U}$-induced action on $\mathfrak{M}_{\sigma_{\sche}}(\lambda)$ and the union in the source ranges over $v \in N(S)$ that are sent to $\lambda$ via the forgetful functor. Moreover, the surjection is étale of degree $2$ at $\sigma_{S}$-stable objects not fixed by the involution $\Pi$. 
	
Similarly, for any $v \in N(S)$, the moduli stack $\mathfrak{M}_{\sigma_{S}}(v)$ is smooth at any point corresponding to a $\sigma_{S}$-stable object. Suppose that $v$ is fixed by the group homomorphism $\Pi_*$, the inflation functor $\sche \rightarrow \cate{D}^b(S)$ induces a surjective morphism 
$$
\bigsqcup_{\sfinf(\lambda) =v} \mathfrak{M}_{\sigma_{\sche}}(\lambda) \rightarrow \mathfrak{M}_{\sigma_{S}}(v)^{\Pi},
$$
where the target is the fixed locus of the $\Pi$-action on $\mathfrak{M}_{\sigma_{S}}(v)$ and the union in the source ranges over $\lambda \in N(\sche)$ such that $\sfinf(\lambda)=v$. 
Moreover, the surjection is étale of degree $2$ at $\sigma_{\sche}$-stable objects not fixed by the involution $\funct{U}$. 
\end{theorem}

Moreover, the proper algebraic stacks $\mathfrak{M}_{\sigma_{\sche}}(\lambda)$ and $\mathfrak{M}_{\sigma_S}(v)$ admit proper good moduli spaces say $M_{\sigma_{\sche}}(\lambda)$ and $M_{\sigma_S}(v)$ (see \cite{Al13} for definition) according to \cite{ALH23,BLMNPS21}. Also, one has the corresponding subspaces $M^{\stable}_{\sigma_{\sche}}(\lambda)$ and $M^{\stable}_{\sigma_S}(v)$ corresponding to the locus containing stable objects. Then the morphisms between moduli stacks descend to morphisms between moduli spaces.

In this article, we work out some explicit morphisms between moduli spaces for two families of Enriques categories over K3 surfaces. The moduli spaces involved on these Enriques categories are usually understood, so the actual work is to study the relevant moduli spaces on K3 surfaces and involutions. But a better knowledge of moduli spaces on the Enriques side can be obtained as well. In the appendix, we discuss the general theory of Enriques categories over K3 surfaces and two types of involutions on the bounded derived categories of K3 surfaces. The main results are summarized in the subsequent subsections.

\subsection{Kuznetsov components of quartic double solids}
A quartic double solid $X$ is a double cover of $\bp^3$ branched over a quartic K3 surface $S\subset\bp^3$ and it admits a natural ample polarization $\sho_X(1)$ which generates $\Pic(X)$. The admissible subcategory 
$$\Kuc(X)=\langle\sho_X,\sho_X(1)\rangle^{\perp}\subset \cate{D}^b(X)$$
is an Enriques category over the K3 surface $S$ given by an involution $\funct{U}$. The residual $\bz/2\bz$-action on $\cate{D}^b(S)\cong \Kuc(X)^{\funct{U}}$ is generated by some involution $\Pi_q$.

It is proved in Corollary~\ref{QDS-stability-condition-exists} that, via the equivariant construction, the induced $\Pi_q$-fixed stability condition $\sigma_S$ on $\cate{D}^b(S)$ of any $\funct{U}$-fixed stability condition $\sigma$ on $\Kuc(X)$ is equal to the same stability condition $\sigma_q$ up to a unique $\grp$-action. In particular, we recover the quartic double solid case of \cite[Theorem A.10]{JLLZ24}.

The numerical Grothendieck group $N(\Kuc(X))$ has rank 2, with a basis $\mu_1$, $\mu_2$ such that
$$\ch(\mu_1)=1-\frac{1}{2}h^2,\quad \ch(\mu_2)=h-\frac{1}{2}h^2-\frac{1}{3}h^3$$
where $h:=c_1(\sho_X(1))$. The group $N(\Kuc(X))$ is invariant under $\funct{U}_*$, while the $\Pi_q$-fixed lattice in $N(S)$ is spanned by the Mukai vectors $u_1=(1,0,-1)$ and $u_2=(1,-H,1)$, where $H$ is the ample polarization on $S$ defined by the restriction $\sho_S(H):=\sho_{\bp^3}(1)|_S$.

Now we take the $\funct{U}$-fixed stability condition $\sigma$ on $\Kuc(X)$ which exactly induces $\sigma_q$ via the equivariant construction. By studying the moduli spaces and morphisms related to $u_1$, we are able to recover some classic constructions in \cite{Bea83,Wel81}.

\begin{theorem}[Subsection \ref{QDS-ss-over-u1}]
The moduli space $M_{\sigma_q}(u_1)$ satisfies the following properties.
\begin{enumerate}
	\item It has finitely many singular points represented by $\sho_L(-2)\oplus\sho_S(-L)$ for each line $L\subset S$.
	\item One has a small contraction $S^{[2]}\rightarrow M_{\sigma_q}(u_1)$ which sends any length two subscheme $Z\subset L$ to the singular point represented by $\sho_L(-2)\oplus\sho_S(-L)$. This contraction induces an isomorphism on the stable locus $M^{\stable}_{\sigma_q}(u_1)$.
	\item Through this contraction, the involution on $M^{\stable}_{\sigma_q}(u_1)$ induced by $\Pi_q$ coincides with the birational Beauville involution on $S^{[2]}$, and the fixed locus $M_{\sigma_q}(u_1)^{\Pi_q}$ is isomorphic to the surface of bitangent lines to $S$.
	\item The moduli space $M_{\sigma}(-\mu_1)$ is the Hilbert scheme of lines on $X$ (see \cite{FP23,PY22}), and the induced surjection $M_{\sigma}(-\mu_1)\rightarrow M_{\sigma_q}(u_1)^{\Pi_q}$ is exactly the one in \cite[Page 15]{Wel81}.
\end{enumerate}
\end{theorem}

Besides, we acquire new knowledge of singular loci in moduli spaces on $\Kuc(X)$. The first one gives a negative answer to a question in \cite[Remark A.16]{FLZ24}.

\begin{theorem}[Subsection \ref{QDS-ss-over-2mu1}]
Suppose that $S$ is a generic quartic K3 surface, then the singular locus in $M^{\stable}_{\sigma}(2\mu_1)$ is the disjoint union of a point and a four-dimensional component.
\end{theorem}

\begin{theorem}[Subsection \ref{QDS-ss-over-mu2}]
The singular locus of $M_{\sigma}(\mu_2)$ is isomorphic to the union of the quartic K3 surface $S$ and finitely many disjoint points. The number of these points can be counted by enumerating certain classes in $\NS(S)$. Suppose in addition that $S$ is generic, then this number is zero and the singular locus of $M_{\sigma}(\mu_2)$ is isomorphic to $S$.
\end{theorem}

Moreover, we can give a birational involution for the ten-dimensional irreducible symplectic manifold first constructed by O'Grady \cite{O'Grady99}.

\begin{theorem}[Subsection \ref{QDS-ss-over-2u1}]
The involution $\Pi_q$ induces a birational involution $\pi_q$ on the singular moduli space $M_H(2u_1)$. This birational involution lifts to a birational involution on the ten-dimensional irreducible symplectic manifold $\tilde{M}_H(2u_1)$. 
\end{theorem}

In \cite[Section 3]{FMOS22}, an involution $\pi_D$ on the irreducible symplectic manifold $M_H(u_1-u_2)$ is investigated. Using wall-crossing, it turns out that $\pi_D$ is closely related to the $\Pi_q$-induced involution $\pi_q$ on the singular moduli space $M_{\sigma_q}(u_1-u_2)$.

\begin{theorem}[Subsection \ref{QDS-ss-u1minusu2}] The following statements are true.
\begin{enumerate}
	\item One has a contraction $M\rightarrow M_{\sigma_q}(u_1-u_2)$ from an irreducible symplectic manifold $M$.
	\item The irreducible symplectic manifold $M$ is birational to $M_H(u_1-u_2)$ via a flop.
	\item The birational involutions on $M$ induced by the involutions $\pi_D$ and $\pi_q$ coincide and can be extended to a regular involution on $M$.
\end{enumerate}
\end{theorem}
 
This relation could be helpful to understand the fixed locus inside $M_{\sigma_q}(u_1-u_2)$ and hence the covering morphism over it. 

\subsection{Kuznetsov components of special Gushel--Mukai threefolds}
A special Gushel--Mukai threefold $X$ is a double cover of a smooth codimension three linear section $W\subset\Gr(2,5)$ branched over a Brill--Noether general K3 surface $S$ of degree $10$. It admits semiorthogonal decompositions
$$\cate{D}^b(X)=\langle\schk_X^1,\sho_X,\shu_X^\vee\rangle\quad\textup{and}\quad \cate{D}^b(X)=\langle\schk_X^2,\shu_X,\sho_X\rangle$$
where $\shu_X$ is the restriction of the tautological vector bundle on $\Gr(2,5)$. The subcategories $\schk_X^1$ and $\schk_X^2$ are equivalent Enriques categories over the same K3 surface $S$ whose associated $\bz/2\bz$-actions are both generated by the double cover involution $X\rightarrow W$. The residual $\bz/2\bz$-actions on $\cate{D}^b(S)$ are generated by some involutions $\Pi_1$ and $\Pi_2$ respectively.

It is proved in Corollary~\ref{GM-stability-condition-exists} that, via the equivariant construction, the induced $\Pi_1$-fixed stability condition $\sigma_S$ on $\cate{D}^b(S)$ of any given $\funct{U}$-fixed stability condition $\sigma$ on $\schk_X^1$ is equal to the same stability condition $\sigma_1$ up to a unique $\grp$-action. In particular, we recover the special Gushel--Mukai threefold case of \cite[Theorem A.10]{JLLZ24}.

The numerical Grothendieck group $N(\schk_X^1)$ has rank 2, with a basis $\kappa_1$, $\kappa_2$ such that
$$\ch(\kappa_1)=1-\frac{1}{5}h^2,\quad \ch(\kappa_2)=2-h+\frac{1}{12}h^3$$
where $h:=c_1(\sho_X(1))$. The group $N(\schk_X^1)$ is invariant under $\funct{U}_*$, while the $\Pi_1$-fixed lattice in $N(S)$ is spanned by the Mukai vectors $w_1=(1,0,-1)$ and $w_2=(2,-H,2)$, where $H$ is the degree $10$ polarization on the Brill--Noether general K3 surface $S$.

Unlike the case for quartic double solids, we will only consider generic special Gushel--Mukai threefolds in the following results, in which case the K3 surface has Picard rank $1$. 

Now we take the $\funct{U}$-fixed stability condition $\sigma$ on $\schk_X^1$ which exactly induces $\sigma_1$ via the equivariant construction. By studying the moduli spaces and morphisms related to $w_1$, we are able to recover the birational involution in \cite[Section 4.3]{O'Grady05} and the covering morphism.

\begin{theorem}
[Subsection \ref{GM-ss-over-w1}]
The moduli space $M_{\sigma_1}(w_1)$ satisfies the following properties.
\begin{enumerate}
	\item It has a unique singular point represented by a direct sum of two stable objects.
	\item One has a contraction $S^{[2]}\rightarrow M_{\sigma_1}(w_1)$ which sends a closed locus of codimension at least $2$ to the unique singular point and restricts to an isomorphism on $M^{\stable}_{\sigma_1}(w_1)$.
	\item Through this contraction, the involution induced by $\Pi_1$ on $M^{\stable}_{\sigma_1}(w_1)\subset S^{[2]}$ coincides with the unique non-trivial birational involution on $S^{[2]}$ constructed in \cite[Section 4.3]{O'Grady05}. 
	\item The fixed locus $M_{\sigma_1}(w_1)^{\Pi_1}$ of the involution on $M_{\sigma_1}(w_1)$ induced by $\Pi_1$ is the union of a singular point and the surface $BTC_S$ of bitangent conics to $S$.
\end{enumerate}
\end{theorem}

A generic special Gushel--Mukai threefold can be canonically associated with a Lagrangian subspace $A$ of the vector space $\land^3 \bc^6$ due to \cite[Theorem 3.10]{DK18}. This gives a 4-dimensional irreducible symplectic manifold $\tilde{Y}_A$ called the \emph{double Eisenbud--Popescu--Walter sextic} (abbr.\ \emph{double EPW sextic}) \cite{O'Grady06}, and a 6-dimensional irreducible symplectic manifold $\tilde{K}_A$ called the \emph{double Eisenbud--Popescu--Walter cube} (abbr.\ \emph{double EPW cube}) \cite{IKKR19}.

The double EPW sextic admits a double covering map $\tilde{Y}_A\rightarrow Y_A$ which is branched along a surface $Y_A^{\geq 2}$, called the \emph{EPW surface}.

\begin{theorem}[Subsection \ref{GM-ss-over-w2}]
The covering $M_{\sigma}(-\kappa_2)\rightarrow M_{\sigma_1}(w_2)^{\Pi_1}$ is described as follows.
\begin{enumerate}
	\item The moduli space $M_{\sigma_1}(w_2)$ is isomorphic to the double EPW sextic $\tilde{Y}_A$ and the involution $\pi_1$ induced by $\Pi_1$ is the covering involution, so that $M_{\sigma_1}(w_2)^{\Pi_1}\cong Y_A^{\geq 2}$.
	\item The moduli space $M_{\sigma}(-\kappa_2)$ is isomorphic to the moduli space $M_X(2,1,5)$ of Gieseker semistable sheaves with $c_0=2,c_1=h,c_2=h^2/2,c_3=0$ due to \cite{JLLZ24}. The induced double covering onto the fixed locus $Y_A^{\geq 2}$ is étale. 
\end{enumerate}
\end{theorem}

\begin{remark}
It is expected that the covering map $M_{\sigma}(-\kappa_2)\rightarrow M_{\sigma_1}(w_2)^{\Pi_1}$ is the same as the one constructed in \cite{DK20}. The author would like to explore this in future work. 
\end{remark}

Similarly, we are able to recover the double EPW cube $\tilde{K}_A$ via studying the moduli spaces and morphisms related to the vector $w_1+w_2$ (one compares with \cite[Theorem 1.2]{FGLZ25}).

\begin{theorem}[Subsection \ref{GM-ss-over-w1plusw2}]
The moduli space $M_{\sigma_1}(w_1+w_2)$ is isomorphic to the double EPW cube $\tilde{K}_A$ and is birational to $S^{[3]}$ via a flopping wall. The induced involution $\pi_1$ on $\tilde{K}_A$ by the autoequivalence $\Pi_1$ is the covering involution of $\tilde{K}_A\rightarrow K_A$ and it gives via $\tilde{K}_A\dashrightarrow S^{[3]}$ the unique non-trivial birational involution on $S^{[3]}$ described in \cite[Example 4.12]{De22}.
\end{theorem}

\subsection{More on Enriques categories over K3 surfaces}
In the Appendix \ref{ST-Hodge-Involution}, we discuss the general theory of Enriques categories over K3 surfaces. One recalls that an Enriques category $\sche$ over a K3 surface $S$ is equivalent to $\cate{D}^b(S)^{\Pi}$, where $\Pi$ is a suitable involution on $\cate{D}^b(S)$. Here we make this precise. At first, an autoequivalence $\Phi$ on $\cate{D}^b(S)$ induces a Hodge isometry on the Mukai lattice $\tilde{H}(S,\bz)$ (see,~e.g.~\cite[Section 10.2]{Hu}), and $\Phi$ is said to be \emph{anti-symplectic} if the induced Hodge isometry does not preserve the symplectic class in $H^2(S,\bz)$.

\begin{proposition}[Proposition~\ref{A_definition-Enriques}]
Consider a $\bz/2\bz$-action on $\cate{D}^b(S)$ for a K3 surface $S$ generated by a non-symplectic involution $\Pi$, then $\sche=\cate{D}^b(S)^{\Pi}$ is an Enriques category.
\end{proposition}

Then it is natural to ask what kinds of non-symplectic autoequivalences on the bounded derived category $\cate{D}^b(S)$ would give Enriques categories. The two families of Enriques categories considered in this article originate from \cite{KP17}. Motivated by their constructions, we consider autoequivalences on $\cate{D}^b(S)$ of a polarized K3 surface $(S,H)$ with the following two forms
$$\Pi=\funct{O}^n[-1]\quad\textup{and}\quad \Pi=\funct{T}_U\circ\funct{O}[-1]$$
where $\funct{O}:=\funct{T}_{\sho_S}(-\otimes\sho_S(H))$ and $U$ is a spherical object on $\cate{D}^b(S)$, corresponding respectively to the scenarios that $\cate{D}^b(X)=\langle\cate{D}^b(S)^{\Pi},\sho_X,\cdots,\sho_X(n)\rangle$ and $\cate{D}^b(X)=\langle\cate{D}^b(S)^{\Pi},\sho_X,U_X\rangle$.

It turns out that these autoequivalences are involutive if and only if we are in one of the two cases discussed in this article, according to Proposition~\ref{A_involution-type-1} and Proposition~\ref{A_involution-type-2}.

\subsection{Conventions and organization}
In this article, schemes and morphisms of schemes are assumed to be over the field of complex numbers. Also, categories and functors are assumed to be additive and complex linear. Although it is not strictly necessary, we will suppose that all categories are essentially small to avoid set-theoretic problems.

Section~\ref{ST-group-actions-and-EnCat} first recalls some generalities about group actions and equivariant categories following \cite{BO23,Ela15}. Then, the section contains a summary of some contents in \cite{BP23,PPZ26}, briefly introducing the notion of a geometric Enriques category over a K3 surface.

Section~\ref{ST-stability-condition-gtc} collects some contents from \cite{PPZ26}, which deals with stability conditions and moduli spaces of semistable objects on Enriques categories. Moreover, equivariant stability conditions and the relation between moduli spaces are discussed following \cite{PPZ26}.

Section~\ref{ST-stuff-on-K3} contains a brief introduction to two notions of stability on K3 surfaces, the relevant moduli spaces and their wall-crossing behaviors.

Section~\ref{ST-Family-QDS} investigates Kuznetsov components of quartic double solids. In particular, some birational involutions of moduli spaces and covering maps between moduli spaces are figured out. Similar considerations are elaborated on in Section~\ref{ST-Family-special-GM} for geometric Enriques categories constructed from special Gushel--Mukai threefolds.

Appendix \ref{ST-Hodge-Involution} discusses the general theory of Enriques categories over K3 surfaces and two types of Hodge isometries on polarized K3 surfaces. 

\subsection{Acknowledgments}
The author would like to thank his advisor Laura Pertusi for her consistent guidance throughout the preparation of this article, thank his advisor Paolo Stellari for suggestions and encouragement, thank Shizhuo Zhang for explanations on his works and remarks on this article, thank Alex Perry for helpful comments, and thank Yu-Wei Fan, Annalisa Grossi, Kuan-Wen Lai for answering the author's questions. Also, the author would like to thank the anonymous referee for careful reading and improvements.

\section{Geometric Enriques categories over K3 surfaces}\label{ST-group-actions-and-EnCat}

\subsection{Group actions and equivariant objects}
We fix a finite group $G$ and a category $\schc$.

\begin{definition}
A \defi{(right) action} of $G$ on $\schc$ is an assignment $\funct{P}_{(-)}\colon G\rightarrow\Aut(\schc)$ and a set of natural transformations $\{\epsilon_{g,h}\colon \funct{P}_{gh}\stackrel{\sim}{\rightarrow}\funct{P}_h\circ \funct{P}_g\,|\,g,h\in G\} $ such that
\begin{equation*}
	\xymatrix{
		\funct{P}_{g_1g_2g_3} \ar[d]_{\epsilon_{g_1,g_2g_3}} \ar[rr]^{\epsilon_{g_1g_2,g_3}} & & \funct{P}_{g_3} \circ \funct{P}_{g_1g_2} \ar[d]^{\funct{P}_{g_3}  \epsilon_{g_1,g_2}} \\
		\funct{P}_{g_2g_3} \circ \funct{P}_{g_1} \ar[rr]^{\epsilon_{g_2,g_3} \funct{P}_{g_1}} & & \funct{P}_{g_3} \circ \funct{P}_{g_2} \circ \funct{P}_{g_1}
	}
\end{equation*}
commutes for all $g_1,g_2,g_3 \in G$.
\end{definition}

Given a group action of $G$ on $\schc$, one can always construct a new category $\schc^G$, called the \emph{equivariant category} of the group action $G$ on $\schc$, whose objects are $G$-equivariant objects and whose morphisms are $G$-equivariant morphisms.

\begin{definition}
A \defi{$G$-equivariant object} of $\schc$ is a pair $(A, \phi)$, where $A$ is an object of $\schc$ and $\phi=\{\phi_g\}_{g\in G}$ is a collection of isomorphisms $\phi_{g}: A \xrightarrow{\sim}\funct{P}_gA$ for $g \in G$, such that the diagram
\begin{displaymath}
	\xymatrix{
		A \ar[rr]^{\phi_{h}} \ar[drrrr]_{\phi_{gh}} & & \funct{P}_h(A) \ar[rr]^{\funct{P}_h(\phi_{g})} 
		& & \funct{P}_h(\funct{P}_g(A)) \\
		& &  & & \funct{P}_{gh}(A) \ar[u]_{\epsilon_{g,h}(A)} 
	}
\end{displaymath}
commutes for all $g, h \in G$. The equivariant structure $\phi$ is often omitted from the notation. 
\end{definition}

\begin{definition}
A \defi{$G$-equivariant morphism} $(A,\phi^A)\rightarrow (B,\phi^B)$ of $G$-equivariant objects of $\schc$ is a morphism $f\colon A\rightarrow B$ in $\schc$ such that the diagram
\begin{displaymath}
	\xymatrix{
		A\ar[d]_f\ar[rr]^{\phi^A_g}&&\funct{P}_g(A)\ar[d]_{\funct{P}_g(f)}\\
		B\ar[rr]^{\phi^B_g}&&\funct{P}_g(B)
}\end{displaymath}
commutes for every $g\in G$.
\end{definition}

\begin{remark}
Equivariant categories in this article are automatically additive and complex linear under our convention that categories and functors are additive and complex linear. 
\end{remark}

\subsection{Reconstruction theorem}
We consider a category $\schc$, a finite abelian group $G$, and the dual group $G^\vee=\Hom(G,\bc^{\times})$ of $G$. Suppose that $G$ acts on $\schc$, then we can define a natural action of $G^\vee$ on $\schc^G$ by twisting: for $\chi\in G^\vee$ and $(A,\phi)\in\schc^G$ we let
$$\psi_{\chi}(A,\{\phi_g\}_{g\in G})=(A,\{\phi_g\}_{g\in G})\otimes\chi=(A,\{\phi_g\cdot\chi(g)\})$$
The equivariant objects $\psi_{\chi_1}(\psi_{\chi_2}(A,\phi))$ and $\psi_{\chi_2\chi_1}(A,\phi)$ are the same, so we can just let the required natural transformations $\psi_{\chi_1}\circ\psi_{\chi_2}\rightarrow\psi_{\chi_2\chi_1}$ be identities.

\begin{theorem}[{{\cite[Theorem 4.2]{Ela15}}}]\label{Ela-reconstruction}
Consider a finite abelian group $G$ acting on an idempotent complete category $\schc$, then there exists an equivalence
$$\funct{X}\colon(\schc^G)^{G^\vee}\cong\schc$$ 
such that the $G$-action is identified with the $(G^\vee)^\vee\cong G$-action.
\end{theorem}

There are mutually left and right adjoint functors: the \emph{forgetful functor} 
$$\funct{Forg}\colon\schc^G\rightarrow\schc$$ 
defined by forgetting the $G$-equivariant structure, and the \emph{inflation functor} 
$$\funct{Inf}\colon\schc\rightarrow\schc^G$$ 
given on objects by sending $A$ to $\bigoplus_{g \in G} \funct{P}_gA$ with the natural $G$-equivariant structure (one can see, for example, \cite[Lemma 3.8]{Ela15}). Moreover, one has the following observation:  

\begin{lemma}\label{Ela-reconstruction-lemma}
Under the assumptions of Theorem~\ref{Ela-reconstruction}, there is an equivalence of functors
$$
\funct{Forg} \cong \funct{Inf} \circ \funct{X} \colon (\schc^G)^{G^{\vee}} \rightarrow \schc^G
$$
where $\funct{Forg}$ is the one for the induced $G^\vee$-action on $\schc^G$.
\end{lemma}

Suppose in addition that both $\schc$ and $\schc^G$ are triangulated categories, then $(\schc^G)^{G^\vee}$ comes with a natural class of exact triangles, consisting of those triangles whose images under 
$$\funct{Forg}\colon (\schc^G)^{G^\vee} \rightarrow \schc^G$$ 
are exact. Then, due to Lemma~\ref{Ela-reconstruction-lemma}, the equivalence $\funct{X}$ in Theorem~\ref{Ela-reconstruction} respects exact triangles and is therefore an equivalence of triangulated categories.

\subsection{Geometric Enriques categories over K3 surfaces}
The theory of Enriques categories is developed in \cite{BP23,PPZ26} using the language of stable infinity categories. As the article only treats the geometric ones over K3 surfaces, we will refrain from introducing the whole story. 

\begin{definition}
The \defi{Serre functor} $\funct{S}_{\schc}$ of a category $\schc$ is an autoequivalence such that there exists a canonical isomorphism of complex vector spaces
$$\Hom(A,B)\cong\Hom(B,\funct{S}_{\schc}A)^\vee$$
for any objects $A,B\in\schc$. The dual on the right-hand side is taken for the vector space structure.
\end{definition}

The Serre functor, once it exists, is unique up to a unique isomorphism of functors.

\begin{definition}
A triangulated category is called \defi{geometric} if it is equivalent to an admissible subcategory of the bounded derived category $\cate{D}^b(X)$ of a smooth proper variety $X$. 
\end{definition}

The category $\cate{D}^b(X)$ always admits a Serre functor $\funct{S}_X\cong(-\otimes\omega_X)[\dim{X}]$, which induces a Serre functor for any geometric triangulated category. 

\begin{definition}
A \defi{(geometric) Enriques category $\sche$ over a K3 surface $S$} is a geometric triangulated category equipped with a $\bz/2\bz$-action generated by $\funct{U}$, such that $\funct{U}\cong\funct{S}_{\sche}\circ[-2]$ is not the identity in $\Aut(\sche)$ and $\sche^{\funct{U}}\cong\cate{D}^b(S)$ as triangulated categories.
\end{definition}

The generator of the dual group action on $\cate{D}^b(S)$ is usually denoted by $\Pi$ in this article.

\begin{example}
The bounded derived category $\cate{D}^b(T)$ of an Enriques surface $T$ is an Enriques category. The group action is generated by $\funct{U}=(-\otimes\omega_T)$. The category $\cate{D}^b(T)^{\funct{U}}$ is equivalent to the bounded derived category $\cate{D}^b(S)$ of its K3 cover, and the residual $\bz/2\bz$-action on $\cate{D}^b(S)$ is induced by the covering involution of the double cover $S\rightarrow T$.
\end{example}

The following examples are applications of \cite{Ku19,KP17} and we refer the readers to \cite{PPZ26} for more details and explanations. 

\begin{example}\label{En_example-1}
A quartic double solid $X$ is a double cover of $\bp^3$ branched over a quartic K3 surface $S\subset\bp^3$. The derived category of $X$ admits a semiorthogonal decomposition
$$\cate{D}^b(X)=\langle\Kuc(X),\sho_X,\sho_X(h)\rangle$$
where $\sho_X(h)$ is the pullback bundle from $\sho_{\bp^3}(1)$. Then the admissible component $\Kuc(X)$ is an Enriques category over the K3 surface $S$ with the group action induced by the double cover involution $X\rightarrow\bp^3$. The residual $\bz/2\bz$-action on $\cate{D}^b(S)\cong \Kuc(X)^{\funct{U}}$ is generated by 
$$\Pi_q:=\funct{T}_{\sho_S}\circ(-\otimes\sho_S(H))\circ\funct{T}_{\sho_S}\circ(-\otimes\sho_S(H))[-1]$$
where $\sho_S(H)=\sho_{\bp^3}(1)|_S$ and $\funct{T}_{\sho_S}$ is the spherical twist at the spherical object $\sho_S$.
\end{example}

\begin{example}\label{En_example-2}
A special Gushel--Mukai threefold $X$ is a double cover of a smooth codimension three linear section $W$ of $\Gr(2,5)$ branched over a smooth divisor $S\in|\sho_W(2D)|$, where $D$ comes from the ample generator of $\Pic(\Gr(2,5))$. By abuse of notation, we use $H$ to denote the ample class on $S$ coming from $\Gr(2,5)$. Then $(S,H)$ is a Brill--Noether general K3 surface of degree 10.

One has two semiorthogonal decompositions
$$\cate{D}^b(X)=\langle\schk_X^1,\sho_X,\shu_X^\vee\rangle\quad\textup{and}\quad \cate{D}^b(X)=\langle\schk_X^2,\shu_X,\sho_X\rangle$$
where $\shu_X$ is the tautological vector bundle pulled back from that of $\Gr(2,5)$. Then the admissible components $\schk_X^1$ and $\schk_X^2$ are equivalent Enriques categories, with the same action induced by the double cover. The corresponding equivariant categories are both equivalent to $\cate{D}^b(S)$.

However, the residual actions are generated by different autoequivalences:
\begin{align*}
\Pi_1:=&\funct{T}_{\sho_S}\circ\funct{T}_{\shu_S^\vee}\circ(-\otimes\sho_S(H))[-1]\\
\Pi_2:=&\funct{T}_{\shu_S}\circ\funct{T}_{\sho_S}\circ(-\otimes\sho_S(H))[-1]
\end{align*}
where $\shu_S$ is the unique $H$-stable spherical vector bundle on $S$ with Mukai vector $(2,-H,3)$.
\end{example}

\section{Stability conditions on geometric triangulated categories}\label{ST-stability-condition-gtc}
\subsection{Stability conditions}
Here we introduce full numerical stability conditions on geometric triangulated categories, originated from the seminal work of Bridgeland \cite{Bri07}.

\begin{definition}
A \defi{stability function} on an abelian category $\scha$ is a group homomorphism $Z\colon K(\scha)\rightarrow\bc$ such that $\Im Z(E) \geq 0$, and in the case that $\Im Z(E) = 0$, we have $\Re Z(E)< 0$.
\end{definition}

The \emph{slope} $\mu_Z\colon K(\scha)\rightarrow\br\cup\{+\infty\}$ associated with a stability function $Z\colon K(\scha)\rightarrow\bc$ is defined by $\mu_Z(E)= -\Re Z(E)/\Im Z(E)$ for $E$ when $\Im Z(E)>0$, and $\mu_Z(E)=+\infty$ otherwise. A non-zero object $E$ is said to be \emph{semistable} (resp.\ \emph{stable}) with respect to the stability function $Z$ if for every proper subobject $A\subset E$ one has $\mu_Z(A)\leq \mu_Z(E/A)$ (resp.\ $\mu_Z(A)<\mu_Z(E/A)$).

It is sometimes more convenient to use a notion of phase to arrange semistable objects associated with a stability function $Z$. The \emph{phase} $\phi_Z\colon K(\scha)\rightarrow(-1,1]$ associated with $Z$ is defined by $\phi_Z(E)=(1/\pi)\arg Z(E)$ for $E$ when $\Im Z(E)>0$, and $\phi_Z(E)=1$ otherwise.

The subscripts for $\mu_Z$ and $\phi_Z$ associated with a stability function $Z$ will usually be omitted once the stability function is clear from the context.

\begin{definition}
Let $Z$ be a stability function on an abelian category $\scha$, then a non-zero object $E$ of $\scha$ is said to admit a \defi{Harder--Narasimhan (HN) filtration} if there exists a filtration
$$0=E_0 \hookrightarrow E_1 \hookrightarrow \dots E_{n-1} \hookrightarrow E_n=E$$
whose factors $A_i:=E_i/E_{i-1} \neq 0$ are semistable with respect to $Z$ and $\mu(A_1) > \dots > \mu(A_n)$. 
\end{definition}

Once a Harder--Narasimhan filtration exists, it is unique up to isomorphisms. 

\begin{definition}
Let $\schc$ be a geometric triangulated category, $N(\schc)$ be its numerical Grothendieck group, and $\bv\colon K(\schc)\rightarrow N(\schc)$ be the obvious map. Then a \defi{(numerical) stability condition} $\sigma$ on $\schc$ is a pair $(\scha, Z)$, where $\scha$ is the heart of a bounded $t$-structure on $\schc$ and $Z\colon N(\schc)\rightarrow \bc$ is a group homomorphism, such that the composition 
$$K(\scha)= K(\schc)\stackrel{\bv}{\rightarrow} N(\schc)\stackrel{Z}{\rightarrow}\bc$$ 
is a stability function on the abelian category $\scha$, and every non-zero object $E$ in $\scha$ admits a (unique) Harder--Narasimhan filtration with respect to the stability function.

A non-zero object $E$ of $\schc$ is said to be \defi{$\sigma$-semistable} (resp.\ \defi{$\sigma$-stable}) if a shift of $E$ is semistable (resp.\ stable) with respect to the stability function $Z\circ\bv$. 
\end{definition}

The slope of $\sigma$-semistable objects can only be defined for those in hearts, so it is often convenient to order the $\sigma$-semistable objects by their phase which is $\phi_{\sigma}(E)=\phi_Z(E[n])-n$ for a $\sigma$-semistable object $E$ with $E[n]$ being in the heart $\scha$.

It induces a collection $\{\schp(\phi)\}_{\phi\in\br}$ of full additive subcategories of $\schc$, such that $\schp(\phi)$ consists of the zero object of $\schc$ and all $\sigma$-semistable objects with phase $\phi$. Moreover, we can define for each interval $I \subset \br$ a subcategory $\schp(I)$ as the extension-closed subcategory of $\schc$ generated by all the subcategories $\schp(\phi)$ with $\phi \in I$. One notices that $\schp((0, 1])= \scha$.

\begin{remark} 
The collection of subcategories $\schp=\{\schp(\phi)\}_{\phi \in \br}$ forms a \emph{slicing}, and the definition of a stability condition can be rephrased as the data of a pair $(\schp, Z)$ of a slicing $\schp$ and a suitably compatible central charge $Z$ on $\schp((0,1])$ according to \cite[Proposition 5.3]{Bri07}.  
\end{remark} 

\subsection{The complex manifold of stability conditions}
A stability condition $\sigma$ is called \emph{full} if there exists a quadratic form $Q$ on $N(\schc)_{\br}$ such that the restriction of $Q$ to the kernel of $Z(-)_{\br}$ is negative definite and $Q(E) \geq 0$ for any semistable object $E$ in $\scha$ with respect to $Z\circ\bv$. 

\begin{remark}
The condition of being full is also called the \emph{support property} and has several equivalent formulations (see, for example, appendices of \cite{BM11,BMS16}). This property is essential for the following wall and chamber structure. Henceforth in this article, the term \emph{stability condition} means a full (numerical) one except for the special cases mentioned.
\end{remark}

The set $\Stab(\schc)$ of all stability conditions on a geometric triangulated category $\schc$ admits a structure of complex manifold due to \cite{Bay19,Bri07}. In particular, a stability condition can be deformed uniquely given a small deformation of the central charge. Moreover, one has

\begin{proposition}[{{\cite[Proposition 2.3]{BM14}, \cite[Section 9]{Bri08}, \cite[Section 3.2]{Li25}}}]
Consider a class $v$ in $N(\schc)$, then there exists a locally finite set of walls (real codimension one submanifolds with boundary) in the complex manifold $\Stab(\schc)$, depending only on $v$, such that:
\begin{itemize}
	\item The set of stable objects of class $v$ does not change in a chamber, and the set of semistable objects of class $v$ does not change for any generic stability condition in a chamber;
	\item When $\sigma$ lies on a single wall in $\Stab(\schc)$, then there is a $\sigma$-semistable object that is unstable in one of the adjacent chambers, and semistable in the other adjacent chamber.
\end{itemize}
In particular, a stability condition $\sigma$ lies on a wall with respect to $v$ only if there exists a strictly $\sigma$-semistable object of class $v$. 
\end{proposition}

\begin{definition}
A stability condition $\sigma$ on $\schc$ is called \emph{generic} with respect to a given numerical class $v\in N(\schc)$ if it does not lie on a wall with respect to $v$.
\end{definition}

According to \cite[Lemma 8.2]{Bri07}, there are two mutually-commutative group actions on the complex manifold $\Stab(\schc)$ as follows.

The group $\Aut(\schc)$ acts on $\Stab(\schc)$ on the left, through
$$\Phi. (\scha, Z):= (\Phi(\scha), Z \circ \Phi_*^{-1})  $$
where $\Phi_*$ denotes the induced automorphism of the group $N(\schc)$.

The universal cover 
$$
\grp  = \left \{ \tilde{g}=(M, f) ~  \middle\vert ~ 
\begin{array}{l} 
	M \in \GLp(2,\br), f \colon \br \rightarrow \br \text{ is an increasing function} \\ 
	\text{such that for all $\phi \in \br$ we have $f(\phi+1) = f(\phi) + 1$} \\
	\text{and $M \cdot e^{i\pi \phi} \in \br_{>0} \cdot e^{i \pi f(\phi)}$} 
\end{array}  \right \}
$$
of the group $\GLp(2,\br)=\{M\in\GL(2,\br)\,|\,\det(M)>0\}$ acts on the right on $\Stab(\schc)$, via 
$$
(\scha, Z) . (M, f):= (\schp((f(0), f(1)]), M^{-1} \circ Z). 
$$ 
where we identify $\bc$ with $\br\oplus\br\sqrt{-1}$ to make sense of the composite $M^{-1}\circ Z\colon N(\schc)\rightarrow\bc$.

\subsection{Moduli stacks of semistable objects}
Let $\sigma=(\scha,Z)$ be a stability condition on a geometric triangulated category $\schd\subset\cate{D}^b(X)$ and $v$ be a class in $N(\schd)$. Then one can define a moduli stack $\mathfrak{M}_{\sigma}(v)$ of $\sigma$-semistable objects with class $v$. A comprehensive reference is \cite{BLMNPS21}.

\begin{definition}
The moduli stack $\mathfrak{M}_{\sigma}(v)\colon(\cate{Sch}/\bc)^{\cate{op}}\rightarrow\cate{Gpds}$ is defined by setting 
$$
\mathfrak{M}_{\sigma}(v)(T)= \left \{ E\in\cate{D}^b(X_T) ~  \middle\vert ~ 
\begin{array}{l} 
	\,E|_{X_t}\textup{ is in }\schd_t\subset\cate{D}^b(X_t)\textup{ and} \\ 
	\textup{is }\sigma\textup{-semistable with class }v \\
	\textup{for any closed point }t\in T 
\end{array}  \right \}
$$
for any locally finitely generated complex scheme $T$. Similarly, one defines $\mathfrak{M}^{\stable}_{\sigma}(v)$.
\end{definition}

A stability condition defined in \cite{BLMNPS21} satisfies additional requirements and will be called \emph{proper} in this article. The key property of being proper is the following. 

\begin{theorem}[{{\cite[Theorem 21.24]{BLMNPS21}}}]
Consider a proper stability condition $\sigma$ on a geometric triangulated category $\schd$, then $\mathfrak{M}_{\sigma}(v)$ is an algebraic stack of finite type over $\bc$ and there is an open immersion $\mathfrak{M}_{\sigma}^{\stable}(v)\subset \mathfrak{M}_{\sigma}(v)$. Moreover, the moduli stack $\mathfrak{M}_{\sigma}(v)$ admits a proper good moduli space $M_{\sigma}(v)$ such that $\mathfrak{M}_{\sigma}^{\stable}(v)$ is a $\mathbb{G}_m$-gerbe over the counterpart $M^{\stable}_{\sigma}(v)\subset M_{\sigma}(v)$.
\end{theorem}

Such a good moduli space (in the sense of \cite{Al13}) say $M_{\sigma}(v)$ is an algebraic space, whose points are $S$-equivalence classes of semistable objects in $\schp((-1,1])$ with numerical class $v$.

\subsection{Invariant stability conditions}
Stability conditions on triangulated categories are also compatible with equivariant constructions. Here we collect some essential results from \cite{PPZ26} but restrict to geometric triangulated categories and (full numerical) stability conditions on them.

Given an autoequivalence $\Phi$ of a triangulated category $\schc$, then a stability condition $\sigma$ is called \emph{$\Phi$-fixed} (resp.\ \emph{$\Phi$-invariant}) if $\Phi.\sigma=\sigma$ (resp.\ $\Phi.\sigma=\sigma.\tilde{g}$) for some $\tilde{g}\in\grp$.

\begin{theorem}[\cite{MMS09,PPZ26,Pol07}]\label{En-equivariant-stability}
Consider a geometric triangulated category $\schc$ and a finite group action $G$ on $\schc$ such that $\schc^G$ is also geometric, and a stability condition $\sigma=(\scha,Z)$ on $\schc$ fixed by the action of any element in $G$, then the following data gives a stability condition $\sigma^G$ on $\schc^G$:
\begin{align*}
\scha^G&=\{(A,\phi)\in\schc^G\,|\,A\in\scha\}\\
Z^G&= Z\circ\forg\colon N(\schc^G)\rightarrow\bc
\end{align*}
where the homomorphism $\forg\colon N(\schc^G)\rightarrow N(\schc)$ is induced by the forgetful functor $\schc^G\rightarrow\schc$. 

Moreover, for an object $(A,\phi)\in\schc^G$, one has
\begin{itemize}
	\item $(A,\phi)$ is $\sigma^G$-semistable if and only if $A$ is $\sigma$-semistable;
	\item $(A,\phi)$ is $\sigma^G$-stable implies that $A$ is $\sigma$-poly-stable;
	\item $A$ is $\sigma$-stable implies that $(A,\phi)$ is $\sigma^G$-stable;
\end{itemize}
where a $\sigma$-poly-stable object is a direct sum of $\sigma$-stable objects of the same phase.
\end{theorem}

As before, the induced stability condition $\sigma^G$ will be denoted by $\sigma^{\Phi}$ if we have chosen a specific generator $\Phi$ of the cyclic group action by $G$.

\begin{remark}
The category $\scha^G$ is equivalent to the category of $G$-equivariant objects for the induced group action on $\scha$ by construction. 
\end{remark}

The $G^\vee$-action on $\schc^G$ induces a trivial group action on $N(\schc^G)$, thus the stability condition $\sigma^G$ constructed from Theorem~\ref{En-equivariant-stability} is fixed by the action of any element of $G^\vee$. Therefore, one can apply the theorem again to obtain a stability condition $(\sigma^G)^{G^{\vee}}$ on $(\schc^G)^{G^{\vee}}$. 

\begin{proposition}[{{\cite[Lemma 4.11]{PPZ26}}}]\label{En-back-to-stability}
In the setup of Theorem~\ref{En-equivariant-stability}, the induced stability condition $(\sigma^{G})^{G^{\vee}}$ is identified with the stability condition  $|G| \cdot \sigma = (\scha, |G|\cdot Z)$ on the triangulated category $\schc$ under the equivalence $\funct{X}\colon (\schc^G)^{G^{\vee}} \xrightarrow{\sim} \schc$ stated in Theorem~\ref{Ela-reconstruction}. 
\end{proposition} 

In addition, we would like to point out the following proposition which might have already been known to experts.

\begin{proposition}\label{En-invariant-grp-action}
In the setup of Theorem~\ref{En-equivariant-stability}, one has $(\sigma.\tilde{g})^G=\sigma^G.\tilde{g}$ for any $\tilde{g}\in\grp$.

\begin{proof}
One can see $(\sigma.\tilde{g})^G=((\schp(f(0),f(1)])^G,M^{-1}\circ Z\circ\forg)=\sigma^G.\tilde{g}$ for $\tilde{g}=(M,f)$.
\end{proof}
\end{proposition}

\subsection{Moduli spaces of semistable objects on Enriques categories}
Choose an Enriques category $(\sche,\funct{U})$ over a K3 surface $S$, assume a $\funct{U}$-fixed stability condition $\sigma\in\Stab(\sche)$, and take the induced stability condition $\sigma^{\funct{U}}$ on $\cate{D}^b(S)$. As we have mentioned before, the induced stability condition $\sigma^{\funct{U}}$ is $\Pi$-fixed. So we deduce the following statement.

\begin{proposition}\label{ME-involution}
In the above setup, the autoequivalence $\funct{U}$ of $\sche$ induces an involution 
$$U\colon \mathfrak{M}_{\sigma}(\lambda)\rightarrow\mathfrak{M}_{\sigma}(\lambda) $$
of stacks, for any $\lambda\in N(\sche)$ fixed by the induced homomorphism $\funct{U}_*\colon N(\sche)\rightarrow N(\sche)$.

Moreover, the autoequivalence $\Pi$ of $\cate{D}^b(S)$ induces an involution 
$$\pi\colon \mathfrak{M}_{\sigma^{\funct{U}}}(v)\rightarrow\mathfrak{M}_{\sigma^{\funct{U}}}(v) $$
of stacks, for any $v\in N(S)$ fixed by the induced homomorphism $\Pi_*\colon N(S)\rightarrow N(S)$.
\end{proposition}

The fixed loci of these induced involutions are the moduli stacks of $\bz/2\bz$-equivariant objects in corresponding moduli stacks. Moreover, one has the following proposition.

\begin{proposition}\label{ME-double-cover}
In the above setup, the forgetful functor $\funct{Forg}$ induces a surjection
\begin{align*}
\bigsqcup_{\forg(v)=\lambda}\mathfrak{M}_{\sigma^{\funct{U}}}(v)\rightarrow\mathfrak{M}_{\sigma}(\lambda)^{\funct{U}}
\end{align*}
for any $\lambda\in N(\sche)$ fixed by the homomorphism $\funct{U}_*\colon N(\sche)\rightarrow N(\sche)$, where $\mathfrak{M}_{\sigma}(\lambda)^{\funct{U}}$ is the fixed locus of the induced $\funct{U}$-action on $\mathfrak{M}_{\sigma}(\lambda)$ and $\forg\colon N(S)\rightarrow N(\sche)$ is the induced map by $\funct{Forg}$.

Similarly, the inflation functor $\funct{Inf}$ induces a surjection
\begin{align*}
	\bigsqcup_{\sfinf(\lambda)=v}\mathfrak{M}_{\sigma}(\lambda)\rightarrow\mathfrak{M}_{\sigma^{\funct{U}}}(v)^{\Pi}
\end{align*}
for any $v\in N(S)$ fixed by the homomorphism $\Pi_*\colon N(S)\rightarrow N(S)$, where $\sfinf\colon N(\sche)\rightarrow N(S)$ is the homomorphism induced by $\funct{Inf}$ and $\mathfrak{M}_{\sigma^{\funct{U}}}(v)^{\Pi}$ is similarly defined.

\begin{proof}
It is a special case of \cite[Proposition 5.7, Remark 5.9]{PPZ26}. 
\end{proof}
\end{proposition}

\begin{remark}
A similar statement can be found in \cite[Proposition 3.9]{BO22}. One should notice that the notation $\mathfrak{M}_{\sigma}(-)^{G}$ has a different meaning in their article.
\end{remark}

The restriction of both surjections to stable loci on sources has degree $2$ cf.\ \cite[Lemma 1]{Plo07}. Suppose in addition that the stability conditions $\sigma$ and $\sigma^{\funct{U}}$ are proper.

\begin{proposition}[{{\cite[Corollary 5.6]{PPZ26}}}]\label{ME-smooth-points}
Consider an Enriques category $(\sche,\funct{U})$ over a K3 surface $S$ and a $\funct{U}$-invariant proper stability condition on $\sche$, then the moduli stack $\mathfrak{M}_{\sigma}(\lambda)$ is smooth at the locus of $\sigma$-stable points that are not fixed by the autoequivalence $\funct{U}$.
\end{proposition}

On the K3 side, there is also a smoothness criterion which asserts that the moduli stack $\mathfrak{M}_{\sigma^{\funct{U}}}(v)$ is smooth at the stable locus $\mathfrak{M}^{\stable}_{\sigma^{\funct{U}}}(v)$. Hence \cite[Proposition 5.7.(2)]{PPZ26} implies

\begin{corollary}\label{ME-double-cover-stable}
In the setup of Proposition~\ref{ME-double-cover} and suppose more that the stability conditions are proper, then restrictions of the surjections to stable loci in the source are degree $2$ onto their images and étale at the points that are not fixed by the corresponding $\bz/2\bz$ action.
\end{corollary}

\section{Stability on K3 surfaces, moduli spaces, and wall-crossing}\label{ST-stuff-on-K3}
In this section, we will summarize necessary backgrounds on two different notions of stability on projective K3 surfaces. At first, we recall the Mukai lattice of a K3 surface. Then we review twisted Gieseker stability of sheaves and their moduli spaces. Finally, we will recollect facts about geometric stability conditions on derived categories of K3 surfaces, including moduli spaces of semistable objects and wall-crossing.

\subsection{The Mukai lattices}
Let $S$ be a smooth projective K3 surface, then there is a polarized weight-two Hodge structure on the integer cohomology ring
$$\tilde{H}(S,\bz)=H^0(S,\bz)\oplus H^2(S,\bz)\oplus H^4(S,\bz)$$
defined by $\tilde{H}^{1,1}(S):=H^{0,0}(S)\oplus H^{1,1}(S)\oplus H^{2,2}(S),\,\tilde{H}^{0,2}(S):=H^{0,2}(S)$ and the Mukai pairing
$$\langle(r_1,\alpha_1,s_1),(r_2,\alpha_2,s_2)\rangle_{\tilde{H}}=\alpha_1.\alpha_2-r_1.s_2-r_2.s_1$$
where the products are taken for cohomological classes. 

The algebraic part of the cohomology $\tilde{H}(S,\bz)$ is precisely the numerical Grothendieck group
$$N(S)=H^0(S,\bz)\oplus \NS(S)\oplus H^4(S,\bz)$$
which inherits the Mukai pairing from $\tilde{H}(S,\bz)$. The Mukai vector 
$$v(E)=\ch(E).\sqrt{\td(S)}=(\rank(E),c_1(E),\rank(E)+c_1(E)^2/2-c_2(E))$$
of a complex $E\in\cate{D}^b(S)$ belongs to the group $N(S)$.

Given a vector $v\in N(S)$, we denote its \emph{orthogonal complement} by 
$$v^{\perp}:=\{w\in \tilde{H}(S,\bz)\,|\,\langle v,w\rangle_{\tilde{H}}=0\}$$
and call $v$ \emph{primitive} if it is not divisible in $N(S)$.

\subsection{Twisted Gieseker stability}
Let $\omega,\beta\in \NS(S)_{\bq}$ with $\omega$ ample, the \emph{$(\omega,\beta)$-twisted Hilbert polynomial} is defined by
$$p_{\omega}^{\beta}(\shf,m):=\int_S(1,-\beta,\frac{\beta^2}{2}).(1,m\omega,\frac{m^2\omega^2}{2}).\ch(\shf).\td(S)$$
for a coherent sheaf $\shf$ on $S$. This gives rise to the notion of $(\omega,\beta)$-twisted Gieseker stability for coherent sheaves, studied first in \cite{MW97}. The following observation asserts that it is the same as Gieseker stability for the cases we concern.

\begin{proposition}\label{K3-variation-of-stability}
Consider a K3 surface $S$, an ample divisor $H$ on $S$, and two algebraic classes $\omega=aH,\beta=bH$ for some $a,b\in\bq$ with $a>0$, then a pure coherent sheaf on $S$ is $H$-Gieseker semistable (resp.\ stable) if and only if it is $(\omega,\beta)$-twisted Gieseker semistable (resp.~stable).

\begin{proof}
One notices that $p_{\omega}^{\beta}(\shf,m)=p_H(\shf,am-b)$ for any coherent sheaf $\shf$ on $S$. 
\end{proof}
\end{proposition}

Now fix a class $v\in N(S)$, then there is a moduli stack $\mathfrak{M}_{\omega}^{\beta}(v)$ of flat families of $(\omega,\beta)$-Gieseker semistable sheaves with Mukai vector $v$. This moduli stack admits a coarse moduli space $M_{\omega}^{\beta}(v)$, which is a projective variety and whose closed points are $S$-equivalence classes of semistable sheaves. The substack $\mathfrak{M}_{\omega}^{\beta,\stable}(v)\subseteq\mathfrak{M}_{\omega}^{\beta}(v)$ parameterizing stable sheaves is open, and is a $\mathbb{G}_m$-gerbe over the open locus $M_{\omega}^{\beta,\stable}(v)\subseteq M_{\omega}^{\beta}(v)$ consisting of $S$-equivalence classes of stable sheaves with class $v$. The superscript $\beta$ in the notation will be omitted when $\beta=0$. In this case, one reduces to the usual Gieseker stability. Nevertheless, Proposition~\ref{K3-variation-of-stability} implies

\begin{corollary}
Consider a K3 surface $S$, an ample divisor $H$ on $S$, and two algebraic classes $\omega=aH,\beta=bH$ for some $a,b\in\bq$ with $a>0$, then $M_{\omega}^{\beta}(v)=M_H(v)$ for any $v\in N(S)$.
\end{corollary}

So one can actually reduce to the usual Gieseker stability in many cases. A nice exposition of basic properties of Gieseker stability is \cite{HL}.

\subsection{Moduli spaces of Gieseker semistable sheaves}
In this article, we will involve some moduli space $M_H(mv_0)$, where $H$ is a primitive ample divisor, $v_0$ is primitive in $N(S)$, and $m=1,2$. Such a moduli space is fairly understood when $H$ is general enough, based on contributions of Huybrechts, Mukai, O'Grady, Yoshioka among many others. 

\begin{theorem}[{{\cite[Proposition 4.12]{Yo01}}}]\label{K3-moduli-sheaf-description-primitive}
Consider a primitive vector $v=(r,\Delta,s)$ with $r>0$ and a general primitive ample class $H$ with respect to $v$ on a K3 surface $S$, then $M_H(v)=M^{\stable}_H(v)$ is a smooth irreducible symplectic manifold of dimension $2n=:\langle v,v\rangle_{\tilde{H}}+2$, deformation equivalent to the Hilbert scheme of $n$ points on $S$.
\end{theorem}

\begin{theorem}[{{\cite[Theorem 1.6]{PR13}}}]\label{K3-moduli-sheaf-description-OG10}
Consider a primitive vector $v=(r,\Delta,s)$ with $r>0$ and $\langle v,v\rangle_{\tilde{H}}=2$, and a general primitive ample class $H$ on a K3 surface $S$ with respect to $2v$, then $M=M_H(2v)$ admits a symplectic resolution $\rho\colon \tilde{M}\rightarrow M$ such that $\tilde{M}$ is a smooth irreducible symplectic manifold of dimension $10$ deformation equivalent to the O'Grady tenfold \cite{O'Grady99}.
\end{theorem}

The notion of being a \emph{general} ample class with respect to $v$ in the theorems above is subtle to define, and we refer the readers to \cite{PR23} for some discussions. On the other hand, moduli spaces of semistable sheaves considered in this article will be on a generic K3 surface $(S,H)$ so that one does not need to worry about generality of $H$, except for the following one.

\begin{example}
Let $v=(1,0,1-n)$ for some $n>0$, then $M_H(v)=M^{\stable}_H(v)$ is isomorphic to the Hilbert scheme $S^{[n]}$ of $n$ points on $S$ for any polarized K3 surface $(S,H)$ and a point of $M_H(v)$ is isomorphic to the ideal sheaf $\shi_Z$ of a length $n$ subscheme $Z\subset S$.
\end{example}

The second integral cohomology of these irreducible symplectic manifolds, from the aspect of moduli spaces of sheaves, is studied thoroughly in \cite{O'Grady97,PR13,PR24,Rap08,Yo01}.

\begin{definition}\label{K3-Mukai-homomorphism}
Let $v\in N(S)$ be a primitive vector in $N(S)$ with $r>0$ and $H$ be a primitive ample divisor on $S$ such that $M=M_H(v)$ is a smooth irreducible symplectic manifold, then we define the \defi{Mukai homomorphism} $\theta_v\colon v^{\perp}\rightarrow H^2(M,\bz)$ by
$$\theta_v(w)=\frac{1}{\rho}[\Phi_{\shu}^{\textup{coho}}(w)]_2$$
where $\Phi^{\textup{coho}}_{\shu}$ is the cohomological Fourier--Mukai transform whose kernel $\shu$ is a quasi-universal family for $M$ with similitude $\rho$, and the operator $[-]_2$ takes components in $H^2(M,\bz)$.
\end{definition}

\begin{proposition}[Yoshioka \cite{Yo01}]\label{K3-Mukai-isomorphism}
Consider a primitive Mukai vector $v=(r,\Delta,s)$ in $N(S)$ with $r>0$ and $\langle v,v\rangle_{\tilde{H}}>0$ and a primitive ample class $H$ such that $M=M_H(v)$ is an irreducible symplectic manifold, then the Mukai homomorphism $\theta_v$ is an isometry of polarized weight-two Hodge structures, where $H^2(M,\bz)$ is endowed with the Beauville-Bogomolov-Fujiki form.
\end{proposition}

In some non-primitive cases, a morphism $\theta_v\colon v^{\perp}\rightarrow H^2(M_H(v),\bz)$ with similar properties is constructed in \cite[Proposition 3.9]{PR13} and \cite[Section 4]{PR24}. We will discuss this later.

\subsection{Geometric stability conditions}
An ample class $\omega\in \NS(S)_{\br}$ defines a notion of slope stability called $\mu_{\omega}$-stability by declaring the slope
$$\mu_{\omega}(\shf):=\frac{c_1(\shf)\cdot\omega}{\rank(\shf)}$$
for a torsion-free sheaf $\shf$ on $S$. By truncating Harder--Narasimhan filtrations at $\beta\cdot\omega$, a pair $\omega,\beta\in\NS(S)_{\br}$ with $\omega$ ample determines a torsion pair $(\scht_{\omega,\beta},\schf_{\omega,\beta})$ on $\cate{Coh}(S)$ such that
\begin{align*}
\scht_{\omega,\beta}&=\{\shf\in\cate{Coh}(S)\,|\,\textup{torsion free part of }\shf\textup{ has }\mu_{\omega}\textup{-semistable HN factors of slope} >\beta\cdot\omega\}\\
\schf_{\omega,\beta}&=\{\shf\in\cate{Coh}(S)\,|\,\shf \textup{ is torsion free and its }\mu_{\omega}\textup{-semistable HN factors have slope} \leq\beta\cdot\omega\}
\end{align*}
Here and henceforth, we adopt the convention that $\mu_{\omega}=+\infty$ if the torsion-free part is zero.

The torsion pair gives a bounded $t$-structure on $\cate{D}^b(S)$ with heart 
$$\scha_{\omega,\beta}=\{E\in \cate{D}^b(S)\,|\, \shh^i(E)=0\textup{ for }i\notin\{-1,0\},\shh^{-1}(E)\in\schf_{\omega,\beta},\shh^0(E)\in\scht_{\omega,\beta}\}$$
and a group homomorphism
$$Z_{\omega,\beta}\colon N(S)\rightarrow\bc,\quad (r,\Delta,s)\mapsto \frac{1}{2}(2\beta\cdot\Delta-2s+r(\omega^2-\beta^2))+(\Delta-r\beta)\cdot\omega\sqrt{-1}$$
Suppose that $Z(\she)\notin\br_{\leq0}$ for any spherical sheaf $\she$ on $S$, the pair $\sigma_{\omega,\beta}=(\scha_{\omega,\beta},Z_{\omega,\beta})$ is a stability condition according to \cite[Lemma 6.2]{Bri08}. Moreover, the skyscraper sheaf $\sho_p$ is always stable of phase $1$ with respect to $\sigma_{\omega,\beta}$ for every $p\in S$. 

\begin{definition}
A stability condition $\sigma$ on $S$ is called \defi{geometric} if the skyscraper sheaf $\sho_p$ is a $\sigma$-stable object of the same phase for every point $p\in S$.
\end{definition}

All the geometric stability conditions on $\cate{D}^b(S)$ (with support property) belong to the same connected component in $\Stab(S)$, say $\Stab^{\circ}(S)$, according to \cite[Theorem 5.36]{Dell25}. Moreover, up to a unique $\grp$ action on $\Stab(S)$, a geometric stability condition on $\cate{D}^b(S)$ goes to a uniquely determined stability condition $\sigma_{\omega,\beta}$ by virtue of \cite[Proposition 10.3]{Bri08}. 

\begin{theorem}[{{\cite[Proposition 14.2]{Bri08}}} and {{\cite[Section 6]{Yo01}}}]\label{K3-large-volume-limits}
Consider a pair $(\omega,\beta)\in \NS(S)_{\bq}^2$ such that $\omega$ is ample and $\sigma_{\omega,\beta}$ is a stability condition, and a Mukai vector $v=(r,\Delta,s)\in N(S)$ with $r>0,(\Delta-r\beta)\cdot\omega>0$ or with $r=0$, then there exists a positive number $t_0$ such that, for any $t\geq t_0$, an object $E$ of $\cate{D}^b(S)$ with vector $v$ is $\sigma_{t\omega,\beta}$-semistable if and only if $E$ is isomorphic to a shift of an $(\omega,\beta)$-twisted semistable sheaf on $S$. Moreover, $M_{\sigma_{t\omega,\beta}}(v)=M_{\omega}^{\beta}(v)$ for $t\gg0$.
\end{theorem}

\subsection{Chambers, walls and wall-crossing}
The walls and wall-crossing behaviors of stability conditions in $\Stab^{\circ}(S)$ for a K3 surface $S$ are studied intensively in \cite{BM14,BM14-MMP,MZ16}. 

Suppose that $W\subset\Stab(S)$ is a wall with respect to a given $v\in N(S)$ with $\langle v,v\rangle_{\tilde{H}}>0$ and $\sigma_0\in W$ is not in any other wall, then $W$ is called a \emph{totally semistable wall} if $M^{\stable}_{\sigma_0}(v)=\varnothing$.

Suppose in addition that $\sigma_{\pm}$ are $v$-generic stability conditions in two different sides of $W$ and sufficiently close to $W$ in $\Stab(S)$. Then \cite[Proposition 5.2]{MZ16} asserts a birational morphism $$\pi_{\pm}\colon M_{\sigma_{\pm}}(v)\rightarrow\overline{M}_{\pm}$$ 
which contracts curves in $M_{\sigma_{\pm}}(v)$ parameterizing $S$-equivalent objects under $\sigma_0$, where $\overline{M}_{\pm}$ is the image of $\pi_{\pm}$ in $M_{\sigma_0}(v)$. Then $W$ is called a \emph{fake wall} if $\pi_{\pm}$ are isomorphisms, a \emph{flopping wall} if $\pi_{\pm}$ are small contractions, and a \emph{divisorial wall} if $\pi_{\pm}$ are divisorial contractions.

The criterion to decide a wall is developed in \cite{BM14-MMP} for primitive vectors and in \cite{MZ16} for special multiples of primitive vectors. Here we only record a part of the criterion for the primitive case as all but one of the explicit examples we will consider here are of this sort.

Suppose that $v$ is a primitive Mukai vector and $W$ is a wall with respect to $v$. Then one defines a rank two hyperbolic lattice $\mathbb{H}_W\subset N(S)$ associated with the wall $W$ to be the set containing all the $w$ with $\Im(Z(w)/Z(v))=0$ for any $\sigma=(\scha,Z)\in W$.

\begin{theorem}[{{\cite[Theorem 5.7]{BM14-MMP}}}]
Consider a primitive $v$ with $2\leq\langle v,v\rangle_{\tilde{H}}\leq 4$, then $W$ is a divisorial wall if and only if there exists a spherical vector $w\in\mathbb{H}_W$ with $\langle w,v\rangle_{\tilde{H}}=0$ or an isotropic vector $w\in\mathbb{H}_W$ with $\langle w,v\rangle_{\tilde{H}}$ being $1$ or $2$; otherwise, $W$ is flopping if there exists a spherical vector $w\in\mathbb{H}_W$ with $0<\langle w,v\rangle_{\tilde{H}}\leq\frac{1}{2}\langle v,v\rangle_{\tilde{H}}$.

\begin{proof}
The only difference from the original theorem is that the condition $2\leq\langle v,v\rangle_{\tilde{H}}\leq 4$ ensures that $v$ cannot be written as the sum of two positive classes. One can see this point from the first paragraph of the argument for \cite[Lemma 6.13]{MZ16}.
\end{proof}
\end{theorem}

\subsection{Moduli spaces of semistable objects}
The moduli stacks $\mathfrak{M}_{\sigma_{\omega,\beta}}(v)$ and $\mathfrak{M}^{\stable}_{\sigma_{\omega,\beta}}(v)$ admit projective moduli spaces according to \cite{BM14,To08}. Similar to Gieseker stability, one has

\begin{theorem}[{{\cite[Corollary 6.9]{BM14}}}]
Consider a primitive vector $v$ and a $v$-generic geometric stability condition $\sigma$ on $\cate{D}^b(S)$, then $M_{\sigma}(v)=M^{\stable}_{\sigma}(v)$ is a smooth irreducible symplectic manifold of dimension $\langle v,v\rangle_{\tilde{H}}+2$, deformation equivalent to the Hilbert scheme of points on $S$.
\end{theorem}

One notices that the construction of Mukai homomorphism extends to moduli spaces of $\sigma$-semistable objects without any obstruction.

\begin{proposition}[{{\cite[Theorem 6.10]{BM14}}}]
Consider a primitive vector $v$ with $\langle v,v\rangle_{\tilde{H}}>0$, and a $v$-generic geometric stability condition $\sigma$ such that $M=M_{\sigma}(v)$ is an irreducible symplectic manifold admitting a quasi-universal family, then the Mukai homomorphism $\theta_{\sigma,v}$ defined in Definition \ref{K3-Mukai-homomorphism} is an isometry of polarized weight-two Hodge structures, where the Hodge structure on $H^2(M,\bz)$ is given by the Beauville-Bogomolov-Fujiki form.
\end{proposition}

Moreover, the Mukai isomorphism is compatible with cohomology. One recalls that any autoequivalence of $\cate{D}^b(S)$ is isomorphic to a Fourier--Mukai transform $\Phi_K$ and the associated cohomological Fourier--Mukai map $\Phi_K^{\textup{coho}}$ is a Hodge self-isometry on $\tilde{H}(S,\bz)$ (see e.g.\ \cite{Hu}).

\begin{proposition}[{{\cite[Proposition 3.5]{Ou18}}}]\label{K3-invariant-lattices-compatible}
Consider an autoequivalence $\Phi_K$ of $\cate{D}^b(S)$ and a primitive vector $v\in N(S)$ with $\langle v,v\rangle_{\tilde{H}}>0$. Let $\sigma$ be a $v$-generic geometric stability condition fixed by the isometry $\Phi^{\textup{coho}}_K$. Suppose that $\Phi_K$ induces an automorphism $\phi$ of $M_{\sigma}(v)$, then
\begin{displaymath}
	\xymatrix{
		v^{\perp}\ar[d]_{\theta_{\sigma,v}}\ar[rr]^{\Phi_K^{\textup{coho}}}&&v^{\perp}\ar[d]_{\theta_{\sigma,v}}\\
		H^2(M_{\sigma}(v),\bz)\ar[rr]^{\phi_*}&&H^2(M_{\sigma}(v),\bz)
	}
\end{displaymath}
is a commutative diagram, where $\phi_*$ is the Hodge isometry induced by $\phi$. 
\end{proposition}

The arguments of the two propositions above also run for the stable locus $M^{\stable}_{\sigma}(v)$ once it has a quasi-universal family and one can find $H^2(M^{\stable}_{\sigma}(v),\bz)\cong H^2(M,\bz)$ for an irreducible symplectic manifold $M$. Besides, Proposition~\ref{K3-invariant-lattices-compatible} holds in some non-primitive cases for a suitably constructed Mukai homomorphism, such as \cite[Lemma 3.11]{PR13} and \cite[Section 4]{PR24}.

\section{Kuznetsov components of quartic double solids}\label{ST-Family-QDS}
In this section, we apply this framework to the Enriques categories over K3 surfaces described in Example~\ref{En_example-1}. In the first three subsections, we work out the (essentially unique) stability condition $\sigma_{\omega_q,\beta_q}$ on $\cate{D}^b(S)$ invariant under the $\bz/2\bz$-action, write down the induced involutions on moduli spaces with respect to $\sigma_{\omega_q,\beta_q}$ in Proposition~\ref{ME-involution}, and write down the explicit covering maps between moduli spaces as in Proposition~\ref{ME-double-cover}. In the remaining subsections, we describe the involutions and covering maps for some explicit numerical classes.

\subsection{Generalities on these Enriques categories}
We start with a K3 surface $S$ polarized by a very ample divisor class $H$ with $H^2=4$. In other words, $S$ is a quartic K3 surface. The double covering $X\rightarrow\bp^3$ branched over $S$ admits a semiorthogonal decomposition
$$\cate{D}^b(X)\cong\langle\Kuc(X),\sho_X,\sho_X(1)\rangle $$
such that $\Kuc(X)$ is an Enriques category over $S$. 

The numerical Grothendieck group $N(\Kuc(X))$ has rank 2, with a basis $\mu_1$, $\mu_2$ such that
$$\ch(\mu_1)=1-\frac{1}{2}h^2,\quad \ch(\mu_2)=h-\frac{1}{2}h^2-\frac{1}{3}h^3$$
and the Euler form is given in this basis by
$$\begin{pmatrix}
	\chi(\mu_1,\mu_1)&\chi(\mu_1,\mu_2)\\\chi(\mu_2,\mu_1)&\chi(\mu_2,\mu_2)
\end{pmatrix}=\begin{pmatrix}
	-1&-1\\-1&-2
\end{pmatrix}$$
Moreover, one notices that the induced action of $\funct{U}$ on $N(\Kuc(X))$ is trivial as the covering involution $X\rightarrow\bp^3$ fixes the polarization class $h$. So $\funct{U}$ induces an involution on moduli spaces for any numerical class in $N(\Kuc(X))$ and any $\funct{U}$-fixed stability condition.

The residual action on $\cate{D}^b(S)\cong \Kuc(X)^{\funct{U}}$ is generated by the autoequivalence
$$\Pi_q=\funct{T}_{\sho_S}\circ(-\otimes\sho_S(H))\circ\funct{T}_{\sho_S}\circ(-\otimes\sho_S(H))[-1]$$
and a direct computation shows that:

\begin{proposition}\label{QDS-induced-Hodge-isometry}
The Hodge isometry $\tau_q$ induced by $\Pi_q$ is given by the mapping
$$\tau_q\colon (r,\Delta,s)\mapsto(-r-2s-\Delta\cdot H,(r+s+\Delta\cdot H)H-\Delta,-2r-s-\Delta\cdot H)$$
\end{proposition}

Moreover, one checks that

\begin{proposition}\label{QDS-fixed-stability-condition}
The pair $\omega_q=\frac{1}{2}H$ and $\beta_q=-\frac{1}{2}H$ is the only pair $(\omega,\beta)$ making
$$Z_{\omega,\beta}(r,\Delta,s)=Z_{\omega,\beta}(\tau_q(r,\Delta,s))$$
hold for any $(r,\Delta,s)\in N(S)$.
	
\begin{proof}
Comparing the imaginary parts, one can obtain an equation
$$2\Delta.\omega_q=(r+s+\Delta\cdot H)(2\omega_q.\beta_q+\omega_q.H)-(\Delta\cdot H)(\omega_q.\beta_q)$$
which directly implies $\omega_q.H=-2\omega_q.\beta_q>0$ and the equation
$$2\Delta.\omega_q=-(\Delta\cdot H)(\omega_q.\beta_q)\qquad(\star)$$
for any $\Delta\in\NS(S)$. The latter equation then gives $\beta_q.H=-2$ by plugging $\Delta=\beta_q$.
		
Comparing the real parts, one can obtain the equation
$$	4\Delta.\beta_q=(2r+2s+\Delta\cdot H)(\beta_q^2-\omega_q^2+2+\beta_q.H)+(\Delta\cdot H)(\beta_q.H)$$
which follows $\omega_q^2=\beta_q^2$ by $\beta_q\cdot H=-2$ and then gives the equation
$$4\Delta.\beta_q=(\Delta\cdot H)(\beta_q.H)=-2(\Delta\cdot H)$$ 
for any $\Delta\in\NS(S)$. So one must have $2\beta_q=-H$. Then $\omega_q^2=\beta_q^2=1$ and the equation $(\star)$ ensures $\omega_q.\beta_q=-1$. So the equation $(\star)$ becomes $2\Delta.\omega_q=\Delta\cdot H$ and we are done.
\end{proof}
\end{proposition}

Together with Theorem~\ref{En-equivariant-stability} and the comments before Proposition~\ref{En-back-to-stability}, one obtains

\begin{corollary}\label{QDS-stability-condition-exists}
Suppose that there is a $\funct{U}$-fixed stability condition $\sigma$ on $\Kuc(X)$ such that $\sigma^{\funct{U}}$ is a geometric stability condition, then $\sigma^{\funct{U}}=\sigma_{\omega_q,\beta_q}.\tilde{g}$ for a uniquely determined $\tilde{g}\in\grp$. In particular, $\sigma_{\omega_q,\beta_q}$ is indeed a stability condition and the involution $\Pi_q$ fixes $\sigma_{\omega_q,\beta_q}$.
\end{corollary}

The stability conditions on $\Kuc(X)$ constructed in \cite[Theorem 6.8]{BLMS23} satisfy the conditions in Corollary~\ref{QDS-stability-condition-exists} according to \cite[Lemma 6.2 and Lemma 6.3]{PPZ26}. In this case, one can use this corollary, Proposition~\ref{En-back-to-stability} and Proposition~\ref{En-invariant-grp-action} to deduce that all $\funct{U}$-fixed stability conditions on $\Kuc(X)$ are in the same $\grp$-orbit. One notices that $(\sigma.\tilde{g}).\tilde{g}=\sigma$ if and only if one has $\tilde{g}=(I,\id)$, so a stability condition on $\Kuc(X)$ is \emph{Serre-invariant} if and only if it is $\funct{U}$-fixed. As a consequence, we recover a special case of \cite[Theorem A.10]{JLLZ24}.

\subsection{The induced involutions}
It is computed in \cite[Lemma 7.4]{BP23} that the invariant lattice of the Hodge isometry $\tau_q\colon \tilde{H}(S,\bz)\rightarrow\tilde{H}(S,\bz)$ induced by $\Pi_q$ is spanned by the Mukai vectors $u_1=(1,0,-1)$ and $u_2=(1,-H,1)$. So one applies Proposition~\ref{ME-involution} to see

\begin{proposition}
Consider a quartic K3 surface $(S,H)$, then $\Pi_q$ induces an involution
$$\pi_q\colon M_{\sigma_{\omega_q,\beta_q}}(au_1+bu_2)\rightarrow M_{\sigma_{\omega_q,\beta_q}}(au_1+bu_2)$$
for each pair $(a,b)$ of integers.
\end{proposition}

The following observation on the autoequivalence $\funct{O}_q:=\funct{T}_{\sho_S}\circ(-\otimes\sho_S(H))$ of $\cate{D}^b(S)$ is used in the subsequent discussions.

\begin{proposition}\label{QDS-rotation-identification}
One has $\funct{O}_q.\sigma_{\omega_q,\beta_q}=\sigma_{\omega_q,\beta_q}.\tilde{g}$ for a unique $\tilde{g}\in\grp$. In particular, the autoequivalence $\funct{O}_q$ induces $M_{\sigma_{\omega_q,\beta_q}}(bu_1-au_2)\cong M_{\sigma_{\omega_q,\beta_q}}(au_1+bu_2)$ for any integers $a$ and $b$.
	
\begin{proof}
Set $\funct{O}_q.\sigma_{\omega_q,\beta_q}=(Z',\scha')$, one computes that
$$Z_{\omega_q,\beta_q}(r,\Delta,s)=-(s+\frac{1}{2}\Delta\cdot H)+(r+\frac{1}{2}\Delta\cdot H)\sqrt{-1}$$
$$Z'=Z_{\omega_q,\beta_q}(\tau(r,\Delta,s))=-(r+\frac{1}{2}\Delta\cdot H)-(s+\frac{1}{2}\Delta\cdot H)\sqrt{-1}$$
where $\tau$ is the Hodge isometry induced by $\funct{O}_q$.

Since $(Z',\scha')$ is also a geometric stability condition, one can deduce $$\funct{O}_q.\sigma_{\omega_q,\beta_q}=\sigma_{\omega_q,\beta_q}.\tilde{g}$$ 
by inspecting the argument for \cite[Proposition 10.3]{Bri08} and comparing $Z'$ with $Z_{\omega_q,\beta_q}$. Then the desired isomorphism follows from $\tau(au_1+bu_2)=bu_1-au_2$.
\end{proof}
\end{proposition}

Obviously, the induced involution $\pi_q$ on the moduli space $M_{\sigma_{\omega_q,\beta_q}}(au_1+bu_2)$ is the same as the induced involution on $M_{\sigma_{\omega_q,\beta_q}}(bu_1-au_2)$. 

\subsection{The covering morphisms in two directions}\label{QDS-covering-morphism}
To explicitly write down the covering maps given in Proposition~\ref{ME-double-cover}, we need to figure out the images $\forg(v)$ and $\sfinf(\mu)$.

\begin{proposition}\label{QDS-image-forg}
Consider the homomorphism $\forg\colon N(S)\cong N(\Kuc(X)^{\funct{U}})\rightarrow N(\Kuc(X))$, then
$$\forg(r,\Delta,s)=(2s+\Delta\cdot H)\mu_1-(r+s+\Delta\cdot H)\mu_2$$
for any Mukai vector $(r,\Delta,s)\in N(S)$. 

\begin{proof}
Set $\forg(r,\Delta,s)=a\mu_1+b\mu_2$. One computes that
$$\sfinf(\forg(r,\Delta,s))=(r,\Delta,s)+\tau_q(r,\Delta,s)=(-2s-\Delta\cdot H,(r+s+\Delta\cdot H)H,-2r-\Delta\cdot H)$$ 
so \cite[Lemma 6.4]{PPZ26} and the adjunction between the forgetful functor and inflation functor give
\begin{align*}
	2b=\chi(\forg(1,0,1),\forg(r,\Delta,s))=-\langle(1,0,1),\sfinf(\forg(r,\Delta,s))\rangle_{\tilde{H}}&=-2r-2s-2\Delta\cdot H\\
	-a=\chi(\forg(0,0,1),\forg(r,\Delta,s))=-\langle(0,0,1),\sfinf(\forg(r,\Delta,s))\rangle_{\tilde{H}}&=-2s-\Delta\cdot H
\end{align*}
with which one computes that $a=2s+\Delta\cdot H$ and $b=-r-s-\Delta\cdot H$.
\end{proof}
\end{proposition}

Since $\funct{U}$ acts trivially on $N(\Kuc(X))$, the composite $\textup{\forg}\circ \textup{\sfinf}\colon N(\Kuc(X))\rightarrow N(\Kuc(X))$ is the multiplication by $2$. One can therefore deduce that

\begin{corollary}\label{QDS-image-inflation}
Consider the homomorphism $\sfinf\colon N(\Kuc(X)) \rightarrow N(\Kuc(X)^{\funct{U}})\cong N(S)$, then
$$\sfinf(a\mu_1+b\mu_2)=-(a+b)u_1+bu_2$$
for any $a,b\in\bz$. 

\begin{proof}
Suppose that $\sfinf(\mu_1)=(r,\Delta,s)$, then $\forg(r,\Delta,s)=2\mu_1$. So one has
$$r+s+\Delta\cdot H=0,\quad 2s+\Delta\cdot H=2,\quad 2rs-\Delta^2=-2$$
by the previous proposition and $-\langle\sfinf(\mu_1),\sfinf(\mu_1)\rangle_{\tilde{H}}=2\chi(\mu_1,\mu_1)=-2$. It follows that
$$(\Delta\cdot H)^2=2\Delta^2$$
and one thus obtains $\Delta=0$ by Hodge index theorem. So $\sfinf(\mu_1)=(-1,0,1)=-u_1$. Then, by a suitable substitution of variables, one can check $\sfinf(\mu_2)=(0,-H,2)=u_2-u_1$.
\end{proof}
\end{corollary}

As all the $\funct{U}$-fixed stability conditions are in the same $\grp$-orbit and their induced equivariant stability conditions are also in the same $\grp$-orbit, one has

\begin{proposition}\label{QDS-double-cover}
Consider the quartic double solid $X$ over a quartic K3 surface $S$ and a $\funct{U}$-fixed stability condition $\sigma$ on $\Kuc(X)$, then one has a surjection
$$M_{\sigma}(-(a+b)\mu_1+b\mu_2)\rightarrow M_{\sigma_{\omega_q,\beta_q}}(au_1+bu_2)^{\Pi_q}$$
for each pair of $a,b\in\bz$. This surjection has degree $2$ at points corresponding to stable objects and is étale at points corresponding to stable objects which are not fixed by $\funct{U}$.
\end{proposition}

\begin{remark}
The autoequivalence $\funct{L}_{\sho_X}(-\otimes\sho_X(1))$ of $\Kuc(X)$ acts on the lattice $N(\Kuc(X))$ by sending $\mu_1$ to $\mu_2-\mu_1$ and $\mu_2$ to $\mu_2-2\mu_1$, which is compatible with the action of $\funct{T}_{\sho_S}(-\otimes\sho_S(H))$ on $N(S)$ through inflation. In fact, the double cover
$$M_{\sigma}(-(b-a)\mu_1-a\mu_2)\rightarrow M_{\sigma_{\omega_q,\beta_q}}(bu_1-au_2)^{\Pi_q}$$
is isomorphic to $M_{\sigma}(-(a+b)\mu_1+b\mu_2)\rightarrow M_{\sigma_{\omega_q,\beta_q}}(au_1+bu_2)^{\Pi_q}$ for each pair of $a,b\in\bz$.
\end{remark}

\subsection{Over the class $\mu_1$}
Let us start with the simplest example of covering maps
$$\bigsqcup_{\forg(r,\Delta,s)=-\mu_1}M_{\sigma_{\omega_q,\beta_q}}(r,\Delta,s)\rightarrow M_{\sigma}(-\mu_1)^{\funct{U}}$$
where $\sigma$ is a $\funct{U}$-fixed stability condition on $\Kuc(X)$ here and henceforth in this section.

\begin{proposition}\label{QDS-(1,-L,0)}
The object $\sho_L(-2)$ is $\sigma_{\omega_q,\beta_q}$-stable for any line $L\subset S$.

\begin{proof}
At first, one notices that $\sho_L(-2)$ is in the heart $\scha_{\omega_q,\beta_q}$. Suppose otherwise that $\sho_L(-2)$ is not $\sigma_{\omega_q,\beta_q}$-stable, then one can choose a short exact sequence
$$0\rightarrow A\rightarrow\sho_L(-2)\rightarrow B\rightarrow0$$
in $\scha_{\omega_q,\beta_q}$ such that $A$ is $\sigma_{\omega_q,\beta_q}$-stable with $\mu(A)\geq\mu(B)$. One has
$$\Im(Z_{\omega_q,\beta_q}(A))=(r+\frac{1}{2}\Delta\cdot H)\geq0\quad\textup{and}\quad \Im(Z_{\omega_q,\beta_q}(B))=(\frac{1}{2}-r-\frac{1}{2}\Delta\cdot H)\geq0$$
with $v(A)=(r,\Delta,s)$. Together with $\mu(A)\geq\mu(B)$, one gets $2r+\Delta\cdot H=0$. 

Suppose that $2r+\Delta\cdot H=0$, then $\phi(A)=1$. According to \cite[Lemma 10.1]{Bri08}, the stable object $A$ is either isomorphic to $\sho_x$ for some point $x\in S$ or isomorphic to $\she[1]$ for some locally free sheaf $\she$. However, one can obtain from the long exact sequence
$$0\rightarrow\shh^{-1}(A)\rightarrow0\rightarrow\shh^{-1}(B)\rightarrow\shh^0(A)\rightarrow\sho_L(-2)\rightarrow\shh^0(B)\rightarrow0$$
that both cases are impossible. It means that $\sho_L(-2)$ is $\sigma_{\omega_q,\beta_q}$-stable. 
\end{proof}
\end{proposition}

On the other hand, the right hand side fixed locus is known.

\begin{proposition}[\cite{FP23,PY22}]\label{QDS-class-mu1}
Consider a Serre invariant stability condition $\sigma$ on $\Kuc(X)$, then any $\sigma$-semistable object with class $-\mu_1$ is stable and isomorphic to the object $\shi_L$ for some line $L\subset X$ up to a shift. In particular, $M_{\sigma}(-\mu_1)$ is isomorphic to the Fano variety $\Sigma(X)$ of lines.
\end{proposition}

\begin{corollary}\label{QDS-mu1-covering}
The covering map for the class $-\mu_1$ is exactly
$$\bigsqcup_{L\subset S\textup{ a line}} M_{\sigma_{\omega_q,\beta_q}}(1,-L,0)\sqcup M_{\sigma_{\omega_q,\beta_q}}(0,L,-1)\rightarrow M_{\sigma}(-\mu_1)^{\funct{U}}$$
Moreover, one has $M_{\sigma_{\omega_q,\beta_q}}(0,L,-1)=\{\sho_L(-2)\}$ and $M_{\sigma_{\omega_q,\beta_q}}(1,-L,0)=\{\sho_S(-L)\}$.

\begin{proof}
Since $\funct{U}$ is nothing but pushforward along the covering involution for $X\rightarrow\bp^3$, the fixed locus $M_{\sigma}(-\mu_1)^{\funct{U}}$ is a finite set consisting of all $\shi_L[1]$ with $L\subset S$. As objects in $M_{\sigma}(-\mu_1)^{\funct{U}}$ are stable, one sees by Corollary~\ref{ME-double-cover-stable} that the number of points in the source is exactly two times of those in $M_{\sigma}(-\mu_1)^{\funct{U}}$. Hence the first assertion follows from Proposition~\ref{QDS-image-forg}, Proposition~\ref{QDS-(1,-L,0)}, and counting points. The second one is due to the observation $\Pi_q\sho_L(-2)\cong\sho_S(-L)$.
\end{proof}
\end{corollary}

The forgetful functor $\funct{Forg}\colon\cate{D}^b(S)\rightarrow\Kuc(X)$ is given by $\funct{L}_{\sho_X}\funct{L}_{\sho_X(1)}(j_*(-)\otimes\sho_X(2))$ for the ramification embedding $j\colon S\hookrightarrow X$ due to \cite[(6.4)]{PPZ26}, so $\funct{Forg}(\sho_L(-2))\cong\shi_L[1]$. In fact, one computes that $\funct{L}_{\sho_X}\funct{L}_{\sho_X(1)}\sho_L\cong\funct{L}_{\sho_X}\sho_L\cong\shi_L[1]$ for any line $L\subset X$.

\subsection{Over the vector $u_1$}\label{QDS-ss-over-u1}
Here we are going to figure out the covering morphism
$$M_{\sigma}(-\mu_1)\rightarrow M_{\sigma_{\omega_q,\beta_q}}(u_1)^{\Pi_q}$$
for any quartic K3 surface. The first thing is to understand the moduli space $M_{\sigma_{\omega_q,\beta_q}}(u_1)$.

\begin{proposition}\label{QDS-u1-strict-semi}
A strictly $\sigma_{\omega_q,\beta_q}$-semistable object $E$ in $\scha_{\omega_q,\beta_q}$ with Mukai vector $u_1$ is isomorphic to an extension of $\sho_S(-L)$ and $\sho_{L}(-2)$ for a line $L\subset S$.

\begin{proof}
One chooses a short exact sequence $0\rightarrow A\rightarrow E\rightarrow B\rightarrow0$ in $\scha_{\omega_q,\beta_q}$ such that $A$ and $B$ are $\sigma_{\omega_q,\beta_q}$-semistable with the same phase. Set $v(A)=(r,\Delta,s)$, then one gets $2r+\Delta\cdot H=1$ and $2s+\Delta\cdot H=-1$ from $\Im(Z_{\omega_q,\beta_q}(A))\geq0,\Im(Z_{\omega_q,\beta_q}(B))\geq0$ and $\mu(A)=\mu(B)$. It means that $A$ is sent to $M_{\sigma}(-\mu_1)$ via $\forg$ due to Proposition~\ref{QDS-image-forg}. One concludes by Corollary~\ref{QDS-mu1-covering}.
\end{proof}
\end{proposition}

\begin{remark}\label{QDS-u1-strict-semi-part}
Let $Z\subset S$ be a length $2$ subscheme such that $Z$ is in a line $L\subset S$, then the short exact sequence $0\rightarrow\sho_S(-L)\rightarrow\shi_Z\rightarrow\sho_L(-2)\rightarrow0$ implies that $\shi_Z$ is strictly $\sigma_{\omega_q,\beta_q}$-semistable.
\end{remark}

It remains to determine the stable objects in the moduli space $M_{\sigma_{\omega_q,\beta_q}}(u_1)$, and the following observation is useful.

\begin{proposition}\label{QDS-u_1-wall-crossing-lemma}
Suppose that $D$ is a class in $\NS(S)$ such that $D^2=-2$ and $D\cdot H=1$, then the stability condition $\sigma_{t\omega_q,\beta_q}$ is $(1,-D,n)$-generic for any $t\geq 1$ and integer $n$. Suppose in addition that $n\geq0$, then the moduli space $M_{\sigma_{t\omega_q,\beta_q}}(1,-D,n)$ is non-empty only when $n=0$.

\begin{proof}
Since the Mukai vector $(1,-D,n)$ is always primitive, it suffices to show that there are no strictly $\sigma_{t\omega_q,\beta_q}$-semistable objects. Otherwise, one has a short exact sequence
$$0\rightarrow A\rightarrow E\rightarrow B\rightarrow0$$
of semistable objects in $\schp(\phi(E))$ such that $v(E)=(1,-D,n)$. Set $v(A)=(r,\Delta,s)$, then
$$\Im(Z_{\omega_q,\beta_q}(A))=(r+\frac{1}{2}\Delta\cdot H)t\geq0\quad\textup{and}\quad \Im(Z_{\omega_q,\beta_q}(B))=(\frac{1}{2}-r-\frac{1}{2}\Delta\cdot H)t\geq0$$
which has no valid solutions as $\mu(A)=\mu(B)$. So the first assertion is proved. The second one follows from $\langle E,E\rangle_{\tilde{H}}=D^2-2rs=-2-2n\leq -2$ for any $E$ in $M_{\sigma_{t\omega_q,\beta_q}}(1,-D,n)$.
\end{proof}
\end{proposition}

Then we can prove the critical statement which, together with Theorem~\ref{K3-large-volume-limits}, implies that the stable objects in $M_{\sigma_{\omega_q,\beta_q}}(u_1)$ are isomorphic to sheaves.

\begin{proposition}\label{QDS-u_1-wall-crossing}
The stability condition $\sigma_{t\omega_q,\beta_q}$ is $u_1$-generic for $t>1$. Moreover, one has $M_{\sigma_{\omega_q,\beta_q}}(u_1)\cong M_H(u_1)$ if and only if $S$ does not contain any line.

\begin{proof}
It suffices to show that there are no strictly $\sigma_{t\omega_q,\beta_q}$-semistable objects. To simplify the notations, we will work on the stability conditions $\sigma_{\omega_q(t),\beta_q}=(Z_t,\scha_{\omega_q,\beta_q})$ for $t>0$ where the ample class $\omega_q(t)$ is defined to be $\sqrt{2t+1}\omega_q$. In this case, one has
$$Z_t(r,\Delta,s)=(tr-\frac{1}{2}\Delta\cdot H-s)+(r+\frac{1}{2}\Delta\cdot H)\sqrt{2t+1}\sqrt{-1}$$
Suppose that an object $E\in\scha_{\omega_q,\beta_q}$ with Mukai vector $u_1$ is strictly $\sigma_{t\omega_q,\beta_q}$-semistable, then there exists an exact sequence $0\rightarrow A\rightarrow E\rightarrow B\rightarrow0$ in $\schp(\phi(E))$ such that $A$ is $\sigma_{t\omega_q,\beta_q}$-stable. One sets $v(A)=(r,\Delta,s)$ and obtains $2r+\Delta\cdot H=1$ and $2s=2r-2+(2r-1)t$.

Since $A$ is $\sigma_{t\omega_q,\beta_q}$-stable, one has $\Delta^2-2rs\geq-2$. Then $(\Delta\cdot H)^2-8rs\geq-8$ by Hodge index theorem. It follows $9\geq (8r^2-4r)t+4r^2-4r$ by plugging the equations for $\Delta\cdot H$ and $s$.

It is clear that $r=0$ is always a solution. In this case, one sees $\Delta\cdot H=1$. So the Hodge index theorem ensures that $\Delta^2$ is either $0$ or $-2$. The possibility $\Delta^2=0$ is ruled out by applying \cite[Theorem 1 (ii)]{Rei88} since $2H$ is very ample. Therefore, one has $v(B)=(1,-\Delta,n)$ for some integer $n$ with $2n=t>0$. It is impossible according to Proposition~\ref{QDS-u_1-wall-crossing-lemma}.

Otherwise, one has either $r\geq 1$ or $r\leq -1$. In the first case, one can see 
$$\frac{9}{8}\geq\frac{9}{8r}-\frac{1}{2}r+\frac{1}{2}\geq \frac{(2r-1)t}{2}>0$$
which implies that $2s-2r+2=(2r-1)t=2$. It follows $r=s=1$ and $\Delta\cdot H=-1$. The fact that $A$ is stable and the Hodge index theorem therefore ensure that $\Delta^2=0$. It is impossible again according to Reider's result \cite[Theorem 1 (ii)]{Rei88}. Now suppose that $r\leq -1$, one has
$$\frac{1}{8}\geq -\frac{9}{8r}+\frac{1}{2}r-\frac{1}{2}\geq \frac{(1-2r)t}{2}>0 $$
which implies that $(2r-1)t=2s-2r+2$ has no valid values. So we are done. 
\end{proof}
\end{proposition}

So stable objects in $M_{\sigma_{\omega_q,\beta_q}}(u_1)$ are Gieseker semistable sheaves. Such a sheaf is isomorphic to some $\shi_Z$ for a length $2$ subscheme $Z\subset S$. Together with Remark~\ref{QDS-u1-strict-semi-part}, one has

\begin{proposition}\label{QDS-u1-stable}
A stable object in the moduli space $M_{\sigma_{\omega_q,\beta_q}}(u_1)$ is isomorphic to $\shi_Z$ for a length $2$ subscheme $Z\subset S$ which is not in any line inside $S$. In particular, $\sigma_{\omega_q,\beta_q}$ is generic with respect to $u_1$ and $M_{\sigma_{\omega_q,\beta_q}}(u_1)$ is isomorphic to $S^{[2]}$ when $S$ does not contain any line.
\end{proposition}

In conclusion, one has a contraction $M_H(u_1)\twoheadrightarrow M_{\sigma_{\omega_q,\beta_q}}(u_1)$ which is an isomorphism over stable locus of $M_{\sigma_{\omega_q,\beta_q}}(u_1)$ and sends every $\shi_Z$ with $Z\subset L$ for some line $L$ contained in $S$ to the point represented by $\sho_S(-L)\oplus\sho_L(-2)$. Moreover, the induced involution 
$$\pi_q\colon M_{\sigma_{\omega_q,\beta_q}}(u_1)\rightarrow M_{\sigma_{\omega_q,\beta_q}}(u_1)$$ 
can be seen as a birational involution on $S^{[2]}$. Since the open subset $M^{\stable}_{\sigma_{\omega_q,\beta_q}}(u_1)\subset S^{[2]}$ has codimension two complement, one has a canonical identification $$H^2(M^{\stable}_{\sigma_{\omega_q,\beta_q}}(u_1),\bz)\cong H^2(S^{[2]},\bz)$$
which allows one to apply Proposition~\ref{K3-invariant-lattices-compatible} (see the comments after that proposition). Together with Proposition~\ref{QDS-induced-Hodge-isometry}, one concludes that

\begin{proposition}
The Hodge isometry $\pi_q^*$ on $H^2(S^{[2]},\bz)$ induced by $\pi_q$ is the reflection in the span of $\theta_{u_1}(1,-H,1)$ in $\NS(S^{[2]})$, where $\theta_{u_1}\colon u_1^{\perp}\rightarrow H^2(S^{[2]},\bz)$ is the Mukai isomorphism.
\end{proposition}

According to \cite[Section 4.1.3]{O'Grady05}, this is the same Hodge isometry of the birational involution constructed by Beauville in \cite{Bea83}. So $\pi_q\colon S^{[2]}\dashrightarrow S^{[2]}$ has to be the Beauville's involution. The fixed locus $\Fix(\pi_q)\subset M_{\sigma_{\omega_q,\beta_q}}(u_1)$ coincides with the surface $BTL_S$ of bitangent lines for our quartic K3 surface $S$.\footnote{Here we adopt the convention in \cite[Page 9]{Wel81} that $BTL_S$ also contains lines inside $S$.} One has $M_{\sigma}(-\mu_1)\cong \Sigma(X)$ for the quartic double solid $X$ over $S$ due to Proposition~\ref{QDS-class-mu1}. Moreover, the resulting covering map $\Sigma(X)\rightarrow BTL_S$ is exactly the one described in \cite[Page 15]{Wel81} geometrically. The only non-trivial part is 

\begin{proposition}\label{QDS-u1-covering-map}
The covering map $\Sigma(X)\rightarrow BTL_S$ sends a line $L\not\subset S$ to the unique length $2$ subscheme $Z\subset S$ such that $L\cap S=2Z$. In particular, two such lines $L_1$ and $L_2$ are sent to the same $Z$ if and only if their images under the branched covering $X\rightarrow\bp^3$ coincide.

\begin{proof}
Suppose that a line $L_1\not\subset S$ is sent to $Z$, then one has $$\funct{Forg}(\shi_Z)=\funct{Forg}(\funct{Inf}(\shi_{L_1}[1]))\cong \shi_{L_1}[1]\oplus\funct{U}(\shi_{L_1}[1])\cong\shi_{L_1}[1]\oplus\shi_{L_2}[1]$$
where $L_2$ is the unique line such that the images of $L_1,L_2$ under $X\rightarrow\bp^3$ coincide. It remains to show $Z\subset L_1$ or $Z\subset L_2$. At first, one considers the short exact sequence
$$0\rightarrow\sho_X(-2)\rightarrow\shi_{Z/X}\rightarrow \shi_{Z/S}\rightarrow0$$
and sees that $\funct{Forg}(\shi_Z)\cong\funct{L}_{\sho_X}\funct{L}_{\sho_X(1)}\shi_{Z/X}(2)$. Moreover, this subscheme $Z$ is not contained in any line inside $S$ according to Proposition~\ref{QDS-u1-stable}. So it uniquely determines a line in $\bp^3$ and two lines in $X$ that intersect with $S$ at $2Z$. Name one of the lines to be $L$, then
$$0\rightarrow\shi_L\rightarrow\shi_{Z/X}\rightarrow\sho_L(-2)\rightarrow0$$
together with $\funct{L}_{\sho_X}\funct{L}_{\sho_X(1)}\sho_L\cong\shi_L[1]$ ensures a morphism $\shi_{L_1}[1]\oplus\shi_{L_2}[1]\rightarrow\shi_L[1]$, which is trivial if and only if $\funct{L}_{\sho_X}\funct{L}_{\sho_X(1)}\shi_L(2)\cong\shi_L\oplus\shi_{L_1}[1]\oplus\shi_{L_2}[1]$. In this case, one has
$$\Hom(\shi_L\oplus\shi_{L_1}[1]\oplus\shi_{L_2}[1],\shi_L[3])\cong\Hom(\shi_L(2),\shi_L[3])\cong\Hom(\shi_L,\shi_L)\cong\bc$$
which is possible only when $L$ is one of $L_1$ or $L_2$ because $\Hom(\shi_L,\shi_L[3])=0$ for any $L\subset X$ and $\Hom(\shi_{L_1}\oplus\shi_{L_2},\shi_L[2])=0$ if $L$ is not one of $L_1$ or $L_2$.

Otherwise, one has a non-trivial morphism $\shi_{L_1}[1]\oplus\shi_{L_2}[1]\rightarrow\shi_L[1]$. It means that $L$ is either $L_1$ or $L_2$ again, because the ideal sheaves of lines are $\sigma$-stable objects.
\end{proof}
\end{proposition}

\subsection{Over the class $2\mu_1$}\label{QDS-ss-over-2mu1}
Here we are going to understand the covering morphism
$$\bigsqcup_{\forg(r,aH,s)=-2\mu_1}M_{\sigma_{\omega_q,\beta_q}}(r,aH,s)\rightarrow M_{\sigma}(-2\mu_1)^{\funct{U}}$$
for generic quartic K3 surfaces.

\begin{proposition}\label{QDS-(1,-H,3)}
Set $u=(-1,H,-3)$, then $M_{\sigma_{\omega_q,\beta_q}}(u)=\{\sho_S(-H)[1]\}$.
	
\begin{proof}
At first, we prove that $\sigma_{\omega_q,\beta_q}$ is $u$-generic. Since $u$ is primitive, it suffices to check that there are no strictly $\sigma_{\omega_q,\beta_q}$-semistable objects. Suppose that $E\in\scha_{\omega_q,\beta_q}$ with Mukai vector $u$ is strictly $\sigma_{\omega_q,\beta_q}$-semistable, then there exists an exact sequence $0\rightarrow A\rightarrow E\rightarrow B\rightarrow0$ in $\scha_{\omega_q,\beta_q}$ such that $A$ and $B$ are $\sigma_{\omega_q,\beta_q}$-semistable with $\mu(A)=\mu(B)$. Then one has
$$\Im(Z_{\omega_q,\beta_q}(A))=r+2a\geq0\quad\textup{and}\quad \Im(Z_{\omega_q,\beta_q}(B))=1-r-2a\geq0$$
after setting $v(A)=(r,aH,s)$. As $\mu(A)=\mu(B)$, it is impossible.

Then it suffices to check that $\sho_S(-H)[1]$ is $\sigma_{\omega_q,\beta_q}$-semistable. Otherwise, one can choose a destabilizing exact sequence 
$$0\rightarrow A\rightarrow\sho_S(-H)[1]\rightarrow B\rightarrow0$$
in $\scha_{\omega_q,\beta_q}$. It follows a long exact sequence of sheaves
$$0\rightarrow\shh^{-1}(A)\rightarrow\sho_S(-H)\rightarrow\shh^{-1}(B)\rightarrow\shh^0(A)\rightarrow0\rightarrow\shh^0(B)\rightarrow0$$
from which one can see immediately $\shh^0(B)=0$. Also, the rank of $\shh^{-1}(A)$ is $1$ as $\shh^{-1}(A)$ is torsion-free. It follows that the rank of $\shh^{-1}(B)$ should equal to the rank of $\shh^0(A)$, which is impossible as $\shh^{-1}(B)$ is torsion-free and admits a surjection $\shh^{-1}(B)\rightarrow\shh^0(A)$.
\end{proof}
\end{proposition}

The space $M_{\sigma}(-2\mu_1)$, that is described in \cite[Appendix A.3]{FLZ24}, contains shifts of Gieseker semistable sheaves and a two-term complex $G_X\cong \funct{L}_{\sho_X}\funct{L}_{\sho_X(1)}\sho_X(-1)[1]$. 

\begin{proposition}
The forgetful image of $\sho_S(-H)[1]$ in $\Kuc(X)$ is isomorphic to $G_X$.
	
\begin{proof}
One has a short exact sequence $0\rightarrow\sho_X(-2)\rightarrow\sho_X\rightarrow j_*\sho_S\rightarrow0$ and the desired forgetful image $\funct{L}_{\sho_X}\funct{L}_{\sho_X(1)}j_*\sho_S(H)$ is computed by applying the functor $\funct{L}_{\sho_X}\funct{L}_{\sho_X(1)}(-\otimes\sho_X(1))$.
\end{proof}
\end{proposition}

The object $\sho_S(-H)[1]$ is stable, so the preimage of the point corresponding to $G_X$ contains exactly two elements by Corollary~\ref{ME-double-cover-stable}. Another part of the preimage is $M_{\sigma_{\omega_q,\beta_q}}(3,-H,1)$.

\begin{proposition}\label{QDS-cover-over-2mu1}
The covering morphism induced by the forgetful functor has the form
$$M_{\sigma_{\omega_q,\beta_q}}(-1,H,-3)\sqcup M_{\sigma_{\omega_q,\beta_q}}(1,0,-1)\sqcup M_{\sigma_{\omega_q,\beta_q}}(3,-H,1)\rightarrow M_{\sigma}(-2\mu_1)^{\funct{U}} $$

\begin{proof}
According to Proposition~\ref{QDS-image-forg}, a vector $u=(r,aH,s)$ sent to $-2\mu_1$ satisfies $s+2a=-1$ and $r+s+4a=0$. So $\langle u,u\rangle_{\tilde{H}}\geq-2$ is equivalent to $|a|\leq 1$. It means that for $|a|>1$ the moduli space $M_{\sigma_{\omega_q,\beta_q}}(u)$ is either empty or consists of only strictly semistable objects. However, one can show as in Proposition~\ref{QDS-(1,-H,3)} that $M_{\sigma_{\omega_q,\beta_q}}(u)$ does not contain strictly semistable objects. 
\end{proof}
\end{proposition}

\begin{remark}
Using a similar argument together with Reider's method, it is possible to write down the covering morphism for any quartic K3 surface $S$. The source will have more summands corresponding to objects determined by lines and conics contained in $S$.
\end{remark}

It remains to understand the images of objects in $M_{\sigma_{\omega_q,\beta_q}}(u_1)$. At first, the object $\shi_Z$ is sent to $\shi_{L_1}[1]\oplus\shi_{L_2}[1]$ once $\shi_Z\in\{\funct{Inf}(\shi_{L_1}[1]),\funct{Inf}(\shi_{L_2}[1])\}$ and by Proposition~\ref{QDS-u1-covering-map} it holds only when $Z\in BTL_S$ and $L_1\cap L_2=2Z$. Any strictly semistable point in $M_{\sigma}(-2\mu_1)^{\funct{U}}$ has this form according to \cite[Appendix A]{FLZ24} and thus is the image of some $\shi_Z$ for $Z\in BTL_S$.

Moreover, the image of $\shi_Z$ in $M_{\sigma}(-2\mu_1)^{\funct{U}}$ is stable for any $Z\notin BTL_S$. Otherwise, there exist two lines $L_1$ and $L_2$ such that $\funct{Forg}(\shi_Z)=\shi_{L_1}[1]\oplus\shi_{L_2}[1]$. However, one has
$$\Hom(\funct{Forg}(\shi_Z),\funct{Forg}(\shi_{Z'}))=\Hom(\shi_Z,\shi_{Z'}\oplus\shi_{Z'})=0$$ 
for any $Z'\in BTL_S$, a contradiction. In conclusion, one has

\begin{proposition}\label{QDS-sharpen-FLZ-appendix}
The stable locus in $M_{\sigma}(-2\mu_1)^{\funct{U}}-\{G_X\}$ is the image of $S^{[2]}-BTL_S$ under the étale double covering. In particular, the locus is four-dimensional.
\end{proposition}

By virtue of Proposition~\ref{ME-smooth-points}, the points in $M_{\sigma}(-2\mu_1)^{\funct{U}}$ are singular. So the singular stable locus in $M_{\sigma}(-2\mu_1)$ is four-dimensional, refuting a guess in \cite[Remark A.16]{FLZ24}.

\begin{corollary}\label{QDS-scheme-structure-singular-locus}
Then the singular locus in $M^{\stable}_{\sigma}(-2\mu_1)$ is the disjoint union of a point and a four-dimensional component.

\begin{proof}
It remains to show that the images of $M_{\sigma_{\omega_q,\beta_q}}(-1,H,-3)$ and $M_{\sigma_{\omega_q,\beta_q}}(1,0,-1)$ under the covering map in Proposition~\ref{QDS-cover-over-2mu1} are topologically disjoint. It indeed follows from the fact that the union of $M_{\sigma_{\omega_q,\beta_q}}(-1,H,-3)$ and $M_{\sigma_{\omega_q,\beta_q}}(1,0,-1)$ in the source is topologically disjoint and both moduli spaces $M_{\sigma_{\omega_q,\beta_q}}(-1,H,-3)$ and $M_{\sigma_{\omega_q,\beta_q}}(1,0,-1)$ are irreducible.
\end{proof}
\end{corollary}

\subsection{On the vector $2u_1$}\label{QDS-ss-over-2u1}
Here we are going to understand the $\Pi_q$-induced involution
$$ M_{\sigma_{\omega_q,\beta_q}}(2u_1)\rightarrow M_{\sigma_{\omega_q,\beta_q}}(2u_1)$$
and the first step is to understand $M_{\sigma_{\omega_q,\beta_q}}(2u_1)$. Still, we will only consider generic quartic surfaces so that the polarization is $2u_1$-generic. In this case, $M_H(2u_1)$ is precisely the one considered in \cite{O'Grady99}. Therefore, any point in the singular locus $M^{\textsf{sin}}_H(2u_1)$ can be represented by $\shi_{Z_1}\oplus\shi_{Z_2}$ for some length 2 subschemes $Z_1,Z_2\subset S$ and \cite{LS06} asserts that blowing up $M_H(2u_1)$ at the reduced singular locus $M^{\textsf{sin}}_H(2u_1)\subset M_H(2u_1)$ gives an irreducible symplectic manifold $\tilde{M}$.

\begin{proposition}
A strictly $\sigma_{\omega_q,\beta_q}$-semistable object $E$ in $M_{\sigma_{\omega_q,\beta_q}}(2u_1)$ is either an extension of two objects in $M_{\sigma_{\omega_q,\beta_q}}(u_1)$ or an extension of $\Pi_q(\sho_S(-H))[1]$ and $\sho_S(-H)[1]$.

\begin{proof}
The object $E$ admits the Harder--Narasimhan filtration 
$$0\rightarrow A\rightarrow E\rightarrow B\rightarrow0$$
where $A,B$ are $\sigma_{\omega_q,\beta_q}$-semistable in $\scha_{\omega_q,\beta_q}$ such that $\mu(A)=\mu(B)$. Set $v(A)=(r,aH,s)$, then one can obtain as before that $r+s+4a=0$ and $s+2a=-1$. So one concludes.
\end{proof}
\end{proposition}

It remains to determine the $\sigma_{\omega_q,\beta_q}$-stable objects in $M_{\sigma_{\omega_q,\beta_q}}(2u_1)$. The idea is similar to that for Proposition~\ref{QDS-u_1-wall-crossing}, but we have to apply \cite[Theorem 5.1 and Theorem 5.3]{MZ16}.

\begin{proposition}
The stability condition $\sigma_{t\omega_q,\beta_q}$ is $2u_1$-generic for any $t> 1$, and is on a non-totally semistable wall with respect to $2u_1$ for $t=1$.

\begin{proof}
To simplify the notations, we will work on $\sigma_{\omega_q(t),\beta_q}=(Z_t,\scha_t)$ for $t\geq0$ where the ample class $\omega_q(t)$ is defined to be $\sqrt{2t+1}\omega_q$. In this case, one has
$$Z_t(r,aH,s)=(tr-2a-s)+(r+2a)\sqrt{2t+1}\sqrt{-1}$$

Let $W$ be a potential wall for $\sigma_{\omega_q(t),\beta_q}$ with respect to $2u_1$, then the associated hyperbolic lattice $H_W$ consists of classes $(r,aH,s)$ satisfying $r+2(t+2)a+s=0$. So at first $W$ can be an actual wall only when $t$ is a rational number. 

At first, we prove that $W$ can never be a totally semistable wall. It suffices to check that there is no effective spherical class $v_0$ in $H_W$ such that $\langle v_0,2u_1\rangle_{\tilde{H}}<0$. Otherwise, one chooses such a class $(r,aH,s)$. Then one has $rs=2a^2+1,r+s=-2(t+2)a,s>r$ and
$$\Re(Z_t(r,aH,s)Z_t(2u_1)^{-1})=\frac{(tr-2a-s)(t+1)+(r+2a)(2t+1)}{2(t^2+4t+2)}>0$$
The first three conditions allow one to get
$$r=-(t+2)a-\sqrt{(t^2+4t+2)a^2-1}\quad\textup{and}\quad s=-(t+2)a+\sqrt{(t^2+4t+2)a^2-1}$$
Therefore, up to some elementary computations, the last inequality becomes
$$(t^2+4t+2)(at+\sqrt{(t^2+4t+2)a^2-1})<0$$
or equivalently $at+\sqrt{(t^2+4t+2)a^2-1}<0$ due to $t^2+4t+2>0$, from which one sees $a<0$ and $(4t+2)a^2<1$. It is impossible.

Then we will show for $v_0=(r,aH,s)$ that $\langle 2u_1,v_0\rangle_{\tilde{H}}=\langle v_0,2u_1\rangle_{\tilde{H}}=2r-2s$ cannot be $2$ when $v_0$ is an isotropic class, and cannot be $0$ or $2$ when $v_0$ is a spherical class. Moreover, we will show that $2r-2s$ cannot be $4$ for any $t>0$. These conclude the proposition.

Suppose that $v_0=(r,aH,s)$ is isotropic, then $r$ and $s$ are integer solutions of the equation
$$x^2+2(t+2)ax+2a^2=0$$ 
So $|r-s|=2\sqrt{(t+2)^2a^2-2a^2}$, and it cannot equal to $1$ for $t\geq0$ and $a\in\bz$.

Suppose that $v_0=(r,aH,s)$ is spherical, then $r$ and $s$ are integer solutions of the equation 
$$x^2+2(t+2)ax+2a^2+1=0$$
So $|r-s|=2\sqrt{(t+2)^2a^2-2a^2-1}$, and it cannot equal to $0$ or $1$ for $t\geq0$ and $a\in\bz$. It is equal to $2$ if and only if $(t+2)^2=4$ and $a^2=1$, corresponding to $t=0$ and $a=\pm 1$. In this case one has spherical classes $(3,- H,1)$ and $(-1,H,-3)$ making $W$ a flopping wall.
\end{proof}
\end{proposition}

Therefore $M_{\sigma_{t\omega_q,\beta_q}}(2u_1)=M_H(2u_1)$ for any $t>1$ due to Theorem~\ref{K3-large-volume-limits}. Moreover, we know that singular points in $M_H(2u_1)$ are contained in $M_{\sigma_{\omega_q,\beta_q}}(2u_1)$ and stable objects in $M_{\sigma_{\omega_q,\beta_q}}(2u_1)$ are stable sheaves. The only semistable point in $M_{\sigma_{\omega_q,\beta_q}}(2u_1)$ irrelevant to $M_H(2u_1)$ is the one represented by $\sho_S(-H)[1]\oplus\Pi_q(\sho_S(-H))[1]$ and is fixed by $\Pi_q$. 

Hence $\Pi_q$ induces a birational involution 
$$\pi_q\colon M_H(2u_1)\dashrightarrow M_H(2u_1)$$  
which preserves $M^{\textsf{sin}}_H(2u_1)$ by construction. So it also admits a birational involution $\tilde{\pi}_q$ on the irreducible symplectic manifold $\tilde{M}$. As we have mentioned, one can work out the invariant lattice of the birational involution $\tilde{\pi}_q$ using the Hodge isometry $(2u_1)^{\perp}\rightarrow H^2(M_H(2u_1),\bz)$.

\begin{proposition}
Consider the blow-up morphism $\rho\colon \tilde{M}\rightarrow M_H(2u_1)$, the exceptional divisor $E\tilde{E}$, the Weil divisor $B\subset M_H(2u_1)$ parameterizing sheaves that are not locally free and the strict transform $\tilde{B}$ of $B$ under $\rho$, then there exists an isometric embedding of Hodge structures
$$(2u_1)^{\perp}\stackrel{\tilde{\theta}}{\cong}H^2(M_H(2u_1),\bz)\stackrel{\rho^*}{\hookrightarrow} H^2(\tilde{M},\bz),\quad (r,\Delta,r)\mapsto \tilde{\theta}(\Delta)+2rc_1(\tilde{B})+rc_1(\tilde{E})$$
where the Hodge structure on $(2u_1)^{\perp}=H^2(S,\bz)\oplus\bz(1,0,1)$ is induced from the Mukai lattice and the right hand side one is given by the Beauville form. Under the Beauville form, one has an orthogonal decomposition $H^2(\tilde{M},\bz)\cong \tilde{\theta}(H^2(S,\bz))\oplus\Lambda$ where $\Lambda$ is represented by the matrix
$$\begin{pmatrix}
-2&3\\3&-6
\end{pmatrix}$$
with the basis $c_1(\tilde{B})$ and $c_1(\tilde{E})$. Moreover, the following diagram 
\begin{displaymath}
	\xymatrix{
		(2u_1)^{\perp}\ar[d]_{\Pi_{q*}}\ar[rr]^{\theta}&&H^2(M_H(2u_1),\bz)\ar[d]_{\pi_{q*}}\ar[rr]^{\rho^*}&&H^2(\tilde{M},\bz)\ar[d]_{\tilde{\pi}_{q*}}\\
		(2u_1)^{\perp}\ar[rr]^{\theta}&&H^2(M_H(2u_1),\bz)\ar[rr]^{\rho^*}&&H^2(\tilde{M},\bz)\\
}\end{displaymath}
is commutative.

\begin{proof}
The construction of $\tilde{\theta}$ is contained in \cite[Theorems 1.1, 3.1 and 4.3]{Rap08}, and the non-trivial part of the commutative diagram is proved in \cite[Lemma 3.11]{PR13}. 
\end{proof}
\end{proposition}

It follows that the invariant lattice of the non-symplectic involution $\tilde{\pi}_q$ is generated by the classes $2c_1(\tilde{B})-\tilde{\theta}(H)$ and $c_1(\tilde{E})$. One checks that, up to an invertible base change, the invariant lattice is isomorphic to $\bz(2)\oplus\bz(-6)$ which is listed in \cite[No.1 in Table 6]{BG25}.

\subsection{Over the class $-\mu_2$}\label{QDS-ss-over-mu2}
In this subsection, we want to understand the covering morphism
$$\bigsqcup_{\forg(r,\Delta,s)=-\mu_2}M_{\sigma_{\omega_q,\beta_q}}(r,\Delta,s)\rightarrow M_{\sigma}(-\mu_2)^{\funct{U}}$$
for any quartic K3 surface. At first, all objects in $M_{\sigma}(-\mu_2)$ are stable due to \cite{APR22,FLZ24}, so objects in the source are also all stable. In particular, similar to Proposition~\ref{QDS-cover-over-2mu1}, one obtains

\begin{proposition}
Consider a class $(r,aH,s)\in N(S)$ with $\forg(r,aH,s)=-\mu_2$, then the moduli space $M_{\sigma_{\omega_q,\beta_q}}(r,aH,s)$ is non-empty only when $a=0,1$.
\end{proposition}

In the non-generic case, a Mukai vector $(r,\Delta,s)$ sent to $-\mu_2$ has the form $(1+k,\Delta,k)$ such that $\Delta\cdot H=-2k$ for some integer $k$. Then the Hodge index theorem gives
$$-2\leq\Delta^2-2rs=\Delta^2-\frac{1}{2}(\Delta\cdot H)^2+\Delta\cdot H\leq-k^2+2k$$
so $k=0,1,2$. Hence $\forg^{-1}(-\mu_2)$ is the union of $\{(1,0,0),(-1,H,-2)\}$ and the following sets
$$
\{(0,\Delta,-1)\,|\,\Delta\cdot H=2\textup{ and }\Delta^2\in\{0,-2\}\}
$$ $$
\{(1,\Delta,0)\,|\,\Delta\cdot H=0\textup{ and }\Delta^2=-2\}\cup\{(-1,\Delta,-2)\,|\,\Delta\cdot H=4\textup{ and }\Delta^2=2\}
$$
Since $H$ is very ample, the case where $\Delta^2=0$ and $\Delta\cdot H=2$ cannot happen. So all the Mukai vectors in the set $\forg^{-1}(-\mu_2)$ besides $(1,0,0)$ and $(-1,H,-2)$ are spherical classes.

\begin{corollary}
The fixed locus $M_{\sigma}(-\mu_2)^{\funct{U}}$ is isomorphic to the union of $M_{\sigma_{\omega_q,\beta_q}}(1,0,0)$ and finitely many reduced points. When $S$ is general, it is isomorphic to $M_{\sigma_{\omega_q,\beta_q}}(1,0,0)$.

\begin{proof}
One notices that the source of the covering morphism over $M_{\sigma}(-\mu_2)^{\funct{U}}$ can be written into a disjoint union $M_1\sqcup M_2$ such that $\Pi_q$ induces an isomorphism between $M_1$ and $M_2$. So the fixed locus $M_{\sigma}(-\mu_2)^{\funct{U}}$ must be isomorphic to $M_1$, as the covering morphism is étale of degree $2$. The discussions above assert that $M_1$ is exactly what we want it to be.
\end{proof}
\end{corollary}

\begin{proposition}\label{QDS-(1,0,0)-moduli-space}
The moduli space $M_{\sigma_{\omega_q,\beta_q}}(1,0,0)$ is identified with $M_H(1,0,0)\cong S$.

\begin{proof}
Since $\sigma_{\omega_q,\beta_q}$ is generic with respect to $(1,0,0)$, the moduli space $M_{\sigma_{\omega_q,\beta_q}}(1,0,0)$ is a projective K3 surface. Also, one checks that $\shi_x$ is $\sigma_{\omega_q,\beta_q}$-semistable for each $x\in S$. It follows an injection $S\cong M_H(1,0,0)\hookrightarrow M_{\sigma_{\omega_q,\beta_q}}(1,0,0)$, which has to be an isomorphism.
\end{proof}
\end{proposition}

In conclusion, we have reduced the description of $M_{\sigma}(-\mu_2)^{\funct{U}}$ to enumerating certain classes in $\NS(S)$, which recovers a part of \cite[Theorem 5.2]{APR22} and sharpens \cite[Lemma 7.5]{FLZ24}.

\subsection{Over the vector $u_1-u_2$}\label{QDS-ss-u1minusu2}
In the end, we will try to understand the induced involution
$$\pi_q\colon M_{\sigma_{\omega_q,\beta_q}}(u_1-u_2)\rightarrow M_{\sigma_{\omega_q,\beta_q}}(u_1-u_2)$$
for generic quartic K3 surfaces. At first, we need to work out $M_{\sigma_{t\omega_q,\beta_q}}(u_1-u_2)$.

\begin{proposition}\label{QDS-(0,H,-2)-wall-crossing}
The stability condition $\sigma_{t\omega_q,\beta_q}$ is $(u_1-u_2)$-generic for $t\in(1,+\infty)\setminus\{\sqrt{3}\}$.
	
\begin{proof}
We work again on $\sigma_{\omega_q(t),\beta_q}=(Z_t,\scha_t)$ for $t\geq0$ where $\omega_q(t):=\sqrt{2t+1}\omega_q$. It suffices to show that any semistable object is stable. Otherwise, there exists a short exact sequence $0\rightarrow A\rightarrow E\rightarrow B\rightarrow0$ in $\scha_{\omega_q,\beta_q}$ such that $A$ is stable and $\mu(A)=\mu(B)$. 
		
Set $v(A)=(r,aH,s)$, then $r=1-2a$ and $s=t(1-2a)-2a$. Since $4a^2-2rs\geq-2$, one obtains that $r=1,a=0$ and then $s=t=1$ or $r=0,a=1$ and then $t=0,s=-2$. 
\end{proof}
\end{proposition}

Together with Theorem~\ref{K3-large-volume-limits}, one concludes from this proposition that $M_{\sigma_{t\omega_q,\beta_q}}(u_1-u_2)$ is identified with $M_H(u_1-u_2)$ for $t>\sqrt{3}$. The wall for $t=\sqrt{3}$ is a flopping wall induced by the spherical vectors $(1,0,1)$ and $(-1,H,-3)$. Set $M:=M_{\sigma_{t\omega_q,\beta_q}}(u_1-u_2)$ for $\sqrt{3}>t>1$, then the exceptional locus of the flop $M_H(u_1-u_2)\dashrightarrow M$ consists of those $i_{C*}\sho_C$ for curves $i_{C}\colon C\subset S$ with $C\in|H|$ and is therefore isomorphic to $\bp^3$ (see \cite[Example 3.21]{FMOS22}).

\begin{proposition}
Consider the Mukai isometry $\theta\colon (u_1-u_2)^{\perp}\rightarrow H^2(M,\bz)$. Then the extreme rays for the movable cone $\Mov(M)$ are generated by $\theta(0,0,-1)$ and $\theta(2,-H,0)$ and the extreme rays for the ample cone $\Amp(M)$ are given by $\theta(2,-H,-2)$ and $\theta(2,-H,0)$.

\begin{proof}
The extreme rays for the closed positive cone $\overline{\Pos(M)}$ are given by the classes $\theta(0,0,-1)$ and $\theta(2,-H,1)$, then one applies \cite[Theorem 12.3]{BM14-MMP} to compute the movable cones. 

Set $v=(r,aH,s)$ and suppose that $v^{\perp}\cap (u_1-u_2)^{\perp}\neq\varnothing$, then $(2m,-mH,n)\in v^{\perp}$ for some $m,n\in\bz$. It is not hard to check that the vector $(2m,-mH,n)$ cannot be spherical. Then $v$ has to be isotropic and satisfies $1\leq 4a+2r\leq 2$. It follows $2a+r=1$ and $2a=rs$. In particular
$$s=\frac{2a^2}{1-2a}\in\bz$$
means that we must have $a=1$ or $a=0$. So $v$ can only be $(1,0,0)$ or $(-1,H,-2)$. The corresponding vector in $v^{\perp}\cap (u_1-u_2)^{\perp}$ is $(2,-H,0)$ for both $(1,0,0)$ and $(-1,H,-2)$. 

Since we have already known the flopping walls, the ample cone $\Amp(M)$ can be computed similarly by \cite[Theorem 12.1]{BM14-MMP} or directly using \cite[Lemma 9.2]{BM14}.
\end{proof}
\end{proposition}

\begin{corollary}
One has $\Aut(M)=\Bir(M)=\bz/2\bz$.

\begin{proof}
One can see $\Aut(M)=\Bir(M)\subseteq\bz/2\bz$ by inspecting the argument of \cite[Theorem 4.9]{De22} and the fact that $\theta(0,0,-1)$ and $\theta(2,-H,0)$ cannot be exchanged by any Hodge isometry on the lattice $\NS(M)$. It is an equality as one can find a non-trivial element in it.

Indeed, the anti-autoequivalence $\deriver\HHom(-,\sho_S(-H))[1]$ on $\cate{D}^b(S)$ induces an involution $\pi_D$ on $M_H(u_1-u_2)$ according to \cite[Section 3]{FMOS22}. Via the flopping $M_H(u_1-u_2)\dashrightarrow M$, it induces a birational involution $\pi_D\colon M\dashrightarrow M$ which is clearly nontrivial. 
\end{proof}
\end{corollary}

Consequently, one has a regular non-trivial involution $\pi_M$ on $M$. Additionally, any stable object $A$ in the flopping locus $\bp^3\subset M$ satisfies $\shh^{-1}(A)\cong\sho_S(-H)$. Therefore, one has

\begin{corollary}
The fixed locus of $\pi_M$ contains the closure $\Omega_0$ of the locus parameterizing sheaves $i_{C*}\shr$ where $C\in|\sho_S(H)|$ is a smooth curve and $\shr$ is a non-trivial square root of $\sho_C$.

\begin{proof}
This follows from \cite[Proposition 4.1]{FMOS22} and our description of the flopping locus.
\end{proof}
\end{corollary}

The wall for $t=1$ is a divisorial wall induced by $(1,0,0)$ and $(-1,H,-2)$ and one has a divisorial contraction $M\rightarrow M_+\subset M_{\sigma_{\omega_q,\beta_q}}(u_1-u_2)$. The contracted locus inside $M$ is isomorphic to a $\bp^2$-bundle over $S\times S$ following \cite[Example 3.11]{FMOS22} and Proposition~\ref{QDS-(1,0,0)-moduli-space}. Moreover, the involution $\pi_q$ induced by $\Pi_q$ admits a non-trivial birational involution on $M$ which can be extended to the same regular involution $\pi_M$ on $M$. The divisorial contracted locus $T\subset M$ contains a stable object $B$ fitting into an exact triangle $$\shi_x(-H)\rightarrow\shi_y\rightarrow B\rightarrow\shi_x(-H)[1]$$ 
for some points $x,y\in S$. One notices that $\Hom(\sho_S,\shi_x)=0$, so $T$ does not intersect the flopping locus $\bp^3\subset M$ and one gets a total point-wisely description of $\pi_M$.

On the other hand, the moduli space $M_{\sigma}(-\mu_2)$ is isomorphic to $X\cup \tilde{\Omega}$ (see \cite{APR22,FLZ24}) such that $X\cap\tilde{\Omega}$ is the fixed locus $M_{\sigma}(-\mu_2)^{\funct{U}}\cong S$. So due to the branched double covering 
$$M_{\sigma}(-\mu_2)\rightarrow M_{\sigma_{\omega_q,\beta_q}}(u_1-u_2)^{\Pi_q}$$
the fixed locus in $M_{\sigma_{\omega_q,\beta_q}}(u_1-u_2)$ is isomorphic to $\bp^3\cup \Omega$ for some 3-dimensional locus $\Omega$ inside $M_{\sigma_{\omega_q,\beta_q}}(u_1-u_2)$ such that the branch locus of the covering is $\bp^3\cap\Omega\cong S$. In addition, the branched double covering component $X\rightarrow\bp^3$ recovers the quartic double solid. 

\begin{remark}
The identification of the component $\Omega\subset M_{\sigma_{\omega_q,\beta_q}}(u_1-u_2)$ is not clear. The main difficulty is to determine which part of $\Omega_0$ is contracted via the divisorial contraction.
\end{remark}

\section{Kuznetsov components of special Gushel--Mukai threefolds}\label{ST-Family-special-GM}
In this section, we apply this framework to the Enriques categories over K3 surfaces described in Example~\ref{En_example-2}. In the first three subsections, we work out the (essentially unique) stability condition $\sigma_1$ on $\cate{D}^b(S)$ invariant under the $\bz/2\bz$-action, write down the induced involutions on moduli spaces with respect to $\sigma_1$ in Proposition~\ref{ME-involution}, and write down the explicit covering maps between moduli spaces as in Proposition~\ref{ME-double-cover}. In the remaining subsections, we describe the involutions and covering maps for some explicit numerical classes.

\subsection{Two constructions of the Enriques categories}
We start with a Brill--Noether general K3 surface $S$ polarized by a very ample divisor class $H$ with $H^2=10$. Such a K3 surface can be realized as a divisor in a codimension 3 linear section $W$ of $\Gr(2,5)$. The special Gushel--Mukai threefold $X$ over $S$ is a double cover $X\rightarrow W$ branched over $S$, and it carries a unique primitive ample class $h$ such that $h^3=10$. In this section, we assume that $X$ and $W$ are smooth.

One has two semiorthogonal decompositions
$$\cate{D}^b(X)=\langle\schk_X^1,\sho_X,\shu_X^\vee\rangle\quad\textup{and}\quad \cate{D}^b(X)=\langle\schk_X^2,\shu_X,\sho_X\rangle$$
where $\shu_X$ is the tautological vector bundle pulled back from that of $\Gr(2,5)$. Then the admissible components $\schk_X^1$ and $\schk_X^2$ are equivalent Enriques categories, with the same action generated by the pushforward $\funct{U}$ along the covering involution of the double covering $X\rightarrow W$. The equivariant categories are both equivalent to $\cate{D}^b(S)$ according to \cite{KP17}.

The numerical Grothendieck group $N(\schk_X^1)$ has rank 2, with a basis $\kappa_1$, $\kappa_2$ such that
$$\ch(\kappa_1)=1-\frac{1}{5}h^2,\quad \ch(\kappa_2)=2-h+\frac{1}{12}h^3$$
and the Euler form is given in this basis by $-I_2$. Equivalently, the numerical Grothendieck group $N(\schk_X^2)$ has rank 2, with a basis $\kappa'_1$, $\kappa'_2$ such that
$$\ch(\kappa'_1)=1-h+\frac{3}{10}h^2+\frac{1}{30}h^3,\quad \ch(\kappa'_2)=-2+h-\frac{1}{12}h^3$$
and the Euler form is also given in this basis by $-I_2$.

Moreover, one notices that the induced actions of $\funct{U}$ on $N(\schk_X^1)$ and $N(\schk_X^2)$ are trivial as the covering involution $X\rightarrow W$ fixes the polarization class $h$. So the covering involution induces an involution on moduli spaces of $\sigma$-semistable objects for any numerical class in $N(\schk_X^1)$ or $N(\schk_X^2)$ and any $\funct{U}$-fixed stability condition $\sigma$.

According to \cite{KP17}, the residual group actions related to $\schk_X^1$ and $\schk_X^2$ are generated by different involutions:
\begin{align*}
	\Pi_1:=&\funct{T}_{\sho_S}\circ\funct{T}_{\shu_S^\vee}\circ(-\otimes\sho_S(H))[-1]\\
	\Pi_2:=&\funct{T}_{\shu_S}\circ\funct{T}_{\sho_S}\circ(-\otimes\sho_S(H))[-1]
\end{align*}
where $\shu_S$ is the unique stable spherical vector bundle in $M_H(2,-H,3)$. In fact, it is isomorphic to the restriction of the tautological vector bundle on $\Gr(2,5)$. Also, one has
$$\Pi_2\cong\funct{O}_{GM}^{-1}\circ\Pi_1\circ\funct{O}_{GM}$$
where $\funct{O}_{GM}:=\funct{T}_{\sho_S}(-\otimes\sho_S(H))$. A direct computation shows that 

\begin{proposition}
The Hodge isometry $\tau_1$ induced by $\Pi_1$ is given by the mapping
$$(r,\Delta,s)\mapsto(-4r-5s-2\Delta\cdot H,(2r+2s+\Delta\cdot H)H-\Delta,-5r-4s-2\Delta\cdot H)$$
and the Hodge isometry $\tau_2$ induced by $\Pi_2$ is given by the mapping
$$(r,\Delta,s)\mapsto(-9r-5s-3\Delta\cdot H,(6r+3s+2\Delta\cdot H)H-\Delta,-20r-9s-6\Delta\cdot H)$$
\end{proposition}

Moreover, one checks as in Proposition~\ref{QDS-fixed-stability-condition} that

\begin{proposition}\label{GM-fixed-stability-condition}
The equation 
$$Z_{\omega,\beta}(r,\Delta,s)=Z_{\omega,\beta}(-4r-5s-2\Delta\cdot H,(2r+2s+\Delta\cdot H)H-\Delta,-5r-4s-2\Delta\cdot H)$$
holds for any $(r,\Delta,s)\in N(S)$ if and only if $\omega=\frac{1}{5}H$ and $\beta=-\frac{2}{5}H$.

Also, the equation 
$$Z_{\omega,\beta}(r,\Delta,s)=Z_{\omega,\beta}(-9r-5s-3\Delta\cdot H,(6r+3s+2\Delta\cdot H)H-\Delta,-20r-9s-6\Delta\cdot H)$$
holds for any $(r,\Delta,s)\in N(S)$ if and only if $\omega=\frac{1}{5}H$ and $\beta=-\frac{3}{5}H$.
\end{proposition}

Henceforth, we will denote $\beta_1=-\frac{2}{5}H,\beta_2=-\frac{3}{5}H$ and $\omega_1=\omega_2=\frac{1}{5}H$. One can check that $\sigma_{\omega_1,\beta_1}$ is a stability condition if and only if $\sigma_{\omega_2,\beta_2}$ is. Similar to Corollary~\ref{QDS-stability-condition-exists}, one has

\begin{corollary}\label{GM-stability-condition-exists}
Consider $i\in\{1,2\}$ and a $\funct{U}$-fixed stability condition $\sigma$ on $\schk_X^i$ such that $\sigma^{\funct{U}}$ is a geometric stability condition, then $\sigma^{\funct{U}}=\sigma_{\omega_i,\beta_i}.\tilde{g}$ for a unique $\tilde{g}\in\grp$. In particular, $\sigma_{\omega_i,\beta_i}$ is indeed a stability condition and the involution $\Pi_i$ fixes $\sigma_{\omega_i,\beta_i}$, if such $\sigma$ exists.
\end{corollary}

Such a $\funct{U}$-fixed stability condition $\sigma$ is constructed in \cite[Theorem 6.8]{BLMS23} and it satisfies the conditions in Corollary~\ref{GM-stability-condition-exists} according to \cite[Lemma 7.1 and Lemma 7.2]{PPZ26}. Therefore, one can use this corollary to recover \cite[Theorem A.10]{JLLZ24} similar to the quartic double solid case.

\subsection{The induced involutions}
One computes that the invariant lattice of the Hodge isometry $\tau_1$ is spanned by the Mukai vectors $w_1=(1,0,-1)$ and $w_2=(2,-H,2)$. Moreover, it is computed in \cite[Lemma 7.4]{BP23} that the invariant lattice of the Hodge isometry $\tau_2$ is spanned by the Mukai vectors $w'_1=(1,-H,4)$ and $w'_2=(-2,H,-2)$. In particular, one can see that

\begin{proposition}\label{GM-rotation-identification}
Consider a Brill--Noether general K3 surface $(S,H)$ of degree $10$, then 
$$M_{\sigma_{\omega_2,\beta_2}}(r,\Delta,s)\cong M_{\sigma_{\omega_1,\beta_1}}(-5r-s-\Delta\cdot H,rH+\Delta,-r)$$
for any $(r,\Delta,s)\in N(S)$. In particular, one has $M_{\sigma_{\omega_1,\beta_1}}(aw_1+bw_2)\cong M_{\sigma_{\omega_2,\beta_2}}(aw'_1+bw'_2)$.

\begin{proof}
Similar to Proposition~\ref{QDS-rotation-identification}, one computes the stability condition $\funct{O}^{-1}_{GM}.\sigma_{\omega_1,\beta_1}$. It turns out that the result $\funct{O}^{-1}_{GM}.\sigma_{\omega_1,\beta_1}$ is nothing but $\sigma_{\omega_2,\beta_2}$ up to an action of $\grp$. So by definition an object $A$ is (semi)stable for $\sigma_{\omega_2,\beta_2}$ if and only if $\funct{O}_{GM}A$ is $\sigma_{\omega_1,\beta_1}$-(semi)stable.

The Hodge isometry induced by the autoequivalence $\funct{O}_{GM}$ is 
$$(r,\Delta,s)\mapsto (-5r-s-\Delta\cdot H,rH+\Delta,-r)$$
and $aw'_1+bw'_2$ is sent to $aw_1+bw_2$ under the induced action of $\funct{O}_{GM}$. One thus concludes.
\end{proof}
\end{proposition}

Also, due to the conjugation of $\Pi_1$ and $\Pi_2$, the involutions and covering morphisms they induce are the same. Similar to the case of quartic double solids, one has

\begin{proposition}\label{GM-induced-involution}
Consider a Brill--Noether general K3 surface $(S,H)$ of degree $10$, then the autoequivalence $\Pi_1$ of $\cate{D}^b(S)$ induces an involution
$$\pi_{1}\colon M_{\sigma_{\omega_1,\beta_1}}(aw_1+bw_2)\rightarrow M_{\sigma_{\omega_1,\beta_1}}(aw_1+bw_2)$$
for each pair of $a,b\in\bz$. Equivalently, the autoequivalence $\Pi_2$ induces an involution
$$\pi_{2}\colon M_{\sigma_{\omega_2,\beta_2}}(aw'_1+bw'_2)\rightarrow M_{\sigma_{\omega_2,\beta_2}}(aw'_1+bw'_2)$$
for each pair of $a,b\in\bz$.
\end{proposition}

\subsection{The covering morphisms in two directions}\label{GM-covering-morphism}
To explicitly write down the covering maps given in Proposition~\ref{ME-double-cover}, we need to figure out the images $\forg_i(v)$ and $\sfinf_i(\mu)$ of the homomorphisms induced by $N(S)\cong   N(\schk_X^{i,\funct{U}})\rightarrow N(\schk_X^i)$ and $N(\schk_X^i) \rightarrow N(\schk_X^{i,\funct{U}})\cong N(S)$.

\begin{proposition}\label{GM-image-forg}
Consider the homomorphisms $\forg_i\colon N(S)\cong N(\schk_X^{i,\funct{U}})\rightarrow N(\schk_X^i)$, then
\begin{align*}
\forg_1(r,\Delta,s)&=(s-r)\kappa_1+(2r+2s+\Delta\cdot H)\kappa_2\\
\forg_2(r,\Delta,s)&=(4r+s+\Delta\cdot H)\kappa'_1-(2r+2s+\Delta\cdot H)\kappa'_2
\end{align*}
for any Mukai vector $(r,\Delta,s)\in N(S)$.
	
\begin{proof}
Set $\forg_1(r,\Delta,s)=a\kappa_1+b\kappa_2$. One computes that
$$\sfinf_1(\forg_1(r,\Delta,s))=(-3r-5s-2\Delta\cdot H,(2r+2s+\Delta\cdot H)H,-5r-3s-2\Delta\cdot H)$$ 
so \cite[Lemma 7.4]{PPZ26} and the adjunction between the forgetful functor and inflation functor give
\begin{align*}
			4b=-\chi(\forg_1(1,0,1),\forg_1(r,\Delta,s))=\langle(1,0,1),\sfinf_1(\forg_1(r,\Delta,s))\rangle_{\tilde{H}}&=8r+8s+4\Delta\cdot H\\
			a+2b=-\chi(\forg_1(0,0,1),\forg_1(r,\Delta,s))=\langle(0,0,1),\sfinf_1(\forg_1(r,\Delta,s))\rangle_{\tilde{H}}&=3r+5s+2\Delta\cdot H
\end{align*}
with which one computes that $a=s-r$ and $b=2r+2s+\Delta\cdot H$.

The result about $\forg_2$ is due to $\funct{L}_{\sho_X}(\schk_X^2\otimes\sho_X(1))\cong\schk_X^1$. Actually, one can compute as in \cite[Lemma 7.4]{PPZ26} that $\forg_2(1,0,1)=5\kappa'_1+4\kappa'_2$ and $\forg_2(0,0,1)=\kappa'_1+2\kappa'_2$.
\end{proof}
\end{proposition}

Since $\funct{U}$ acts trivially on $N(\schk_X^i)$ for $i=1,2$, the composite $\forg_i\circ \sfinf_i\colon N(\schk_X^i)\rightarrow N(\schk_X^i)$ is the multiplication by $2$. One can therefore deduce as Corollary~\ref{QDS-image-inflation} that

\begin{corollary}\label{GM-image-inflation}
One has
$$\sfinf_1(a\kappa_1+b\kappa_2)=-aw_1-bw_2\quad\textup{and}\quad \sfinf_2(a\kappa'_1+b\kappa'_2)=-aw'_1-bw'_2$$
for any integers $a$ and $b$. 
\end{corollary}

Therefore, one can see in the same way as Proposition~\ref{QDS-double-cover} that

\begin{proposition}
Consider a Brill--Noether general K3 surface $(S,H)$ of degree $10$ and the special Gushel--Mukai threefold $X\rightarrow W$ branched over $S$.

Then a $\funct{U}$-fixed stability condition $\sigma$ on $\schk_X^1$ admits a surjection
$$M_{\sigma}(-a\kappa_1-b\kappa_2)\rightarrow M_{\sigma_{\omega_1,\beta_1}}(aw_1+bw_2)^{\Pi_1}$$
for any integers $a,b$. This surjection has degree $2$ at points corresponding to stable objects and is étale at points corresponding to stable objects in $M_{\sigma}(-a\kappa_1-b\kappa_2)$ which are not fixed by $\funct{U}$.

Equivalently, a $\funct{U}$-fixed stability condition $\sigma$ on $\schk_X^2$ admits a surjection
$$M_{\sigma}(-a\kappa'_1-b\kappa'_2)\rightarrow M_{\sigma_{\omega_2,\beta_2}}(aw'_1+bw'_2)^{\Pi_2}$$
for any integers $a,b$. This surjection has degree $2$ at points corresponding to stable objects and is étale at points corresponding to stable objects in $M_{\sigma}(-a\kappa'_1-b\kappa'_2)$ which are not fixed by $\funct{U}$.
\end{proposition}

\subsection{Over the class $\kappa_1$}
Let us start with the simplest example of covering maps
$$\bigsqcup_{\forg_1(r,\Delta,s)=-\kappa_1}M_{\sigma_{\omega_1,\beta_1}}(r,\Delta,s)\rightarrow M_{\sigma}(-\kappa_1)^{\funct{U}}$$
where $\sigma$ is a $\funct{U}$-fixed stability condition on $\schk_X^1$ here and (on $\schk_X^2$ by abuse of notation) henceforth.

\begin{proposition}\label{GM-(1,-C,0)}
The object $\sho_C(-2)$ is $\sigma_{\omega_1,\beta_1}$-stable for any conic $C\subset S$.
	
\begin{proof}
At first, $\sho_C(-2)\in\scha_{\omega_1,\beta_1}$ and has Mukai vector $(0,C,-1)$. Suppose otherwise that the object $\sho_C(-2)$ is not $\sigma_{\omega_1,\beta_1}$-stable, then one can choose a short exact sequence
$$0\rightarrow A\rightarrow\sho_C(-2)\rightarrow B\rightarrow0$$
in $\scha_{\omega_1,\beta_1}$ such that $A$ is $\sigma_{\omega_1,\beta_1}$-stable with $\mu(A)\geq\mu(B)$. One has
$$\Im(Z_{\omega_1,\beta_1}(A))=(\frac{4}{5}r+\frac{1}{5}\Delta\cdot H)\geq0\quad\textup{and}\quad \Im(Z_{\omega_1,\beta_1}(B))=(\frac{2}{5}-\frac{4}{5}r-\frac{1}{5}\Delta\cdot H)\geq0$$
with $v(A)=(r,\Delta,s)$. Together with $\mu(A)\geq\mu(B)$, one gets $4r+\Delta\cdot H=0$ or $1$. 

Suppose that $4r+\Delta\cdot H=0$, then $\phi(A)=1$. According to \cite[Lemma 10.1]{Bri08}, the stable object $A$ is either isomorphic to $\sho_x$ for some point $x\in S$ or isomorphic to $\she[1]$ for some locally free sheaf $\she$. However, one can obtain from the long exact sequence
$$0\rightarrow\shh^{-1}(A)\rightarrow0\rightarrow\shh^{-1}(B)\rightarrow\shh^0(A)\rightarrow\sho_C(-2)\rightarrow\shh^0(B)\rightarrow0$$
that both cases are impossible. 

Suppose that $4r+\Delta\cdot H=1$, then the long exact sequence above will force $\shh^0(B)$ to be a torsion. It follows that $\shh^{-1}(B)=0$ or equivalently $r=0$. Otherwise, $\shh^{-1}(B)$ will be the torsion-free part of $\shh^0(A)$ which is impossible by construction. On the other hand, one can use \cite[Proposition 4.8]{Bo24} to see that $A$ and $B$ are spherical objects once $r=0$. It follows
$$\Delta^2-2rs=(C-\Delta)^2-2r(s+1)=-2$$ 
so that $C\cdot\Delta=-1-r=-1$. It is impossible as $\Delta$ should be effective in this case.
\end{proof}
\end{proposition}

On the other hand, the right hand side fixed locus is known.

\begin{theorem}[{\cite[Section 7.1]{JLLZ24}}]\label{GM-moduli-kappa1}
Consider a special Gushel--Mukai threefold $X\rightarrow W$ and a $\funct{U}$-fixed stability condition $\sigma$ on $\schk_X^1$, then any object in $M_{\sigma}(-\kappa_1)$ is $\sigma$-stable. Such a stable object is either isomorphic to a distinguished two-term complex $G_X$ or to a shift $\shi_C[1]$ of the ideal sheaf of a conic $C\subset X$ that is not the preimage of a line in $W$.\footnote{A more explicit schematic description of $M_{\sigma}(-\kappa_1)$ can be obtained by applying \cite[Theorem 7.3]{DK24}. However, that description will not be used in this paper so we choose not to elaborate on it.}
\end{theorem}

In particular, the fixed locus $M_{\sigma}(-\kappa_1)^{\funct{U}}$ consists precisely of the distinguished object $G_X$ and the shifts $\shi_C[1]$ with $C\subset S$. Similar to Corollary~\ref{QDS-mu1-covering}, one concludes that 

\begin{corollary}\label{GM-covering-w1}
The covering map for the class $-\kappa_1$ is decomposed as two coverings
$$M_{\sigma_{\omega_1,\beta_1}}(-2,H,-3)\sqcup M_{\sigma_{\omega_1,\beta_1}}(3,-H,2)\rightarrow\{G_X\}$$
and
$$\bigsqcup_{C\subset S\textup{ a conic}} M_{\sigma_{\omega_1,\beta_1}}(1,-C,0)\sqcup M_{\sigma_{\omega_1,\beta_1}}(0,C,-1)\rightarrow M_{\sigma}(-\kappa_1)^{\funct{U}}-\{G_X\}$$
where $M_{\sigma_{\omega_1,\beta_1}}(0,C,-1)=\{\sho_C(-2)\}$ and $M_{\sigma_{\omega_1,\beta_1}}(1,-C,0)=\{\sho_S(-C)\}$ for each conic $C\subset S$.

\begin{proof}
The only thing that needs to be checked is the isomorphism $$\funct{Forg}(\sho_C(-2))\cong\funct{L}_{\sho_X}\funct{L}_{\shu_X^\vee}(\sho_C)\cong\funct{L}_{\sho_X}\shi_C \cong \shi_C[1]$$
which is due to \cite[(7.2)]{PPZ26} and the fact that $\shi_C$ is in $\schk_X^1=\langle\sho_X,\shu_X^\vee\rangle^{\perp}\subset\cate{D}^b(X)$.
\end{proof}
\end{corollary}

\begin{example}
Let $Z\subset S$ be a length $2$ subscheme lying in a conic $C\subset S$, then the short exact sequence $0\rightarrow\sho_S(-C)\rightarrow\shi_Z\rightarrow\sho_C(-2)\rightarrow0$ implies that $\shi_Z$ is strictly $\sigma_{\omega_1,\beta_1}$-semistable.
\end{example}

\subsection{Over the vector $w_1$}\label{GM-ss-over-w1}
Here we are going to understand the covering morphism
$$M_{\sigma}(-\kappa_1)\rightarrow M_{\sigma_{\omega_1,\beta_1}}(w_1)^{\Pi_1}$$
for generic Brill--Noether general K3 surfaces of degree 10. At first, one observes that

\begin{proposition}\label{GM-(-2,H,-3)}
One has $M_{\sigma_{\omega_1,\beta_1}}(-2,H,-3)=\{\shu_S[1]\}$ and $M_{\sigma_{\omega_1,\beta_1}}(3,-H,2)=\{\shv_S\}$ where $\shu_S$ is the restricted tautological bundle on $S$ and $\shv_S\cong\funct{T}_{\sho_S}\shu^\vee_S[-1]$ is in $M_H(3,-H,2)$.
	
\begin{proof}
One checks that the stability condition $\sigma_{\omega_1,\beta_1}$ is $(-2,H,-3)$-generic and that the object $\shu_S[1]$ is $\sigma_{\omega_1,\beta_1}$-semistable. So the first assertion $M_{\sigma_{\omega_1,\beta_1}}(-2,H,-3)=\{\shu_S[1]\}$ follows. Then one has $M_{\sigma_{\omega_1,\beta_1}}(3,-H,2)=\{\shv_S\}$ by applying $\Pi_1$ and the fact that $\shu_S^\vee\cong\shu_S(H)$.

The remaining statement $M_{\sigma_{\omega_1,\beta_1}}(3,-H,2)=M_H(3,-H,2)$ follows from Theorem~\ref{K3-large-volume-limits} and the observation that $\sigma_{t\omega_1,\beta_1}$ is $(3,-H,2)$-generic for any $t\geq 1$.
\end{proof}
\end{proposition}

\begin{corollary}
Any strictly $\sigma_{\omega_1,\beta_1}$-semistable object with Mukai vector $w_1$ in $\scha_{\omega_1,\beta_1}$ is an extension of the objects $\shu_S[1]$ and $\shv_S$ in $\scha_{\omega_1,\beta_1}$. 

\begin{proof}
Choose a strictly $\sigma_{\omega_1,\beta_1}$-semistable object $E$ with Mukai vector $w_1$ in $\scha_{\omega_1,\beta_1}$, then one finds a short exact sequence $0\rightarrow A\rightarrow E\rightarrow B\rightarrow0$ in $\scha_{\omega_1,\beta_1}$ with $A$ stable and $\mu(A)=\mu(B)$. 

Set $v(A)=(r,aH,s)$, then one gets
$$\Im(Z_{\omega_1,\beta_1}(A))=\frac{4}{5}r+2a\geq0\quad\textup{and}\quad \Im(Z_{\omega_1,\beta_1}(B))=\frac{4}{5}-\frac{4}{5}r-2a\geq0$$
ensuring that $2r+5a=1$. Hence $\Im(Z_{\omega_1,\beta_1}(A))=\Im(Z_{\omega_1,\beta_1}(B))\neq0$ and then $\Re(Z_{\omega_1,\beta_1}(A))=\Re(Z_{\omega_1,\beta_1}(B))$, from which one obtains $3r+20a+5s=-1$. Then $s-r=-1$ and $r+s+5a=0$, so by Proposition~\ref{GM-image-forg} one has $\forg_1(r,aH,s)=-\kappa_1$. One concludes by Corollary~\ref{GM-covering-w1}.
\end{proof}
\end{corollary}

Then we come to figure out the stable objects in $M_{\sigma_{\omega_1,\beta_1}}(w_1)$. The spirit is the following assertion, analogous to Proposition~\ref{QDS-u_1-wall-crossing}.

\begin{proposition}
The stability condition $\sigma_{t\omega_1,\beta_1}$ is $w_1$-generic for any $t> 1$. 
	
\begin{proof}
To simplify the notations, we will work on stability conditions $\sigma_{\omega_1(t),\beta_1}$ for $t\geq0$ where the ample class $\omega_1(t)$ is $\sqrt{5t+1}\omega_1$. The central charge of $\sigma_{\omega_1(t),\beta_1}$ is
$$Z_t(r,aH,s)=((t-\frac{3}{5})r-4a-s)+(\frac{4}{5}r+2a)\sqrt{5t+1}\sqrt{-1}$$
		
Since $w_1$ is primitive, it suffices to check that there are no strictly $\sigma_{t\omega_1,\beta_1}$-semistable objects with Mukai vector $w_1$. Suppose otherwise that such an object $E$ exists, then up to a shift one finds an exact sequence $0\rightarrow A\rightarrow E\rightarrow B\rightarrow0$ in $\scha_{\omega_1,\beta_1}$ with $A$ stable and $\mu(A)=\mu(B)$. 
		
Set $v(A)=(r,aH,s)$, then one gets
$$\Im(Z_t(A))=(\frac{4}{5}r+2a)\sqrt{5t+1}\geq0\quad\textup{and}\quad \Im(Z_t(B))=(\frac{4}{5}-\frac{4}{5}r-2a)\sqrt{5t+1}\geq0$$
ensuring that $2r+5a=1$. So $\Im(Z_t(A))=\Im(Z_t(B))\neq0$ and then $\Re(Z_t(A))=\Re(Z_t(B))$, from which one obtains $(10t-6)r-40a-10s=5t+2$. Therefore, one sees
$$r=\frac{1}{2}-\frac{5}{2}a\quad\textup{and}\quad s=-\frac{1}{2}-\frac{5}{2}a-\frac{5t}{2}a$$
		
Since $A$ is stable, one has $10a^2-2rs\geq -2$. The equations above give
$$ 10a^2-2(\frac{5}{2}a-\frac{1}{2})(\frac{1}{2}+\frac{5}{2}a+\frac{5t}{2}a)\geq -2$$
or equivalently $ ta+1\geq 5ta^2+a^2$. It is always true when $a=0$, but then $r\notin\bz$. Also, one can see directly that $a\leq -1$ is impossible. Suppose that $a\geq 1$, then 
$$1\geq a^2+t(5a^2-a)\geq a^2+t$$
which implies that there are no valid $t>0$.
\end{proof}
\end{proposition}

In conclusion, one has a contraction $S^{[2]}\rightarrow M_{\sigma_{\omega_1,\beta_1}}(w_1)$ which sends a closed subset to the unique semistable point. According to \cite[Theorem 5.7]{BM14-MMP}, the contracted locus has codimension at least $2$. So the involution $\pi_{1}$ of $M_{\sigma_{\omega_1,\beta_1}}(w_1)$ given in Proposition~\ref{GM-induced-involution} can be seen as a birational involution on $S^{[2]}$. Since $\Bir(S^{[2]})=\bz/2\bz$ (see \cite[Appendix B]{DM19}), this birational involution is exactly the one constructed in \cite[Section 4.3]{O'Grady05} (see also \cite[Example 4.16]{De22}). It also means that $M^{\stable}_{\sigma_{\omega_1,\beta_1}}(w_1)$ is contained in the regular locus $\{\shi_Z\,|\,Z\not\subset L\subset W\}$ of the birational involution.

The fixed locus $\Fix(\pi_{1})\subset M_{\sigma_{\omega_1,\beta_1}}(w_1)$ is isomorphic to the union of a distinguished point and the surface $BTC_S$ of bitangent conics. Moreover, the covering morphism 
$$M_{\sigma}(-\kappa_1)\rightarrow \Fix(\pi_{1})$$ 
sends the point $G_X$ to the distinguished point in the fixed locus. Moreover, a stable object $\shi_C[1]$ is sent to the unique $\shi_Z$ such that $2Z=C\cap S$ by an argument similar to that for Proposition~\ref{QDS-u1-covering-map}.

\subsection{Over the vector $w_2$}\label{GM-ss-over-w2}
Here we are going to understand the covering morphism
$$M_{\sigma}(-\kappa_2)\rightarrow M_{\sigma_{\omega_1,\beta_1}}(w_2)^{\Pi_1}$$
for generic Brill--Noether general K3 surfaces of degree 10.

\begin{proposition}\label{GM-w2-generic}
The stability condition $\sigma_{t\omega_2,\beta_2}$ is $w'_2$-generic for any $t\geq 1$.
	
\begin{proof}
As before, suppose that there exists a strictly semistable object $E\in\scha_{\omega_2,\beta_2}$ together with a short exact sequence $0\rightarrow A\rightarrow E\rightarrow B\rightarrow0$ in $\schp(\phi(E))$ such that $A$ is $\sigma_{t\omega_2,\beta_2}$-stable.

Set $v(A)=(r,aH,s)$, then one gets $1\geq 3r+5a\geq0$ from $\Im Z(A)\geq0$ and $\Im Z(B)\geq0$. It is impossible for integers $r$ and $a$, because $\mu(A)=\mu(B)$.
\end{proof}
\end{proposition}

Therefore one can use Theorem~\ref{K3-large-volume-limits} and Proposition~\ref{GM-rotation-identification} to see that
$$ M_{\sigma_{\omega_1,\beta_1}}(w_2)\cong M_{\sigma_{\omega_2,\beta_2}}(w'_2)\cong M_H(w_2)$$
is exactly the irreducible symplectic fourfold introduced in \cite[Example 5.17]{Mu88} and described with more details in \cite[Section 5.4]{O'Grady05}. Moreover, the fourfold $M_{\sigma_{\omega_1,\beta_1}}(w_2)$ is isomorphic to the double EPW sextic $\tilde{Y}_A$ for a Lagrangian subspace $A\subset\land^3H^0(X,\shi_S(2))$ according to \cite{O'Grady06}. The induced non-symplectic involution $\pi_{1}$ of $M_{\sigma_{\omega_1,\beta_1}}(w_2)$ is the covering involution and the fixed locus, called the EPW surface, is a smooth surface of general type and shares some invariants with the surface of bitangent lines to a generic quartic K3 according to \cite{Bea11}. 

\begin{remark}
A priori, the Lagrangian subspace $A$ constructed in \cite[Theorem 5.2]{O'Grady06} might not coincide with the one $A(S)$ given in \cite[Theorem 3.10]{DK18}. But the EPW sextic $Y_{A(S)}$ coincides with the EPW sextic $Y_A$ according to \cite[Section 5.4]{O'Grady05} and \cite[Corollary 3.19]{DK18}. Therefore the two Lagrangian subspaces give the same double EPW sextic.
\end{remark}

The description of $\pi_1$ can be found in \cite[Section 4.2.2]{O'Grady05} and \cite[Theorem 3.21 (1)]{Ma02} but here for simplicity we will only define it for generic points. Such a point is represented by a locally free sheaf $\she$, which fits into a short exact sequence $0\rightarrow\she\rightarrow\sho_S^{\oplus4}\rightarrow\shg\rightarrow0$ such that $\shg$ is a locally free sheaf in $M_H(2,H,2)$ and $\sho_S^{\oplus 4}\rightarrow\shg$ is the evaluation map. Then the dual locally free sheaf $\shg^\vee$ is in $M_H(w_2)=M_{\sigma_{\omega_1,\beta_1}}(w_2)$ and the image $\pi_1([\she])$ is precisely $[\shg^\vee]$.

\begin{proposition}[{{\cite[Theorem 8.9]{JLLZ24}}}]
Consider a special Gushel--Mukai threefold $X$ and a $\funct{U}$-fixed stability condition $\sigma$ on $\schk_X^1$, then $M_{\sigma}(-\kappa_2)$ is isomorphic to the moduli space $M_X(2,1,5)$ of Gieseker semistable sheaves on $X$ with $c_0=2,c_1=h,c_2=h^2/2,c_3=0$ via the projection functor $\funct{L}_{\sho_X}\funct{L}_{\shu_X^\vee}\colon\cate{D}^b(X)\rightarrow\schk_X^1$. More explicitly, any object in $M_{\sigma}(-\kappa_2)$ is the image of a unique object in $M_X(2,1,5)$ up to shift via the projection functor.
\end{proposition}

The moduli space $M_{\sigma}(\kappa_2)$ has no fixed objects by \cite[Proposition 8.1]{JLLZ24}. So the induced double covering $M_{\sigma}(-\kappa_2)\rightarrow M_{\sigma_{\omega_1,\beta_1}}(w_2)^{\Pi_1}$ is étale everywhere. Here we describe this étale covering on the dense subset $U\subset M_X(2,1,5)\cong M_{\sigma}(-\kappa_2)$ consisting of locally free sheaves.

\begin{proposition}
A locally free sheaf $\shf$ in $M_X(2,1,5)$ is sent to the restriction $j^*\shh$ on the ramification divisor $j\colon S\subset X$ via $M_X(2,1,5)\cong M_{\sigma}(-\kappa_2)\rightarrow M_{\sigma_{\omega_1,\beta_1}}(w_2)^{\Pi_1}$, where the locally free sheaf $\shh\cong\funct{L}_{\sho_X}\funct{L}_{\shu_X^\vee}\shf[-1]\cong\funct{L}_{\sho_X}\shf[-1]$ is an element in $M_X(2,-1,5)$. 

\begin{proof}
Since we are working on a generic special Gushel--Mukai threefold $X\rightarrow W$, the divisor $S$ is a generic element in $|\sho_X(1)|$. So the locally free sheaf $j^*\shh$ on $S$ is again Gieseker semistable due to Maruyama's theorem \cite{Ma80} and hence an element in $M_{\sigma_{\omega_1,\beta_1}}(w_2)=M_H(w_2)$. 

Then we compute $\funct{Forg}(j^*\shh)$. One notices that
\begin{align*}
\Hom(\funct{Forg}(j^*\shh),\shh[1])&=\Hom(j_*j^*\shh(1),\shh[1])\\
&\cong\Hom(\shh,j_*j^*\shh[2])\\
&\cong\Hom(j^*\shh,j^*\shh[2])\cong\bc
\end{align*}
so there exists a non-trivial morphism $\funct{Forg}(j^*\shh)\rightarrow \shh[1]$. It follows that $\funct{Forg}(j^*\shh)$ is a strictly semistable object and can only be $\shh[1]\oplus\funct{U}(\shh)[1]$ as it is fixed by $\funct{U}$. One hence concludes $j^*\shh=\funct{Inf}(\shh[1])$ by  $\Hom(\funct{Inf}(\shh[1]),j^*\shh)=\Hom(\shh[1],\shh[1]\oplus\funct{U}\shh[1])=\bc$.
\end{proof}
\end{proposition}

\begin{remark}
Since the open subset $U\subset M_X(2,1,5)\cong M_{\sigma}(-\kappa_2)$ is dense, the description actually determines the étale covering $M_{\sigma}(-\kappa_2)\rightarrow M_{\sigma_{\omega_1,\beta_1}}(w_2)^{\Pi_1}$. It is expected that this étale double covering coincides with the one in \cite[Theorem 5.2]{DK20}.
\end{remark}

\subsection{Over the vector $w_1+w_2$}\label{GM-ss-over-w1plusw2}
Here we are trying to understand the induced involution
$$\pi_1\colon M_{\sigma_{\omega_1,\beta_1}}(w_1+w_2)\rightarrow M_{\sigma_{\omega_1,\beta_1}}(w_1+w_2)$$
for generic Brill--Noether general K3 surfaces of degree 10. The moduli space $M_{\sigma}(-\kappa_1-\kappa_2)$ is studied in 
\cite{Zhang20} and \cite[Appendix B]{FGLZ25} but its structure is not totally clear. One can first see that $M_{\sigma}(-\kappa_1-\kappa_2)$ has only stable objects and by Proposition~\ref{GM-image-forg} that $M_{\sigma}(-\kappa_1-\kappa_2)^{\funct{U}}=\varnothing$.

\begin{proposition}
One has $M_{\sigma_{t\omega_1,\beta_1}}(w_1+w_2)\cong M_H(w_1+w_2)$ for any $t\geq1$.

\begin{proof}
It suffices to show that the stability condition $\sigma_{t\omega_1,\beta_1}$ is $(w_1+w_2)$-generic for $t\geq 1$, and the argument is similar to that of Proposition~\ref{GM-w2-generic}.
\end{proof}
\end{proposition}

In particular, $M_{\sigma_{\omega_1,\beta_1}}(w_1+w_2)$ is an irreducible symplectic manifold. Moreover, using the identification in Proposition~\ref{GM-rotation-identification}, one can see 

\begin{proposition}
The variety $M_{\sigma_{\omega_1,\beta_1}}(w_1+w_2)=M_{\sigma_{\omega_2,\beta_2}}(w'_1+w'_2)$ is birational to $S^{[3]}$. 

\begin{proof}
To simplify the notations, we will work on stability conditions $\sigma_{\omega_2(t),\beta_2}$ for $t\geq0$ where the ample class $\omega_2(t)$ is defined to be $\sqrt{5t+1}\omega_2$. Let $Z_t$ be its central charge, then
$$Z_t(r,aH,s)=((t-\frac{8}{5})r-6a-s)+(\frac{6}{5}r+2a)\sqrt{5t+1}\sqrt{-1}$$

Choose a potential strictly semistable object $E$ in $\scha_{\omega_2,\beta_2}$ with Mukai vector $-w'_1-w'_2$ and a short exact sequence $0\rightarrow A\rightarrow E\rightarrow B\rightarrow0$ in $\schp(\phi(E))$ such that $A$ is stable. Then
$$\Im(Z_t(A))=(\frac{6}{5}r+2a)\sqrt{5t+1}\geq0\quad\textup{and}\quad\Im(Z_t(B))=(\frac{6}{5}-\frac{6}{5}r-2a)\sqrt{5t+1}\geq0$$
which implies that $3r+5a=1$ or $3r+5a=2$.

Suppose that $3r+5a=1$, then $2\Re(Z_t(A))=\Re(Z_t(B))$ or equivalently $$24r+90a+15s+2+5t-15tr=0$$
which allows one to obtain
$$r=-\frac{5}{3}a+\frac{1}{3}\quad\textup{and}\quad s=-\frac{2}{3}-\frac{10}{3}a-\frac{5t}{3}a$$
Since $A$ is stable, one has $10a^2-2rs\geq -2$. It follows that $t(5a-25a^2)\geq5a^2-11$. One can rule out the solution $a=0$ as $r\in\bz$. The left hand side is negative for any $a\neq0$ and $t>0$, but the right hand side is negative only when $a^2=1$. Since $r\in\bz$, one must have $a=-1$. Then 
$$r=2,\quad s=\frac{8+5t}{3},\quad -30t\geq -6$$ 
and one can see that the only possibility is $s=3,t=1/5$.

Suppose that $3r+5a=2$, then $\Re(Z_t(A))=2\Re(Z_t(B))$ or equivalently $$24r+90a+15s+4+10t-15tr=0$$
which allows one to obtain
$$r=-\frac{5}{3}a+\frac{2}{3}\quad\textup{and}\quad s=-\frac{4}{3}-\frac{10}{3}a-\frac{5t}{3}a$$
Since $A$ is stable, one has $10a^2-2rs\geq -2$. It follows that $t(10a-25a^2)\geq5a^2-17$. One can rule out the solution $a=0$ as $r\in\bz$. The left hand side is negative for any $a\neq0$, but the right hand side is negative only when $a^2=1$. Since $r\in\bz$, one must have $a=1$. Then 
$$r=-1,\quad s=-\frac{14+5t}{3},\quad -15t\geq -12$$ 
and one can see that the only possibilities are $s=-5,t=1/5$ and $s=-6,t=4/5$.

In the latter case, one has an exact sequence $0\rightarrow A\rightarrow E\rightarrow B\rightarrow0$ of semistable objects such that $v(B)=(2,-H,4)$. It is not hard to check as before that $\sigma_{\omega_2(t),\beta_2}$ is $(2,-H,4)$-generic for any $t\geq0$. So the moduli space $M_{\sigma_{\omega_2(t),\beta_2}}(2,-H,4)$ has to be empty.

In conclusion, $\sigma_{\omega_2(t),\beta_2}$ is $(-w'_1-w'_2)$-generic for $t\in[0,+\infty)\setminus\{1/5\}$ and the only wall is a flopping wall induced by Mukai vectors $(2,-H,3)$ and $(-1,H,-5)$. Therefore, one obtains the desired birational map $M_{\sigma_{\omega_2,\beta_2}}(-w'_1-w'_2)\dashrightarrow M_H(-w'_1-w'_2)=S^{[3]}$.
\end{proof}
\end{proposition}

The last thing to mention is that the irreducible symplectic manifold $M_{\sigma_{\omega_1,\beta_1}}(w_1+w_2)$ is nothing but the double EPW cube $\tilde{K}_A$. In fact, one has a Hodge isometry
$$H^2_{prim}(\tilde{Y}_A,\bz)\cong H^2_{prim}(M_H(w_1+w_2),\bz)\cong \{w_1,w_2\}^{\perp}\subset\tilde{H}(S,\bz)$$
where the primitive cohomology structures are given by the ample class which generates the sublattice $\bz[\theta_{w_2}(w_1)]$ for $\tilde{Y}_A=M_H(w_2)$ and the ample class generating $\bz[\theta_{w_1+w_2}(w_1-w_2)]$ for $M_H(w_1+w_2)$. The existence of these ample classes is due to the following observation

\begin{proposition}
The class $\theta_{w_2}(w_1)$ is ample on $\tilde{Y}_A=M_H(w_2)$ up to a minus sign and the class $\theta_{w_1+w_2}(w_1-w_2)$ is ample on $M_H(w_1+w_2)$ up to a minus sign. 

\begin{proof}
Since the fixed lattice of $\Pi_1^*$ on $\tilde{H}(S,\bz)$ is spanned by $w_1$ and $w_2$, the fixed lattice of the Hodge isometry $\pi_1^*$ induced by the involution $\pi_1$ on $H^2(\tilde{Y}_A,\bz)$ is generated by $\theta_{w_2}(w_1)$ due to Proposition~\ref{K3-invariant-lattices-compatible}. It follows that either $\theta_{w_2}(w_1)$ or $-\theta_{w_2}(w_1)$ is an ample class. The argument for one of the classes $\pm\theta_{w_1+w_2}(w_1-w_2)$ being ample on $M_H(w_1+w_2)$ is identical. 
\end{proof}
\end{proposition}

So one has a birational map $M_{\sigma_{\omega_1,\beta_1}}(w_1+w_2)\dashrightarrow \tilde{K}_A$ according to \cite[Theorem 1.1]{KKM24} and the global Torelli theorem. On the other hand, there exists a unique non-trivial birational involution $\pi$ on $S^{[3]}$ according to \cite{BC22} and, by virtue of \cite[Example 4.12]{De22}, $S^{[3]}$ has a unique non-trivial irreducible symplectic birational model $M$ on which $\pi$ is regular. 

\begin{corollary}
One has $M_{\sigma_{\omega_1,\beta_1}}(w_1+w_2)\cong \tilde{K}_A$, and the involution on it induced by $\Pi_1$ induces the unique element in $\Bir(S^{[3]})$ via the flop.

\begin{proof}
The involution on $M_{\sigma_{\omega_1,\beta_1}}(w_1+w_2)$ induced by $\Pi_1$ will give a non-trivial birational involution on $S^{[3]}$, and the canonical involution on $\tilde{K}_A$ constructed in \cite[Section 3]{IKKR19} also induces a non-trivial birational involution on $S^{[3]}$ according to \cite[Section 5]{IKKR19}. Then one can conclude using uniqueness of the birational model declared before this corollary.
\end{proof}
\end{corollary}

{\appendix	

\section{Enriques categories and two types of Hodge isometries}\label{ST-Hodge-Involution}
\subsection{The general considerations}
In general, suppose that one has an involutive Hodge isometry $\tau$ on $\tilde{H}(S,\bz)$ which fixes an object in $\Stab^{\circ}(S)$ and let $\Pi$ be the involution given by Proposition~\ref{Appendix-Huybrechts-conway}. Once $\Pi$ fixes a stable object, one can apply \cite[Corollary 4.11]{BO23} to conclude that $\Pi$ induces a $\bz/2\bz$ action on $\cate{D}^b(S)$. Then we can use the following observation, together with \cite[Lemma 6.5]{BO23}, to decide whether the equivariant category $\cate{D}^b(S)^{\Pi}$ is Enriques.

\label{A_definition-Enriques}
\begin{proposition}
Consider a $\bz/2\bz$-action on $\cate{D}^b(S)$ for a K3 surface $S$ generated by a non-symplectic involution $\Pi$, then the equivariant category $\sche=\cate{D}^b(S)^{\Pi}$ admits a Serre functor isomorphic to the involution $\funct{U}[2]$, where $\funct{U}$ is a generator of the residual $\bz/2\bz$-action on $\sche$.
	
\begin{proof}
According to \cite[Proposition 5.2]{BO23}, the Serre functor $\funct{S}_{\sche}$ exists and is given by 
$$(A,\phi)\mapsto(\funct{S}A,\hat{\phi})\quad\textup{and}\quad f\mapsto\funct{S}(f)$$
on objects and morphisms, such that the new linearization $(\phi')_g$ is defined by the morphisms
$$\hat{\phi}_1\colon \funct{S}A\stackrel{\funct{S}\phi_1}{\longrightarrow}\funct{S}A\quad\textup{and}\quad \hat{\phi}_{-1}\colon \funct{S}A\stackrel{\funct{S}\phi_{-1}}{\longrightarrow}\funct{S}\Pi A \stackrel{t_A}{\longrightarrow} \Pi\funct{S} A$$
where $t\colon \funct{S}\Pi\stackrel{\sim}{\rightarrow}\Pi\funct{S}$ is the canonical isomorphism described in \cite[Lemma 5.1]{BO23}. 
		
Since $\Pi$ is a non-symplectic involution, it sends a generator $\alpha=(\funct{Id}_S\stackrel{\sim}{\rightarrow}\funct{S}[-2])$ in the Hochschild homology group $\Hoch_2(S)=H^{2,0}(S)$ to $-\alpha$. In other words, one has $\Pi\alpha\Pi^{-1}=-\alpha$ or equivalently the following commutative diagram
\begin{displaymath}
			\xymatrix{
				\Pi\circ\cate{Id}\ar[rr]^{\id}\ar[d]_{\Pi\alpha}&&\cate{Id}\circ\Pi\ar[d]_{-\alpha\Pi}\\
				\Pi\funct{S}[-2]\ar[rr]^{t^{-1}[-2]}&&\funct{S}\Pi[-2]
			}
\end{displaymath}
of natural transformations, where the canonical commutativity with $[2]$ is implicitly used. 
		
One can thus get a commutative diagram
\begin{displaymath}
			\xymatrix{
				A\ar[rr]^{\phi_{-1}[2]}\ar[d]_{\alpha_A}&&\Pi A\ar[rr]^{-\id_{\Pi A}}\ar[d]_{(\alpha\Pi)_A}&&\Pi A\ar[d]_{(\Pi\alpha)_A}\\
				\funct{S}A[-2]\ar[rr]^{\funct{S}\phi_{-1}[-2]}&&\funct{S}\Pi A[-2]\ar[rr]^{t_A[-2]}&&\Pi \funct{S} A[-2]
			}
\end{displaymath}
for any $(A,\phi)\in \sche$, where the canonical commutativity with $[2]$ is implicitly used again. 
		
Moreover, one automatically has the following commutative diagram
\begin{displaymath}
			\xymatrix{
				A\ar[rr]^{\phi_1}\ar[d]_{\alpha_A}&& A\ar[rr]^{\id_{A}}\ar[d]_{\alpha_A}&& A\ar[d]_{\alpha_A}\\
				\funct{S}A [-2]\ar[rr]^{\funct{S}[-2]\phi_1}&&\funct{S} A [-2]\ar[rr]^{\id_{\funct{S}A [-2]}}&&\funct{S} A [-2]
			}
\end{displaymath}
for any $(A,\phi)\in \sche$. According to \cite[Section 5]{BO22}, the functor $\funct{S}_{\sche}[-2]$ acts on morphisms through the functor $\funct{S}[-2]$ on $\cate{D}^b(S)$, which is canonically trivial. Since the functor $\funct{U}$ acts also trivially on morphisms, the above two commutative diagrams define an isomorphism $\funct{U}\stackrel{\sim}{\rightarrow}\funct{S}_{\sche}[-2]$.
\end{proof}
\end{proposition}

\begin{remark}
One should notice that this argument also works in a more general setting.
\end{remark}

The next natural question is how we can find a non-symplectic involution on $\cate{D}^b(S)$. Once such an involution on $\cate{D}^b(S)$ is given, one can use the criterion \cite[Lemma 3.7]{BO23} to see whether it will act on $\cate{D}^b(S)$ and then give an Enriques category over $S$.

To find non-symplectic involutions, one needs to look more closely at $\Aut(\cate{D}^b(S))$. Let us define
$$\Aut^{\circ}(\cate{D}^b(S))\subset\Aut(\cate{D}^b(S))$$
to be the subgroup of all autoequivalences fixing $\Stab^{\circ}(S)$. One has a natural surjection
$$\hbar\colon\Aut^{\circ}(\cate{D}^b(S))\rightarrow\Aut(\tilde{H}(S,\bz))_+$$
where the right hand side is the group of Hodge isometries preserving the (natural) orientation of positive planes $P\subset\tilde{H}^{1,1}(S,\br)$. Moreover, one has further a subgroup
$$\Aut(\cate{D}^b(S),\sigma)\subset\Aut^{\circ}(\cate{D}^b(S))$$
which consists of autoequivalences that fix a given stability condition $\sigma\in\Stab^{\circ}(S)$.

\begin{proposition}[{{\cite[Proposition 1.4]{Hu16}}}]\label{Appendix-Huybrechts-conway}
Consider a stability condition $\sigma=(\scha,Z)\in\Stab^{\circ}(S)$, then the surjective homomorphism $\hbar$ induces an isomorphism
$$\Aut(\cate{D}^b(S),\sigma)\stackrel{\sim}{\rightarrow}\Aut(\tilde{H}(S,\bz),Z)$$
where the right hand side group consists of Hodge isometries fixing the stability function $Z$.
\end{proposition}

In this case, the first step is to find involutive Hodge isometries on $\tilde{H}(S,\bz)$. The current examples of Enriques categories over K3 surfaces are mostly due to \cite{KP17}. 

The involution on $\cate{D}^b(S)$ has two forms
$$\Pi=\funct{O}^n[-1]\quad\textup{and}\quad \Pi=\funct{T}_U\circ\funct{O}[-1]$$
where $\funct{O}:=\funct{T}_{\sho_S}(-\otimes\sho_S(H))$ and $U$ is a spherical object on $\cate{D}^b(S)$, corresponding respectively to the scenarios that $\cate{D}^b(X)=\langle\cate{D}^b(S)^{\Pi},\sho_X,\cdots,\sho_X(n)\rangle$ and $\cate{D}^b(X)=\langle\cate{D}^b(S)^{\Pi},\sho_X,U_X\rangle$.

In the following, we study the Hodge isometries induced by $\funct{O}^n[-1]$ and $\funct{T}_U\circ\funct{O}[-1]$ and find all the possible induced involutions among them. Here is a summary of the results.

\begin{proposition}\label{A_involution-type-1}
Consider a polarized K3 surface $(S,H)$ of degree $2d$ and the Hodge isometry $\tau_d$ on $\tilde{H}(S,\bz)$ induced by $\funct{O}_d=\funct{T}_{\sho_S}(-\otimes\sho_S(H))$, then $-\tau_d^n$ is an involution for some $n\neq0$ if and only if the pair $(d,n)$ is an element in $\{(1,3k),(2,2k),(3,3k)\,|\,k\in\bz,k\neq0\}$
\end{proposition}

The solutions $(1,3)$ and $(2,2)$ correspond to existing examples, while the solution $(3,3)$ does not correspond to an involution on the category $\cate{D}^b(S)$ in the generic case.

\begin{proposition}\label{A_involution-type-2}
Consider a polarized K3 surface $(S,H)$ of degree $2d$, a spherical object $U$ of $\cate{D}^b(S)$ with Mukai vector $(a,bH,c)$ for some $a>0$, and the Hodge isometry $\tau_U$ on $\tilde{H}(S,\bz)$ induced by the autoequivalence $\funct{T}_U\circ\funct{T}_{\sho_S}\circ(-\otimes\sho_S(H))[-1]$ of $\cate{D}^b(S)$, then $\tau_U$ is an involution if and only if $d=5$ and the tuple $(a,b,c)$ is an element in the set
$$\{(F_n+F_{n-2},-F_n,F_{n+2}+F_{n}),(F_{n+2}+F_n,-F_n,F_n+F_{n-2})\,|\,n\geq 0\textup{ even}\}$$
where $\{F_n\}$ is a sequence defined by $F_n=1$ for $n\leq 1$ and $F_{n+2}=F_{n+1}+F_n$ for $n\geq0$.	
\end{proposition}

The tuple $(2,-1,3)$ corresponds to the involution $\Pi_2$ and the Enriques category $\schk_X^2$, while the tuple $(3,-1,2)$ corresponds to $\Pi_1$ and $\schk_X^1$. Moreover, each tuple $(a,b,c)$ gives rise to an involution on the derived category $\cate{D}^b(S)$ of a generic K3 surface of degree $10$.

\subsection{The power of rotation functors}
Consider a polarized K3 surface $(S,H)$ with $H$ ample and $H^2=2d$, one can define the functor $\funct{O}_d:=\funct{T}_{\sho_S}\circ(-\otimes\sho_S(H))$ and consider the induced Hodge isometry $\tau_d$. In this subsection, we will work out all the $(d,n)$ such that $\tau^{2n}_d=\id$.

\begin{remark}
This is motivated by the fact that most examples of Enriques categories come from \cite{KP17}. In many cases, the dual group action on $\cate{D}^b(S)$ is generated by an autoequivalence with Hodge isometry $-\tau_d^n$ according to \cite[Theorem 7.7]{KP17}
\end{remark}

The Hodge isometry $\tau_d$ on $\tilde{H}(S,\bz)$ induced by $\funct{O}_d$ can be represented by a matrix
$$O_d=\begin{pmatrix}
	-d&-H&-1\\H&1&0\\ -1&0&0
\end{pmatrix}$$
and the power $\tau_d^{2n}$ is the identity if and only if $O_d^{2n}$ is equal to the identity matrix $I_3$. One notices that the determinant of $O_d$ is $-1$, so it is reasonable to only consider the even powers.

\begin{example}\label{example-R-power-2}
One computes that $O^2_d$ is the matrix 
$$\begin{pmatrix}
		(d-1)^2&(d-1)H&d\\
		-(d-1)H&-H* H+1&-H\\ 
		d&H&1
\end{pmatrix}$$
where $H* H$ is the operator sending $\Delta\in H^2(S,\bz)$ to $(H\cdot \Delta)H$. It cannot be $I_3$ for any $d$.
\end{example}

\begin{example}\label{example-R-power-4}
One computes that $O^4_d$ is the matrix 
$$\begin{pmatrix}
		(d^2-3d+1)^2&(d-2)(d^2-3d+1)H&(d-2)^2d\\
		-(d-2)(d^2-3d+1)H&-(d-2)^2H* H+1&-(d-1)(d-2)H\\ 
		(d-2)^2d&(d-1)(d-2)H&(d-1)^2
\end{pmatrix}$$
which equals the identity matrix only when $d=2$.
\end{example}

To make the computation easier, we instead consider the following matrix of functions
$$O(x)=\begin{pmatrix}
	-x&-\sqrt{2x}&-1\\\sqrt{2x}&1&0\\ -1&0&0
\end{pmatrix}$$
then $O_d^{2n}=I_3$ only if $O^{2n}(d)=I_3$. Indeed, it is an equivalent condition. 

\begin{proposition}\label{example-R-power-6-degree-1}
Suppose that $O^{2n}(d)=I_3$ for some $n\in\bz$, then $O_d^{2n}= I_3$.
	
\begin{proof}
Since $O^{2n}(d)=I_3$, the matrix $O_d^n$ must be of the form
$$A=\begin{pmatrix}
			1&0&0\\0&a(H*H)+b&0\\0&0&1
\end{pmatrix}$$
where $a,b$ are integers such that $2ad+b=1$ and $|b|=1$. Once $a=0$, one has $A=I_3$. 

Otherwise, if $a\neq0$, the only possible case is $(a,b,d)=(1,-1,1)$, and one can check in this case that $A^2=I_3$ and $O_1^3=-A$. Set $n=3k+m$ for the residue $m\in\{0,1,2\}$, then due to $O_1^6=I_3$ one has $O_1^{2m}=I_3$. However, it is impossible for $m=1,2$ by plugging $d=1$ in Example~\ref{example-R-power-2} and Example~\ref{example-R-power-4}. Then $n=3k$, which is clearly impossible. So we are done.
\end{proof}
\end{proposition}

One computes that the characteristic polynomial of the matrix $O(x)$ is 
$$(\lambda+1)(\lambda^2+(x-2)\lambda+1)$$
and the eigenvalues of $O(x)$ are $-1$ and $\frac{1}{2}(2-x\pm\sqrt{x^2-4x})$. One notices for $x>4$ that
$$0>\frac{1}{2}(2-x+\sqrt{x^2-4x})>\frac{1}{2}(2-x-\sqrt{x^2-4x})$$
which follows immediately that $O^n(x)$ cannot be the identity matrix for $x>4$.

Then it remains to deal with the cases $d=1,2,3,4$. At first, one observes that $O^n(4)$ can never equal $I_3$ for any $n>0$, by the Jordan normal form of $O(4)$. Also, one computes that $O^6(3)=I_3$. Combining with Example~\ref{example-R-power-2}, Example~\ref{example-R-power-4}, and Proposition~\ref{example-R-power-6-degree-1}, we conclude 

\begin{theorem}
The only integer pairs $(d,n)$ with $d>0$ and $n\neq0$ such that $O_d^{2n}=I_3$ are $(1,3k),(2,2k)$ and $(3,3k)$ for any non-zero $k\in\bz$.
\end{theorem}

In other words, all the possible non-trivial Hodge involutions on $\tilde{H}(S,\bz)$ one can get of the form $\tau_d^n$ are exactly $\pm\tau_1^3,\pm\tau_2^2,\pm\tau_3^3$.

\begin{example}
Let $(S,H)$ be a polarized K3 surface of degree $2$, then $-\tau_1^3$ is the map
$$(r,\Delta,s)\mapsto(r,(\Delta\cdot H)H-\Delta,s)$$
One checks that $-\tau_1^3$ fixes the stability functions $Z_{aH,bH}$, so by Proposition~\ref{Appendix-Huybrechts-conway} there is an involution $\Pi$ realizing $-\tau_1^3$, unique up to isomorphism. In fact, such a K3 surface admits a double cover $S\rightarrow\bp^2$ branched over a smooth sextic and the involution $\Pi$ is isomorphic to the pushforward along the covering involution of this branched cover.
\end{example}

\begin{remark}
This example is described in \cite[Example 3.12]{PPZ26}. One notices that the involution in \cite[Example 3.11]{PPZ26} can be seen as a variant of this example.
\end{remark}

\begin{example}
Let $(S,H)$ be a polarized K3 surface of degree $4$, then $-\tau_2^2$ is the map
$$(r,\Delta,s)\mapsto(-r-2s-\Delta\cdot H,(r+s+\Delta\cdot H)H-\Delta,-2r-s-\Delta\cdot H)$$
According to Proposition~\ref{QDS-fixed-stability-condition}, $-\tau_2^2$ only preserves the function $Z_{\omega_q,\beta_q}$. Also, one checks that the pair $\sigma_{\omega_q,\beta_q}=(\schp_{\omega_q,\beta_q},Z_{\omega_q,\beta_q})$ is a stability condition if and only if one cannot pick a class $\Delta\in\NS(S)$ with $\Delta\cdot H=-2$ and $\Delta^2=0$. It implies that $H$ is very ample by Reider's Method. So the unique involution realizing $-\tau_2^2$ is isomorphic to $\Pi_q$ due to Proposition~\ref{Appendix-Huybrechts-conway}.
\end{example}

One checks that the remaining Hodge involutions $\pm\tau_3^3$ do not fix any geometric stability function on $\cate{D}^b(S)$. Therefore, provided that the polarized K3 surface $S$ has Picard rank $1$, there are no involutions of $\cate{D}^b(S)$ realizing them by virtue of \cite[Theorem 4.13]{FL23}.

\subsection{The rotation functor with a spherical twist}
Consider a polarized K3 surface $(S,H)$ of degree $2d$ and a spherical object $U$ in $\cate{D}^b(S)$ such that $c_1(U)\in\bz[H]$. We will check when the Hodge isometry induced by the autoequivalence $\funct{T}_U\circ\funct{T}_{\sho_S}\circ(-\otimes\sho_S(H))$
is an involution. 

One writes $v(U)=(a,bH,c)$ for integers $a,b,c$ with $a>0$ and $2b^2d=2ac-2$. Then the Hodge isometry $\tau_U$ induced by that functor is given by the matrix 
$$Q_U=\begin{pmatrix}
	1-ac& abH& -a^2\\-bcH& b^2H*H+1& -abH\\ -c^2& bcH& 1-ac
\end{pmatrix}\begin{pmatrix}
	-d&-H&-1\\H&1&0\\ -1&0&0
\end{pmatrix}$$
where $a,d>0$ and $b$ are all independent integers.

To make the computation easier, we will use the matrix
$$
Q_{a,b}(x)=\begin{pmatrix}
	-b^2x&ab\sqrt{2x}&-a^2\\
	-\frac{b(b^2x+1)}{a}\sqrt{2x}&2b^2x+1&-ab\sqrt{2x}\\ 
	-\frac{(b^2x+1)^2}{a^2}&\frac{b(b^2x+1)}{a}\sqrt{2x}&-b^2x
\end{pmatrix}\begin{pmatrix}
	-x&-\sqrt{2x}&-1\\\sqrt{2x}&1&0\\ -1&0&0
\end{pmatrix}
$$
so that $Q_U^n=I_3$ implies that $Q_{a,b}(d)^n=I_3$. A long computation reveals that $Q_{a,b}(x)$ is 
$$\frac{1}{a^2}\begin{pmatrix}
	a^2(bx+a)^2&a^2b(bx+a)\sqrt{2x}&a^2b^2x\\
	a(b^2x+ab+1)(bx+a)\sqrt{2x}&2abx(b^2x+ab+1)+a^2&ab(b^2x+1)\sqrt{2x}\\
	x(b^2x+ab+1)^2&(b^2x+ab+1)(b^2x+1)\sqrt{2x}&(b^2x+1)^2
\end{pmatrix}$$
and the characteristic polynomial (with variable $\lambda$) of it is exactly
$$(\lambda-1)(\lambda^2-\frac{(abx+b^2x+a^2+1)^2-2a^2}{a^2}\lambda+1)$$
Therefore, the matrix $Q_{a,b}(x)$ has a simple eigenvalue $1$ and two complicated ones
$$\frac{(abx+b^2x+a^2+1)^2-2a^2\pm\sqrt{(abx+b^2x+a^2+1)^2((abx+b^2x+a^2+1)^2-4a^2)}}{2a^2}$$ 
which might not be real. To simplify notations, we will set $t:=abx+b^2x+a^2+1$.

Suppose that $t^2> 4a^2$, then all the eigenvalues are real and distinct. One sees that
$$\frac{t^2-2a^2+\sqrt{t^2(t^2-4a^2)}}{2a^2}>\frac{t^2-2a^2-\sqrt{t^2(t^2-4a^2)}}{2a^2}>0$$
so it is impossible to find an integer $m>0$ such that
$$(\frac{t^2-2a^2+\sqrt{t^2(t^2-4a^2)}}{2a^2})^m=(\frac{t^2-2a^2-\sqrt{t^2(t^2-4a^2)}}{2a^2})^m$$

Then we come to the case $0\leq t^2\leq 4a^2$. In this case, one can write
$$\frac{t^2-2a^2\pm\sqrt{t^2(t^2-4a^2)}}{2a^2}=e^{\pm \theta\sqrt{-1}}$$
where $\theta\in[0,2\pi)$. Moreover, one must have $\theta\in\bq\cdot\pi$, provided that 
$$e^{-m\theta\sqrt{-1}}=e^{m\theta\sqrt{-1}}=1$$ 
for some integer $m>0$. On the other hand, the following statement is known.

\begin{proposition}
One has $\bq\cap\{\cos(\theta)\,|\,\theta\in\bq\cdot\pi\}=\{0,\pm 1,\pm\frac{1}{2}\}$.
\end{proposition}

It follows immediately that $t^2\in\{0,a^2,2a^2,3a^2,4a^2\}$. Since $a$ is an integer and we want $t$ to be an integer, only $t^2=4a^2$ or $t^2=0$ can happen. 

Suppose that $t^2=4a^2$, then the eigenvalues of $Q_{a,b}(x)$ are all $1$. So it is either the identity matrix $I_3$ itself or its power can never be $I_3$. Apparently, it is not $I_3$ as 
$$a^2b^2x=0\quad\textup{and}\quad x(b^2x+ab+1)^2=0$$
cannot hold at the same time, for integers $a,x>0$ and $b$.

Suppose that $t^2=0$, then the eigenvalues of $Q_{a,b}(x)$ are $1$ and two $-1$. Also, we have
$$(a+b)bx+a^2+1=0$$
which has no solutions or a unique solution $x_0$ for $x$. One notices that
$$Q_{a,b}(x_0)=\begin{pmatrix}
	\frac{(ab-1)^2}{(a+b)^2}&\frac{b(ab-1)}{a+b}\sqrt{-2\frac{a^2+1}{b(a+b)}}&-\frac{b(a^2+1)}{a+b}\\
	\frac{(b^2+1)(ab-1)}{(a+b)^2}\sqrt{-2\frac{a^2+1}{b(a+b)}}&-\frac{2(a^2+1)(b^2+1)}{(a+b)^2}+1&-\frac{b(ab-1)}{a+b}\sqrt{-2\frac{a^2+1}{b(a+b)}}\\
	-\frac{(a^2+1)(b^2+1)^2}{b(a+b)^3}&-\frac{(b^2+1)(ab-1)}{(a+b)^2}\sqrt{-2\frac{a^2+1}{b(a+b)}}&\frac{(ab-1)^2}{(a+b)^2}
\end{pmatrix}$$

\begin{proposition}
The Jordan normal form of $Q_{a,b}(x_0)$ is a diagonal matrix.
	
\begin{proof}
It suffices to determine the space of eigenvectors for $-1$. One notices that 
$$Q_{a,b}(x_0)+I_3=\begin{pmatrix}
			\frac{(a^2+1)(b^2+1)}{(a+b)^2}&\frac{b(ab-1)}{a+b}\sqrt{-2\frac{a^2+1}{b(a+b)}}&-\frac{b(a^2+1)}{a+b}\\
			\frac{(b^2+1)(ab-1)}{(a+b)^2}\sqrt{-2\frac{a^2+1}{b(a+b)}}&-\frac{2(ab-1)^2}{(a+b)^2}&-\frac{b(ab-1)}{a+b}\sqrt{-2\frac{a^2+1}{b(a+b)}}\\
			-\frac{(a^2+1)(b^2+1)^2}{b(a+b)^3}&-\frac{(b^2+1)(ab-1)}{(a+b)^2}\sqrt{-2\frac{a^2+1}{b(a+b)}}&\frac{(a^2+1)(b^2+1)}{(a+b)^2}
		\end{pmatrix}$$
		has rank $1$, so we are done.
\end{proof}
\end{proposition}

\begin{remark}
One notices that
$$(0,b(a+b)\sqrt{-2\frac{a^2+1}{b(a+b)}},-2(ab-1))\quad\textup{and}\quad(2(ab-1),(b^2+1)\sqrt{-2\frac{a^2+1}{b(a+b)}},0)$$
are two linearly independent eigenvectors for $-1$.
\end{remark}

In this case, any integer solution $(a,b,d)$ for the equation
$$(a+b)bd+a^2+1=0$$
such that $a>0,0>b>-a$ and $d>0$, ensures that $Q_{a,b}(d)^2=I_3$. It also ensures that 
$$c=\frac{b^2d+1}{a}=-a-bd\in\bz$$
i.e.\ the vector $(a,bH,c)$ is really a spherical class.

The equation can be rewritten into the form
$$d=2+\frac{(a+b)^2+(-b)^2+1}{(a+b)(-b)}$$
where $-b>0,a+b>0$. Then we can conclude by the following elementary fact

\begin{proposition}
The only possible integer $t$ such that $x^2+y^2+1=txy$ for two positive integers $x\geq y$ is $t=3$. In this case, the solutions $(x,y)$ are $(F_n,F_{n-2})$ for even $n\geq0$ where the sequence $\{F_n\}$ is defined by $F_n=1$ for $n\leq 1$ and $F_{n+2}=F_{n+1}+F_n$ for $n\geq0$.
	
\begin{proof}
It is clear that $t\geq 3$; otherwise we will have $(x-y)^2+1\leq 0$ which is absurd. Also, one can see that $(x,y,t)=(1,1,3)$ is the only solution for $x^2+y^2+1=txy$ with $x=y$.
		
Now let us fix a solution $(x_0,y_0,t_0)$, then the equation $x^2-t_0y_0x+y_0^2+1=0$ has two roots $x_0$ and $t_0y_0-x_0$. The latter root can be rewritten as
$$t_0y_0-x_0=\frac{x_0^2+y_0^2+1}{x_0}-x_0=\frac{y_0^2+1}{x_0}$$
which is less than $y_0$ unless $(x_0,y_0,t_0)=(1,1,3)$. In any case, $(y_0,t_0y_0-x_0,t_0)$ is also a solution of our equation. It means that, for any solution $(x_0,y_0,t_0)$, the sequence
$$x_{n+1}=y_n\quad\textup{and}\quad y_{n+1}=t_0y_n-x_n$$
will go to $(x_m,y_m,t_0)=(1,1,3)$ in finite steps. In particular, the only possible integer $t$ to make the equation hold is $3$, and all the possible solutions $(x,y,t)$ are generated by
$$(x_0,y_0)=(1,1)\quad\textup{and}\quad (x_{n+1},y_{n+1})=(3x_n-y_n,x_n)\,\,\textup{for }n\geq0$$
which are precisely the ones we described.
\end{proof}
\end{proposition}

So $d=5$ and the set of solutions $(a,b)$ for the equation consists of $(F_{n+2}+F_n,-F_n)$ and $(F_n+F_{n-2},-F_n)$ for even $n\geq0$. This allows us to conclude that

\begin{theorem}
Consider a polarized K3 surface $(S,H)$, a spherical object $U$ of $\cate{D}^b(S)$ with Mukai vector $(a,bH,c)$ for some $a>0$, and the Hodge isometry $\tau_U$ on $\tilde{H}(S,\bz)$ induced by the autoequivalence $\funct{T}_U\circ\funct{T}_{\sho_S}\circ(-\otimes\sho_S(H))$, then $\tau_U^2=\id$ if and only if $d=5$ and $(a,b,c)$ is in 
$$\{(F_n+F_{n-2},-F_n,F_{n+2}+F_{n}),(F_{n+2}+F_n,-F_n,F_n+F_{n-2})\,|\,n\geq 0\textup{ even}\}$$
where $\{F_n\}$ is a sequence defined by $F_n=1$ for $n\leq 1$ and $F_{n+2}=F_{n+1}+F_n$ for $n\geq0$.
\end{theorem}

\begin{example}
Let $(S,H)$ be a Brill--Noether general K3 surface of degree 10, then $-\tau_U$ is
$$(r,\Delta,s)\mapsto(-9r-5s-3\Delta\cdot H,(6r+3s+2\Delta\cdot H)H-\Delta,-20r-9s-6\Delta\cdot H)$$
for $v(U)=(2,-H,3)$. It fixes the stability condition $\sigma_{\omega_2,\beta_2}$ and can be realized as an involution on $\cate{D}^b(S)$ isomorphic to $\Pi_2$. Similarly, $-\tau_U$ corresponds to $\Pi_1$ for $v(U)=(3,-H,2)$.
\end{example}

In fact, each such $-\tau_U$ fixes the stability function $Z_{\omega,\beta}$ of a $t$-structure $\sigma_{\omega,\beta}$ (which is a priori not necessarily a stability condition). One can check that $\sigma_{\omega,\beta}$ is a stability condition under certain extra numerical conditions (for example, $S$ has Picard rank $1$). Once $\sigma_{\omega,\beta}$ is indeed a stability condition, the Hodge isometry $-\tau_U$ can be lifted to an involution on $\cate{D}^b(S)$ according to Proposition~\ref{Appendix-Huybrechts-conway}. These involutions are conjugate to each other by \cite[Theorem 5.5]{FL23}, so they all induce a group action on $\cate{D}^b(S)$ and the equivariant categories are equivalent.

\begin{example}
The solution $(F_{n+2}+F_n,-F_n,F_n+F_{n-2})$ corresponds to $\tau_U$ with representation
$$Q_U=\begin{pmatrix}
	-5\frac{F_n^3+F_n}{F_{n+2}}-1&F_n(4F_n-F_{n+2})H&5F_n^2\\
	A_n(-5\frac{F_n^3+F_n}{F_{n+2}})H&(A_nF_n(4F_n-F_{n+2})-\frac{1}{5})H*H+1&A_n(5F_n^2) H\\
	B_n(-5\frac{F_n^3+F_n}{F_{n+2}})&B_nF_n(4F_n-F_{n+2}) H&B_n(5F^2_n)-1
\end{pmatrix}$$
where the auxiliary sequences $A_n$ and $B_n$ are defined by
$$A_n=-\frac{4F_n-F_{n+2}}{5F_n}\quad\textup{and}\quad B_n=\frac{F_n^2+1}{F_nF_{n+2}}$$
One computes that $Z_{\omega,\beta}(r,\Delta,s)=Z_{\omega,\beta}(-\tau_U(r,\Delta,s))$ holds for any $(r,\Delta,s)$ only when
$$\omega_n=\frac{1}{5F_n}H\quad\textup{and}\quad \beta_n=A_nH$$
and, suppose in addition that $(S,H)$ is a generic polarized K3 surface of degree $10$, one checks by definition that $\sigma_{\omega_n,\beta_n}$ is a stability condition.\footnote{The fixed function $Z_{\omega,\beta}$ associated with the solution $(F_n+F_{n-2},-F_n,F_{n+2}+F_n)$ is left to readers. In the generic case, one can check that the pair $(\omega,\beta)$ really gives a stability condition $\sigma_{\omega,\beta}$ as well.} Moreover, one has $-\tau_U(2,-H,2)=(2,-H,2)$.
\end{example}

}

\end{document}